\documentclass[a4paper,reqno,11pt]{amsart}

\usepackage{txfonts} 
\usepackage{amssymb, amsfonts, amsmath, amsthm}
\usepackage{mathrsfs} 
\usepackage{stmaryrd} 
\usepackage{esint} 
\usepackage{extarrows} 
\usepackage{chemarrow} 
\usepackage[normalem]{ulem}   
\usepackage{soul}             
\usepackage{xcolor} 

\usepackage[numbers,sort&compress]{natbib}

\usepackage{booktabs, tabularx}
\usepackage{tikz} 
\usetikzlibrary{calc} 

\usepackage[left=1in, right=1in, top=1in, bottom=1in]{geometry}

\newtheorem{Theorem}{Theorem}[section]
\newtheorem{Lemma}{Lemma}[section]
\newtheorem{Proposition}{Proposition}[section]
\newtheorem{Remark}{Remark}[section]
\newtheorem{Corollary}{Corollary}[section]
\newtheorem{Definition}{Definition}[section]

\numberwithin{equation}{section}

\usepackage{hyperref}
\hypersetup{
    bookmarksnumbered=true,
    colorlinks=false, 
    pdfborder={0 0 1},
    linkcolor=blue,       
    citecolor=green,      
    urlcolor=magenta,     
}

\usepackage{bookmark} 

\DeclareMathOperator{\diver}{div}
\DeclareMathOperator{\cur}{curl}

\definecolor{deepblue}{RGB}{25,25,112}      
\definecolor{deepgreen}{RGB}{0,100,0}       
\definecolor{deepred}{RGB}{139,0,0}         
\definecolor{deepmaroon}{RGB}{128,0,0}      

\newcommand{\na}{\nabla}

\newcommand{\al}{\alpha}
\newcommand{\be}{\beta}

\newcommand{\om}{\omega}

\newcommand{\de}{\delta}
\newcommand{\si}{\sigma}
\newcommand{\pa}{\partial}

\newcommand{\va}{\varepsilon}
\newcommand{\De}{\Delta}

\def\r3{\mathbb{R}^3}

\makeatletter
\@namedef{subjclassname@2020}{\textup{2020} Mathematics Subject Classification}
\makeatother

\allowdisplaybreaks

\def\XXint#1#2#3{{\setbox0=\hbox{$#1{#2#3}{\int}$}
    \vcenter{\hbox{$#2#3$}}\kern-.6\wd0}}

\begin{document}
\bibliographystyle{plain}

\title[I\MakeLowercase{ncompressible viscoelastic equations}]
{G\MakeLowercase{lobal well-posedness of incompressible viscoelastic flows with large velocity oscillations}}

\author{Yong Wang}
\address{Yong Wang \newline School of Mathematical Sciences, Key Laboratory of Mathematical Foundation and Computation in Nonlinear Sciences of Guangdong Higher Education Institutes, and Guangdong Key Laboratory of Applied Mathematics and Artificial Intelligence, South China Normal University, Guangzhou, Guangdong 510631, China}
\email{wangyongxmu@163.com}

\author{Zhouping Xin}
\address{Zhouping Xin \newline The Institute of Mathematical Sciences, The Chinese University of Hong Kong, Shatin, N.T., Hong Kong, China}
\email{zpxin@ims.cuhk.edu.hk}

\author{Jianquan Yang$^\ast$}
\address{Jianquan Yang \newline School of Mathematical Sciences, Key Laboratory of Mathematical Foundation and Computation in Nonlinear Sciences of Guangdong Higher Education Institutes, and Guangdong Key Laboratory of Applied Mathematics and Artificial Intelligence, South China Normal University, Guangzhou, Guangdong 510631, China}
\email{yangjianquanscnu@163.com}

\author{Yanting Yuan}
\address{Yanting Yuan \newline School of Mathematical Sciences, Key Laboratory of Mathematical Foundation and Computation in Nonlinear Sciences of Guangdong Higher Education Institutes, and Guangdong Key Laboratory of Applied Mathematics and Artificial Intelligence, South China Normal University, Guangzhou, Guangdong 510631, China}
\email{ytyuan2404@163.com}

\thanks{$^\ast$Corresponding author: Jianquan Yang (yangjianquanscnu@163.com).}

\begin{abstract}
We investigate the global well-posedness of classical solutions to the Cauchy problem of the incompressible viscoelastic equations in both two and three dimensions.
We first establish some new blow-up criteria, which reveal that the deformation gradient of the incompressible viscoelastic fluid will
become unbounded at some points in space
if the classical solution develops singularities in finite time.
It is essentially shown that for the two-dimensional (2D) Cauchy problem,
the classical solution exists globally if the deformation gradient is bounded from above; the same conclusion holds for the 3D case, provided the velocity or pressure remains bounded in certain Serrin-type norms.
These results extend the well-known Serrin-type blow-up criteria in terms of the velocity or pressure for the 3D incompressible Navier-Stokes equations to the 3D viscoelastic case and generalize the results for periodic problems in 
[Feng-Zhu-Zi, J. Funct. Anal., 2017].

Motivated by these blow-up criteria, we prove the global-in-time existence of classical solutions in $H^m$ ($m\geq2$) class under the assumptions that the initial velocity is small in the natural energy space $L^2$ and the initial deformation gradient is close to the identity matrix in $L^2\cap{\dot{B}_{p,1}^{\frac{n}{p}}}$, where $p\in(2,\infty)$ for $n=2$ and $p\in(3,6]$ for $n=3$.
In particular, the oscillations of the initial velocity can be arbitrarily large.
Our work improves the previous results
[Chen-Zhang, Comm. Partial Differential Equations, 2006;
Lei-Liu-Zhou, Arch. Ration. Mech. Anal., 2008;
Lei-Zhou, SIAM J. Math. Anal., 2005;
Lin-Liu-Zhang, Comm. Pure Appl. Math., 2005]
on global classical solutions for small $H^2$ perturbations near a constant equilibrium.
It should be remarked that the effective viscous fluxes play a key role in analyzing the initial layer and closing the required energy estimates.
\end{abstract}

\maketitle

\setcounter{tocdepth}{3}


\pagenumbering{arabic}

\section{Introduction}

In this paper, we investigate the following incompressible viscoelastic system \cite{Liu-Walkington2001,Lin-Liu-Zhang2005,Lei-Liu-Zhou2008,Giga-Kirshtein-Liu2018}:
    \begin{equation}
		\label{1.1}
		\begin{cases}
			\ u_t+u\cdot\na u-\mu\Delta u+\nabla P
            =\diver(\mathbb{F}\mathbb{F}^{\top}),\\
			\ \mathbb{F}_t+u\cdot\nabla\mathbb{F}=\nabla u\mathbb{F},\\
            \ \diver u=0
		\end{cases}
	\end{equation}
for $(x,t)\in\mathbb{R}^n\times\mathbb{R}^+$ with $n=2, 3$.
The unknown variables $u=u(x,t)\in\mathbb{R}^n$, $P=P(x,t)$, and $\mathbb{F}=\mathbb{F}(x,t)\in\mathbb{M}^{n\times n}$ (the set of $n\times n$ matrices with positive determinants) denote the velocity field, the pressure, and the deformation gradient, respectively.
The constant $\mu>0$ is the fluid viscosity.
The superscript $\top$ means the transpose.
We will consider the Cauchy problem of the system \eqref{1.1} with the initial and far-field conditions
    \begin{equation}
		\label{1.2}
        \begin{cases}
			\ (u,\mathbb{F})(x,t)\big|_{t=0}
            =(u_0,\mathbb{F}_0)(x),& x\in\mathbb{R}^n,\\
			\ (u,\mathbb{F})(x,t)
            \to (0,\mathbb{I}), & \text{as}\ |x|\to\infty,\ t>0,
		\end{cases}
	\end{equation}
where $\mathbb{I}$ is the identity matrix.

There is a vast literature on the existence, uniqueness and large-time behavior of solutions to the system \eqref{1.1}.
For the 2D Cauchy problem, Lin, Liu, and Zhang \cite{Lin-Liu-Zhang2005} studied an auxiliary vector field instead of the transport variable $\mathbb{F}$ and proved the global existence of classical solutions under small initial data; see also \cite{Lei-Zhou2005} for a similar result obtained via the incompressible limit method.
For the 3D Cauchy case, motivated by the ideas in \cite{Sideris-Thomases2005} for the incompressible elastodynamics, Chen and Zhang \cite{Chen-Zhang2006} established the global well-posedness near a constant equilibrium state by introducing a curl free quantity $\mathbb{F}^{-1}-\mathbb{I}$.
Observing that $\diver\mathbb{F}$ provides some weak dissipation and $\cur\mathbb{F}$ is a high-order nonlinear term, Lei, Liu, and Zhou \cite{Lei-Liu-Zhou2008} showed the existence of global classical solutions near equilibrium
for Cauchy and periodic problems in both $2D$ and $3D$.
For the convenience of the reader, we summarize the above results as follows.
Further details can be found in \cite{Lin-Liu-Zhang2005,Lei-Zhou2005,Chen-Zhang2006,Lei-Liu-Zhou2008}.
    \begin{Proposition}
	\label{prop1.0}
	Let $n=2$ or $3$, and $m\geq2$ be a positive integer.
Assume that the initial data $(u_0,\mathbb{F}_0-\mathbb{I})\in H^m(\mathbb{R}^n)$ satisfying
    \begin{equation}
		\label{1.4-0321}
		\begin{cases}
            \ \diver u_0=0,\\
            \ \det\mathbb{F}_0=1,\\
            \ \diver \mathbb{F}_0^{\top}=0,\\
			\ \mathbb{F}_{0}^{lk}\na_l\mathbb{F}_{0}^{ij}
            =\mathbb{F}_{0}^{lj}\na_l\mathbb{F}_{0}^{ik},\quad
            1\leq i,j,k\leq n.
		\end{cases}
	\end{equation}
Then there exists a small constant $\va_0>0$ depending only on $n$ and $\mu$, such that if
    \begin{equation}
		\label{prop1.1-2-0321}
		\|u_{0}\|_{H^2(\mathbb{R}^n)}^2+\|\mathbb F_{0}-\mathbb I\|_{H^2(\mathbb{R}^n)}^2\leq \va_0,
	\end{equation}
then the Cauchy problem \eqref{1.1}-\eqref{1.2} admits a unique global classical solution $(u,\mathbb{F})$ satisfying
    \begin{equation*}
		\|u\|_{H^2(\mathbb{R}^n)}^2+\|\mathbb F-\mathbb I\|_{H^2(\mathbb{R}^n)}^2
        +\int_0^T \|\na u\|_{H^2(\mathbb{R}^n)}^2\,dt
        \leq C_0\Big(\|u_{0}\|_{H^2(\mathbb{R}^n)}^2+\|\mathbb F_{0}-\mathbb I\|_{H^2(\mathbb{R}^n)}^2\Big)
	\end{equation*}
for any $T>0$.
The above constant $C_0>0$ depends only on $n$ and $\mu$.
\end{Proposition}
Later, by means of the rotation-strain decomposition of the deformation gradient, Lei \cite{Lei2010,Lei2014} proved that the smallness of the initial strain tensor suffices to guarantee the global regularity of the system \eqref{1.1}.
In particular, the initial deformation gradient is allowed to be large.
We also refer to \cite{Lin-Zhang2008,He-Xu2010} for studies on the initial-boundary value problem, and to \cite{Qian2010,Zhang-Fang2012} for results in the critical Besov framework.
Indeed, the aforementioned works make sufficient use of the structure assumptions $\eqref{1.4-0321}_2$-$\eqref{1.4-0321}_4$, which play an important role in deriving the dissipation of the deformation gradient $\mathbb{F}$.
It is worth noting that Zhu \cite{Zhu2018} established the global existence of small solutions to the 3D incompressible Oldroyd-B model with only velocity dissipation by developing a time-weighted energy framework; see also \cite{Chen-Zhu2023} for the 2D case.
Furthermore, by regarding the ``nonlinear term'' $\diver(\mathbb{F}\mathbb{F}^{\top})$ as a ``linear term'', the global small classical solutions to the 3D incompressible viscoelastic equations without structure assumptions was also addressed in \cite{Zhu2018}.
As for the incompressible elastodynamic system (i.e., $\mu=0$ in \eqref{1.1}),
the global existence of classical solutions with small initial displacements for the 3D Cauchy problem was established by Sideris and Thomases \cite{Sideris-Thomases2005,Sideris-Thomases2007}.
As pointed out in \cite{Lei-Sideris-Zhou2015}, the existence problem of the 2D incompressible nonlinear elastodynamics is more complicated, since quadratic nonlinearities---even under the null condition---have at most critical time decay.
Later, by exploiting the inherent strong null structure of nonlinearities in the Lagrangian coordinate, Lei \cite{Lei2016} obtained the global-in-time classical solutions for small initial data in the 2D case; see also \cite{Wang2017} for a different approach working in the frequency space.
For more studies on the system \eqref{1.1} and closely related models, we refer the reader to \cite{Lions-Masmoudi2000,Qian-Zhang2010,Lin2012,Xu-Zhang-Zhang2013,
Masmoudi2013,Cai-Lei-Lin-Masmoudi2019,D-M-S2020,Huang-Yao-You2025} and the references cited therein.

To the best of our knowledge, the global well-posedness of the incompressible viscoelastic equations \eqref{1.1} with large initial data remains a major open problem listed in \cite{Hu-Lin-Liu2018}.
%
%
This motivates us to investigate the singularity structure and blow-up mechanism of solutions to the system \eqref{1.1}.
In fact, there are several results on the blow-up criteria of classical (or smooth) solutions to \eqref{1.1}; see \cite{Lin-Liu-Zhang2005,Chen-Zhang2006,Lei-Liu-Zhou2008,Yuan-Li2013,
Yuan2013,Feng-Zhu-Zi2017} for instance.
More precisely, let $0<T^*<\infty$ be the maximum time of existence of classical solutions. Then the blow-up criteria can be summarized as follows.
    \begin{itemize}
        \item Lin-Liu-Zhang (\cite{Lin-Liu-Zhang2005}, 2D Cauchy problem), Lei-Liu-Zhou (\cite{Lei-Liu-Zhou2008}, 2D and 3D Cauchy problems)
        \begin{equation*}
		  \int_0^{T^*}\|\na u\|_{H^2(\mathbb{R}^n)}^2\,dt=\infty;
	    \end{equation*}
        \item Chen-Zhang (\cite{Chen-Zhang2006}, 3D Cauchy problem), Yuan (\cite{Yuan2013}, 2D and 3D Cauchy problems)
        \begin{equation*}
		  \int_0^{T^*}\|\na u\|_{L^\infty(\mathbb{R}^n)}\,dt=\infty;
	    \end{equation*}
        \item Yuan-Li (\cite{Yuan-Li2013}, 2D and 3D Cauchy problems)
        \begin{equation*}
		  \begin{cases}
			\ \displaystyle\int_0^{T^*}
            \|\na \mathbb{F}\|_{BMO(\mathbb{R}^2)}\,dt=\infty,&n=2,\\
            \ \displaystyle\int_0^{T^*}\Big(\|\na u\|_{BMO(\mathbb{R}^3)}+\|\na \mathbb{F}\|_{BMO(\mathbb{R}^3)}^2\Big)\,dt=\infty,&n=3.
		\end{cases}
	    \end{equation*}
    \end{itemize}
The above regularity criteria depend crucially on some norms of $\na u$ or $\na \mathbb{F}$.
Indeed, if the initial deformation gradient $\mathbb{F}_0$ is bounded, the boundedness of $\|\na u\|_{L^1(0,T;L^\infty(\mathbb{R}^n))}$ together with the transport equation $\eqref{1.1}_2$ yields the uniform upper bound of $\mathbb{F}$; see
\cite{Chen-Zhang2006,Yuan2013}.
Inspired by the works \cite{Huang-Li-Xin2011,Sun-Wang-Zhang2011} on the compressible Navier-Stokes equations, Feng, Zhu, and Zi \cite{Feng-Zhu-Zi2017} proved the following Serrin-type (see \cite{Serrin1962}) criterion of strong solutions to system \eqref{1.1}
in the periodic domain $\mathbb{T}^n$ ($n=2,3$).
    \begin{itemize}
        \item Feng-Zhu-Zi (\cite{Feng-Zhu-Zi2017}, 2D and 3D periodic problems)
        \begin{equation}
		    \label{1.5}
            \begin{cases}
			\ \lim\limits_{T\to T^*}
            \|\mathbb{F}\|_{L^\infty(0,T;L^\infty(\mathbb{T}^2))}=\infty,&n=2,\\
            \ \lim\limits_{T\to T^*}\left(\|\mathbb{F}\|_{L^\infty(0,T;L^\infty(\mathbb{T}^3))}
            +\|u\|_{L^{s}(0,T;L^{r}(\mathbb{T}^3))}\right)=\infty,&n=3,
		\end{cases}
	    \end{equation}
        where $r$ and $s$ satisfy $\frac{2}{s}+\frac{3}{r}\leq1$ for $3< r\leq \infty$.
    \end{itemize}
However, as noted by the authors in \cite{Feng-Zhu-Zi2017}, their methods are not applicable to the Cauchy problem or the initial-boundary value problem.
Consequently, to establish Serrin-type blow-up criteria for the Cauchy problem of \eqref{1.1}, some new ideas are needed.
We also refer to \cite{Chemin-Masmoudi2001,Kupferman-Mangoubi-Titi2008,Lei-Masmoudi-Zhou2010} for the blow-up criteria to the viscoelastic system of Oldroyd-type, and to \cite{Hu-Hynd2013} for a Beale-Kato-Majda-type criterion for the inviscid case (i.e., $\mu=0$ in \eqref{1.1}).

In this paper, we are concerned with the global well-posedness of classical solutions to the two- and three-dimensional incompressible viscoelastic system \eqref{1.1}.
A deep understanding of the mechanism behind the formation of singularities is therefore needed.
Here, we recall some important results on the classical incompressible Navier-Stokes equations
    \begin{equation}
		\label{1.1-0338}
		\begin{cases}
			\ u_t+u\cdot\na u-\mu\Delta u+\nabla P=0,\\
            \ \diver u=0,
            \quad (x,t)\in\mathbb{R}^n\times\mathbb{R}^+,\ n=2,3,
		\end{cases}
	\end{equation}
which is the subsystem of \eqref{1.1} obtained by ignoring elastic effects.
It is well-known that in two dimensions, the system \eqref{1.1-0338} is globally well-posed for smooth initial data (see \cite{Ladyzhenskaya1969}).
As for the three-dimensional case, it is also well-known that the system \eqref{1.1-0338} has at least a Leray-Hopf weak solution provided the initial kinetic energy is finite (see \cite{Leray1934,Hopf1951}).
Over the years, many sufficient conditions have been obtained under which Leray-Hopf solutions are smooth.
In particular, Serrin et al. showed that a weak solution $(u,P)$ to the 3D incompressible Navier-Stokes system is regular if the velocity $u$ or the pressure $P$ satisfies one of the following conditions:
    \begin{itemize}
        \item Serrin \cite{Serrin1962}, Sohr \cite{Sohr1983}
        \begin{equation*}
        u\in L^s(0,T;L^r(\mathbb{R}^3))
        \quad\text{with}\quad
        \frac{2}{s}+\frac{3}{r}\leq1,\
        3< r\leq \infty
	    \end{equation*}
        (see also \cite{E-S-S2003} for the end point result $(r,s)=(3,\infty)$ by Escauriaza, Seregin, and \v{S}ver\'{a}k);
        \item Beir\~{a}o da Veiga \cite{B1995}, Galdi \cite{Galdi2000}
        \begin{equation*}
        \na u\in L^s(0,T;L^r(\mathbb{R}^3))
        \quad\text{with}\quad
        \frac{2}{s}+\frac{3}{r}\leq2,\
        \frac{3}{2}< r\leq \infty;
	    \end{equation*}
        \item Chae-Lee \cite{Chae-Lee2001}, Berselli-Galdi \cite{Berselli-Galdi2002}
        \begin{equation*}
        P\in L^s(0,T;L^r(\mathbb{R}^3))
        \quad\text{with}\quad
        \frac{2}{s}+\frac{3}{r}\leq2,\
        \frac{3}{2}< r\leq \infty;
	    \end{equation*}
        \item Zhou \cite{Zhou2006,Zhou2006-2}
        \begin{equation*}
        \na P\in L^s(0,T;L^r(\mathbb{R}^3))
        \quad\text{with}\quad
        \frac{2}{s}+\frac{3}{r}\leq3,\
        1< r\leq \infty.
	    \end{equation*}
    \end{itemize}
Indeed, these Serrin-type regularity criteria are of great significance for studying well-posedness and singularity of system \eqref{1.1-0338}; see \cite{Fujita-Kato1964,Giga1986,Kato1984} for example.
For further research on the incompressible Navier-Stokes system, the reader can refer to \cite{Caffarelli-Kohn-Nirenberg1982,Ponce-Racke-Sideris-Titi1994,
Buckmaster-Vicol2019,A-B-C2022,Coiculescu-Palasek2026,C-M-Z2026,A-B-2026} and the references therein.
Turning to the incompressible viscoelastic system \eqref{1.1}, one finds that elastic effects introduce significant challenges.
Since $\eqref{1.1}_2$ is a hyperbolic equation with no damping or dissipative mechanism for the deformation gradient, the strong coupling between the velocity field and the deformation gradient---e.g., the influence of the nonlinear term $\diver(\mathbb{F}\mathbb{F}^{\top})$ on the velocity field---poses some essential difficulties.
To overcome these difficulties, we will make full use of the structure of the equations \eqref{1.1}.

Motivated by the works \cite{Serrin1962,Berselli-Galdi2002} for the incompressible Navier-Stokes equations, we first establish some new blow-up criteria of classical solutions to the system \eqref{1.1} in both two and three dimensions.
Then,
by combining the criteria we derived with some ideas from the studies on the compressible Navier-Stokes \cite{Huang-Li-Xin2012} and viscoelastic \cite{Hu2018} equations,
we prove the global-in-time existence of classical solutions under certain smallness assumptions on the initial data.
In particular, we greatly relax the condition \eqref{prop1.1-2-0321}.
Before stating our main results, let us introduce the notations and conventions that will be used throughout this paper.

\medskip

\textbf{Notations and conventions}.
Let $n=2,3$.
The symbol $\na=\pa_x=(\pa_1,\cdots,\pa_n)^\top$ represents the gradient operator and $\na_j=\pa_{x_j}\ (j=1,\cdots,n)$.
For a column vector $v=(v_1,,\cdots,\pa_n)^\top$, $(\na v)_{ij}=\na _jv_i$ and
    \begin{equation*}
        \cur v=\na v-(\na v)^\top.
	\end{equation*}
The material derivative for a function $f(x,t)$ is defined by
    \begin{equation}
        \label{n-m}
        \frac{D}{Dt}f=\dot{f}\triangleq f_t+u\cdot\nabla f,
	\end{equation}
where $u$ is the velocity field.
The Helmholtz projection operator is given by
    \begin{equation}
        \label{n-p}
        \mathcal{P}\triangleq I_d+\na(-\De)^{-1}\diver,
	\end{equation}
where $I_d$ is the identity mapping.
For any smooth divergence-free vector $\psi$ and scalar function $\varphi$, we have $\mathcal{P}\psi=\psi$ and $\mathcal{P}\na\varphi=0$.
And we denote
    \begin{equation}
        \label{n-q}
        \mathcal{Q}=I_d-\mathcal{P}=-\na(-\De)^{-1}\diver.
	\end{equation}
The Frobenius inner product of two matrices $\mathbb{A},\ \mathbb{B}\in \mathbb{R}^{n\times n}$ is denoted by $\mathbb{A}:\mathbb{B}\triangleq\sum_{i,j=1}^{n}\mathbb{A}^{ij}\mathbb{B}^{ij}$.
Particularly, $|\mathbb{A}|^2=\mathbb{A}:\mathbb{A}$.
We define the Fourier transform $\mathcal{F}$ and the inverse Fourier transform $\mathcal{F}^{-1}$ of a function $f(x)$ as follows:
    \begin{equation*}
		\mathcal{F}(f)(\xi)=\hat{f}(\xi)
        \triangleq
        \int_{\mathbb{R}^n}e^{-ix\cdot\xi}f(x)\,dx
        \quad\text{and}\quad
        \mathcal{F}^{-1}(f)(\xi)=\check{f}(\xi)
        \triangleq
        \int_{\mathbb{R}^n}e^{ix\cdot\xi}f(x)\,dx.
	\end{equation*}
For $1\leq p\leq\infty$, $k\in\mathbb{N}$ and $\eta>0$, the Lebesgue spaces, the homogeneous and inhomogeneous Sobolev spaces are denoted as follows:
    \begin{equation*}
        \begin{cases}
        \ L^p=L^p(\mathbb{R}^n),\quad
        D^{k,p}=\{u\in L^1_{\text{loc}}\mid\|\na^k u\|_{L^p}<\infty\},\\
        \ W^{k,p}=L^p\cap D^{k,p},\quad H^k=W^{k,2},\\
        \ \displaystyle\dot{H}^\eta=\dot{H}^\eta({\mathbb{R}^n})=
        \bigg\{f:\mathbb{R}^n\to\mathbb{R}\ \bigg|\ \|f\|_{\dot{H}^\eta}^2
        =\int_{\mathbb{R}^n}|\xi|^{2\eta}|\hat{f}(\xi)|^2\,d\xi<\infty\bigg\}.
		\end{cases}
	\end{equation*}
For $s\in \mathbb R$ and $1\leq p,r\leq\infty$, we denote by $\dot{B}_{p,r}^{s}=\dot B_{p,r}^s(\mathbb R^n)$ the homogeneous Besov space;
see Definition \ref{def2.1} in Section \ref{sec2} for the details.
The spaces involving time $L^p(0,T;X)$ ($1\leq p\leq\infty$) comprise all measurable functions $f: [0,T]\to X$ with $\|f\|_{L^{p}(0, T;X)}\triangleq \left(\int_0^T\|f(t)\|^{p}_{X} dt\right)^{\frac{1}{p}}<\infty$ for $1\leq p<\infty$ and $\|f\|_{L^{\infty}(0, T;X)}\triangleq \underset{0\leq t\leq T}{\operatorname{ess} \sup}\ \|f(t)\|_{X}<\infty$.
The usual Hardy and bounded mean oscillation ($BMO$) spaces are denoted by $\mathcal{H}^1(\mathbb{R}^n)$ and $BMO(\mathbb{R}^n)$, respectively (cf. \cite{Grafakos2014-M}).
Note that $BMO(\mathbb{R}^n)$ is the dual space of $\mathcal{H}^1(\mathbb{R}^n)$.
For simplicity, we write
$\|f\|_{X_1\cap X_2}\triangleq\|f\|_{X_1}+\|f\|_{X_2}$ for two Banach spaces $X_1$ and $X_2$.


\medskip



We shall study the system \eqref{1.1} with the constraints
    \begin{equation}
		\label{1.4-1224}
		\begin{cases}
            \ \det \mathbb{F}=1,\\
			\ \mathbb{F}^{lk}\na_l\mathbb{F}^{ij}
            =\mathbb{F}^{lj}\na_l\mathbb{F}^{ik},\quad
            1\leq i,j,k\leq n,
		\end{cases}
	\end{equation}
where and hereafter the Einstein summation convection is used.
It follows from Lemma \ref{le2.2} that the constraints in \eqref{1.4-1224} hold for all time $t>0$ if they are initially valid.
Moreover, by virtue of Lemma \ref{leH}, the constraints in \eqref{1.4-1224} imply for all $t\geq0$,
    \begin{equation}
        \label{1.4.1-1224}
        \diver\mathbb{F}^{\top}=0.
	\end{equation}


We first recall the local existence of classical solutions to the problem \eqref{1.1}-\eqref{1.2}.
The proof of the following proposition can be found in \cite{Lin-Liu-Zhang2005,Lei-Liu-Zhou2008}.

\begin{Proposition}(\textbf{Local Well-posedness}).
	\label{prop1.1}
	Let $n=2$ or $3$, and $m\geq2$ be a positive integer.
Assume that the initial data $(u_0,\mathbb{F}_0-\mathbb{I})\in H^m(\mathbb{R}^n)$ satisfying
    \begin{equation}
		\label{1.4}
		\begin{cases}
            \ \diver u_0=0,\\
            \ \det\mathbb{F}_0=1,\\
			\ \mathbb{F}_{0}^{lk}\na_l\mathbb{F}_{0}^{ij}
            =\mathbb{F}_{0}^{lj}\na_l\mathbb{F}_{0}^{ik},\quad
            1\leq i,j,k\leq n.
		\end{cases}
	\end{equation}
Then there exists a positive time $T$ depending only on $n$, $m$, $\mu$, $\|u_0\|_{H^2(\mathbb{R}^n)}$, and $\|\mathbb{F}_{0}-\mathbb{I}\|_{H^2(\mathbb{R}^n)}$, such that the Cauchy problem \eqref{1.1}-\eqref{1.2} admits a unique classical solution $(u,\mathbb{F})$ in $\mathbb{R}^n\times[0,T]$ satisfying
    \begin{equation}
		\label{prop1.1-2-0108}
		\begin{cases}
			\ \pa_t^\iota\na_x^\nu u
            \in L^\infty(0,T;H^{m-2\iota-|\nu|}(\mathbb{R}^n))\cap L^2(0,T;H^{m-2\iota-|\nu|+1}(\mathbb{R}^n)),\\
			\ \pa_t^\iota\na_x^\nu (\mathbb{F}-\mathbb{I})
            \in L^\infty(0,T;H^{m-2\iota-|\nu|}(\mathbb{R}^n))
		\end{cases}
	\end{equation}
for all positive integer $\iota$ and multiindex $\nu$ with $2\iota+|\nu|\leq m$.
\end{Proposition}

Then, the first main result on the blow-up criteria of classical solutions can be stated as follows.
\begin{Theorem}(\textbf{Blow-up Criteria}).
	\label{th1.1}
Let $n=2$ or $3$, $m\geq2$ be a positive integer.
Assume that the initial data $(u_0,\mathbb{F}_0-\mathbb{I})\in H^m(\mathbb{R}^n)$ satisfying \eqref{1.4}.
Let $(u,\mathbb{F})$ be the classical solution to the Cauchy problem \eqref{1.1}-\eqref{1.2} satisfying \eqref{prop1.1-2-0108}. If $0<T^{*}<\infty$ is the maximal time of existence, then for $n=2$,
    \begin{equation}
		\label{th1.1.2}
		\lim\limits_{T\to T^*}
            \|\mathbb{F}\|_{L^\infty(0,T;L^\infty(\mathbb{R}^2))}=\infty,
	\end{equation}
and for $n=3$,
    \begin{equation}
		\label{th1.1.2-2}
		\lim\limits_{T\to T^*}\left(\|\mathbb{F}\|_{L^\infty(0,T;L^\infty(\mathbb{R}^3))}
            +\|u\|_{L^{s_1}(0,T;L^{r_1}(\mathbb{R}^3))}\right)=\infty
        \quad\text{with}\quad
        \frac{2}{s_1}+\frac{3}{r_1}\leq1,\
        3< r_1\leq \infty,
	\end{equation}
and
    \begin{equation}
		\label{th1.1.2-1}
		\lim\limits_{T\to T^*}\left(\|\mathbb{F}\|_{L^\infty(0,T;L^\infty(\mathbb{R}^3))}
            +\|P\|_{L^{s_2}(0,T;L^{r_2}(\mathbb{R}^3))}\right)=\infty
        \quad\text{with}\quad
        \frac{2}{s_2}+\frac{3}{r_2}\leq2,\
        \frac{3}{2}< r_2\leq \infty.
	\end{equation}
\end{Theorem}

We make some remarks on Theorem \ref{th1.1} below.
\begin{Remark}\label{re1.1}
According to the criterion \eqref{th1.1.2}, the 2D incompressible viscoelastic flows will not develop singularities in finite time if the deformation gradient is bounded from above.
For the 3D case, the same conclusion holds provided that the velocity or pressure remains bounded in certain Serrin-type norms; see \eqref{th1.1.2-2} and \eqref{th1.1.2-1}.
Theorem \ref{th1.1} indicates that controlling the upper bound of the deformation gradient $\mathbb{F}$ is a key issue in establishing the global well-posedness of system \eqref{1.1}, which is consistent with the insight of Hu-Lin-Liu in \cite{Hu-Lin-Liu2018}.
\end{Remark}

\begin{Remark}
Theorem \ref{th1.1} also applies to periodic problems as considered in \cite{Feng-Zhu-Zi2017}.
Moreover, we claim that the criterion \eqref{th1.1.2-1} implies \eqref{th1.1.2-2} except for the case $r=\infty$.
In fact, the pressure could be expressed as (see also \eqref{3.4.4})
    \begin{equation*}
		P
        =(-\De)^{-1}\diver\diver(u\otimes u)
        -(-\De)^{-1}\diver\diver((\mathbb{F}-\mathbb{I})(\mathbb{F}^\top-\mathbb{I})),
	\end{equation*}
which, together with
the $L^p$-boundedness of the Riesz transform for $1<p<\infty$ (cf. \cite{Grafakos2014-C}),
gives
    \begin{equation*}
		\|P\|_{L^p(\mathbb{R}^3)}
        \leq C\Big(\|u\|_{L^{2p}(\mathbb{R}^3)}^2
        +\|\mathbb{F}-\mathbb{I}\|_{L^{2p}(\mathbb{R}^3)}^2\Big).
	\end{equation*}
It thus holds that
    \begin{align*}
        \|P\|_{L^{s_2}(0,T;L^{r_2}(\mathbb{R}^3))}
        \leq&\ C\Big(\|u\|_{L^{2s_2}(0,T;L^{2r_2}(\mathbb{R}^3))}
        +\|\mathbb{F}-\mathbb{I}\|_{L^{2s_2}(0,T;L^{2r_2}(\mathbb{R}^3))}\Big)\\
        \leq&\ C\Big(\|u\|_{L^{2s_2}(0,T;L^{2r_2}(\mathbb{R}^3))}
        +\|\mathbb{F}-\mathbb{I}\|_{L^{2s_2}(0,T;L^2(\mathbb{R}^3)\cap L^\infty(\mathbb{R}^3))}\Big)
	\end{align*}
for $r_2$ and $s_2$ satisfying $\frac{2}{s_2}+\frac{3}{r_2}\leq2$ and
$\frac{3}{2}< r_2<\infty$.
Besides, the elementary energy estimate of the incompressible viscoelastic system \eqref{1.1} gives an estimate on the $L^{\infty}(0,T;L^2(\mathbb{R}^3))$-norm of $\mathbb{F}-\mathbb{I}$
(see Lemma \ref{le3.1}).
Therefore, the criterion \eqref{th1.1.2-1} is stronger than \eqref{th1.1.2-2} when
$r$ is finite, and thus improves $\eqref{1.5}_2$.
In summary, Theorem \ref{th1.1} generalizes the results of \cite{Feng-Zhu-Zi2017} on periodic problems.
\end{Remark}

\begin{Remark}
The criteria \eqref{th1.1.2-2}-\eqref{th1.1.2-1} can be regarded as Serrin-type blow-up criteria (in terms of the velocity or pressure)
for the Cauchy problem of the 3D incompressible viscoelastic equations.
\end{Remark}

The second main result is concerned with the global-in-time existence of classical solutions.
\begin{Theorem}(\textbf{Global Well-posedness}).
	\label{th1.2}
Let $n=2$ or $3$, $m\geq2$ be a positive integer, and $M>0$ (not necessarily small).
Suppose $\eta$ and $p$ are two given constants satisfying
    \begin{align}
        &n=2:\eta\in(0,1],\ p\in(2,\infty);\label{th1.2.0}\\
        &n=3:\eta\in\Big(\frac{8}{15},1\Big],\ p\in(3,6].\label{th1.2.0-0413}
	\end{align}
Assume that the initial data $(u_0,\mathbb{F}_0-\mathbb{I})\in H^m(\mathbb{R}^n)$ satisfying \eqref{1.4} and
    \begin{equation}
        \label{th1.2.1}
        \|u_0\|_{\dot{H}^\eta(\mathbb{R}^n)}\leq M.
    \end{equation}
Then there exists a suitably small constant $\de>0$ depending only on $n$, $m$, $M$, $\eta$, $p$, and $\mu$, such that if
    \begin{equation}
        \label{th1.2.2}
        \mathcal{E}_0
        \triangleq\|u_{0}\|_{L^2(\mathbb{R}^n)}^2
        +\|\mathbb F_{0}-\mathbb I\|_{L^2(\mathbb{R}^n)\cap{\dot{B}_{p,1}^{\frac{n}{p}}}(\mathbb{R}^n)}^2
        \leq \de,
    \end{equation}
then the Cauchy problem \eqref{1.1}-\eqref{1.2} admits a unique global classical solution $(u,\mathbb{F})$ satisfying
    \begin{gather}
        \sup\limits_{0\leq t\leq T}\int_{\mathbb{R}^{n}}\Big(|u|^2+|\mathbb F-\mathbb I|^2\Big)\,dx
        +2\mu\int_0^T\int_{\mathbb{R}^{n}}|\nabla u|^2\,dxdt
        \leq \mathcal{E}_0,\nonumber\\
        \sup\limits_{0\leq t\leq T}\int_{\mathbb{R}^{n}}\Big(\si(t)|\nabla u|^2
        +\si(t)^n|\mathcal{P}\dot{u}|^2\Big)\,dx
        +\int_0^T\int_{\mathbb{R}^n}\Big(\si(t)|\mathcal{P}\dot{u}|^2
        +\si(t)^n|\nabla\mathcal{P}\dot{u}|^2\Big)\,dxdt
        \leq \mathcal{E}_0^\frac{1}{4},\label{th1.2-A1}\\
        \sup\limits_{0\leq t\leq \si(T)}\int_{\mathbb{R}^{n}}|u|^q\,dx
        +\int_0^{\si(T)}\int_{\mathbb{R}^n}\Big(|\na u|^2|u|^{q-2}
        +\big|\na|u|^{\frac{q}{2}}\big|^2\Big)\,dxdt
        \leq \mathcal{E}_0^\al,\nonumber\\
        \sup\limits_{0\leq t\leq T}
        \|\mathbb{F}-\mathbb I\|_{\dot B_{p,1}^{\frac{n}{p}}(\mathbb{R}^n)}^2
        \leq \mathcal{E}_0^\be,\nonumber
    \end{gather}
and \eqref{prop1.1-2-0108} for any $T>0$, where $\si(t)=\min\{1,t\}$, and
    \begin{align*}
        &n=2:q=\frac{4}{2-\eta}\in(2,4],\
        \al=\frac{1}{2(2-\eta)}\in\Big(\frac{1}{4},\frac{1}{2}\Big],\
        \be=\frac{1}{16};\\
        &n=3:q=3,\
        \al=\frac{3(2\eta-1)}{8\eta}\in\Big(\frac{3}{64},\frac{3}{8}\Big],\
        \be=\frac{(15\eta-8)}{192\eta}\in\Big(0,\frac{7}{192}\Big].
    \end{align*}
Moreover, the following large-time behavior holds:
    \begin{equation}
        \label{th1.2.3}
        \lim_{t\to\infty}\big(\|\mathbb F(\cdot,t)-\mathbb I\|_{L^{p_1}(\mathbb{R}^n)}
        +\|u(\cdot,t)\|_{L^{p_2}(\mathbb{R}^n)}
        +\|\na u(\cdot,t)\|_{L^{p_3}(\mathbb{R}^n)}\big)=0
    \end{equation}
for all $p_1\in(2,\infty)$, $p_2\in(2,\infty]$, and
    \begin{equation*}
        \begin{cases}
			\ p_3\in[2,\infty),&n=2,\\
            \ p_3\in[2,6),&n=3.
		\end{cases}
	\end{equation*}
\end{Theorem}

A few remarks on Theorem \ref{th1.2} are in order.
\begin{Remark}
We point out that
the $L^2(\mathbb{R}^n)$-norm of $(u_{0}, \mathbb F_{0}-\mathbb I)$
in \eqref{th1.2.2} is motivated the basic energy law (see Lemma \ref{le3.1}), while
the $\dot{B}_{p,1}^{\frac{n}{p}}(\mathbb{R}^n)(\subset L^\infty(\mathbb{R}^n))$-norm of $\mathbb F_{0}-\mathbb I$
takes care of the concentration of the deformation gradient (see Theorem \ref{th1.1}).
It follows from the classical interpolation theory (see Lemma \ref{ri}) that for $p$ satisfying \eqref{th1.2.0}-\eqref{th1.2.0-0413},
    \begin{equation*}
        \|\mathbb F_{0}-\mathbb I\|_{\dot{B}_{p,1}^{\frac{n}{p}}(\mathbb{R}^n)}
        \leq C\|\mathbb F_{0}-\mathbb I\|_{W^{1,p}(\mathbb{R}^n)}
        \leq C\|\mathbb F_{0}-\mathbb I\|_{H^2(\mathbb{R}^n)},
	\end{equation*}
which indicates that the smallness assumption \eqref{th1.2.2} on the initial data is much weaker than \eqref{prop1.1-2-0321} proposed in Proposition \ref{prop1.0}.
In particular, the oscillation of the initial velocity $u_0$ could be arbitrarily large in the sense of $L^\infty(\mathbb{R}^n)$-norm.
For example,
    \begin{equation}
        \label{re1.4-0413}
        u_0(x)=\va^\al\phi\left(\frac{x}{\va^\be}\right)\sin\left(\frac{x_1}{\va}\right),
        \quad -\frac{n\be}{2}<\al<-\frac{n\be}{3},\ 0<\be\leq1,
	\end{equation}
where $\phi(\cdot)\in\mathcal{S}(\mathbb R^n)$ is a Schwartz function
and $\va>0$ is small enough.
It is easy to verify that $u_0(x)$ in \eqref{re1.4-0413} belongs to $H^m(\mathbb{R}^n)$ ($m\geq2$) space and satisfies that
    \begin{equation*}
        \|u_0\|_{L^2(\mathbb{R}^n)}\ll1
	\end{equation*}
for sufficiently small $\va>0$.
Thus, Theorem \ref{th1.2} significantly improves the previous results \cite{Lin-Liu-Zhang2005,Lei-Zhou2005,Chen-Zhang2006,Lei-Liu-Zhou2008}, where they proved that the system \eqref{1.1} is globally well-posed for small initial data in the $H^2(\mathbb{R}^n)$ space.
\end{Remark}

\begin{Remark}
In our case,
the boundedness assumption \eqref{th1.2.1} on the $\dot{H}^\eta(\mathbb{R}^n)$-norm of the initial velocity is used to establish some improved short-time estimates of the velocity (see Lemmas \ref{le4.2} and \ref{le4.4}), which are important for dealing with the initial layer and for closing the required energy estimates.
Note that Theorem \ref{th1.2} allows one to construct global solutions for the highly oscillating initial velocity that may have a large norm in the homogeneous Sobolev space $\dot{H}^\frac{1}{2}(\mathbb{R}^n)$ or in $L^3(\mathbb{R}^n)$ (see \eqref{re1.4-0413} for instance).
This is in sharp contrast to the works of Fujita-Kato \cite{Fujita-Kato1964} and Kato \cite{Kato1984}, in which the global strong solutions to the 3D incompressible Navier-Stokes equations were obtained under the assumption that the initial data are small in $\dot{H}^\frac{1}{2}(\mathbb{R}^3)$ or $L^3(\mathbb{R}^3)$.
\end{Remark}

\begin{Remark}
In \cite{Hu-Lin2016}, Hu and Lin proved the global existence of weak solutions to the 2D incompressible viscoelastic system
under the assumptions that the initial deformation gradient is close to the identity matrix in $L^2(\mathbb{R}^2)\cap L^\infty(\mathbb{R}^2)$ and the initial velocity is small in $L^2(\mathbb{R}^2)$ and bounded in $L^p(\mathbb{R}^2)$ for some $p>2$.
However, the uniqueness of such solutions is unclear.
Based on the weak-strong uniqueness result in \cite[Theorem 2.2]{Hu-Wu2015}, our Theorem \ref{th1.2} can be regarded as the regularity and uniqueness theory of the weak solutions constructed in \cite{Hu-Lin2016}
when the initial perturbation of the deformation gradient is suitably small in the $\dot{B}_{p,1}^{\frac{2}{p}}(\mathbb{R}^2)$-norm.
\end{Remark}

In this paper, we study the Cauchy problem of the incompressible viscoelastic system \eqref{1.1} in both two and three dimensions.
The first novelty is that we establish some new blow-up criteria \eqref{th1.1.2}-\eqref{th1.1.2-1} for classical solutions, revealing that
the formation of finite-time singularities in incompressible viscoelastic flows can cause the deformation gradient to become unbounded.
Indeed, we extend the well-known Serrin-type criteria in terms of the velocity or pressure for the 3D incompressible Navier-Stokes equations to the 3D viscoelastic case; see \eqref{th1.1.2-2} and \eqref{th1.1.2-1}.
As a byproduct, we generalize the results for periodic problems in  \cite{Feng-Zhu-Zi2017}.
The second novelty is that we prove the global existence of classical solutions under the conditions that the initial velocity is small in the standard energy space $L^2(\mathbb{R}^n)$ and the initial deformation gradient is close to the identity matrix in $L^2(\mathbb{R}^n)\cap\dot{B}_{p,1}^{\frac{n}{p}}(\mathbb{R}^n)$, where $p$ satisfies \eqref{th1.2.0}-\eqref{th1.2.0-0413}.
In particular, the initial velocity may have large oscillations; see \eqref{re1.4-0413} for example.
Our work improves the previous results \cite{Lin-Liu-Zhang2005,Lei-Zhou2005,Chen-Zhang2006,Lei-Liu-Zhou2008} on global classical solutions for small $H^2(\mathbb{R}^n)$ perturbations near a constant equilibrium.
We then make some comments on the analysis in this paper.

{\it Comments on the proof of Theorem \ref{th1.1}}.
The main ingredient in proving Theorem \ref{th1.1} is to derive the estimates on the $L^\infty(0,T;L^4(\mathbb{R}^n))$-norm of $\na \mathbb{F}$ and the $L^1(0,T;L^\infty(\mathbb{R}^n))$-norm of $\na u$, which are crucial for establishing high-order energy estimates.
Indeed, a key observation is that the pressure $P$ satisfies the following Poisson equation (see \eqref{3.4.4})
    \begin{equation*}
        \De P
        =-\diver\diver(u\otimes u)
        +\diver\diver((\mathbb{F}-\mathbb{I})(\mathbb{F}^\top-\mathbb{I})),
	\end{equation*}
which, combined with
the boundedness of the Riesz transform on the $L^k(\mathbb{R}^n)$ space for $k\in(1,\infty)$, yields
    \begin{equation}
		\label{1.21}
        \|P\|_{L^k(\mathbb{R}^n)}
        \leq C\Big(\|u\|_{L^{2k}(\mathbb{R}^n)}^2
        +\|\mathbb{F}-\mathbb{I}\|_{L^{2k}(\mathbb{R}^n)}^2\Big).
	\end{equation}
By exploiting the relation \eqref{1.21}, we obtain an estimate on the $L^\infty(0,T;L^4(\mathbb{R}^n))$-norm of the velocity $u$ under the conditions that
the left-hand sides of \eqref{th1.1.2} (for $n=2$) and \eqref{th1.1.2-1} (for $n=3$) are finite; see Lemma \ref{le3.2}.
Note that the norm $\|u\|_{L^\infty(0,T;L^4(\mathbb{R}^3))}$ falls into Serrin's class (see \eqref{1.22}).
Therefore, to obtain the criteria \eqref{th1.1.2}-\eqref{th1.1.2-1}, it suffices to establish suitable a priori estimates for the classical solutions
under the assumption that
    \begin{equation}
		\label{1.22}
        \begin{cases}
			\ \lim\limits_{T\to T^*}
            \|\mathbb{F}\|_{L^\infty(0,T;L^\infty(\mathbb{R}^2))}<\infty,&n=2,\\
            \ \lim\limits_{T\to T^*}\left(\|\mathbb{F}\|_{L^\infty(0,T;L^\infty(\mathbb{R}^3))}
            +\|u\|_{L^{s_1}(0,T;L^{r_1}(\mathbb{R}^3))}\right)<\infty,&n=3
		\end{cases}
	\end{equation}
for any $r_1$ and $s_1$ satisfying
    \begin{equation*}
        \frac{2}{s_1}+\frac{3}{r_1}\leq1,\quad
        3< r_1\leq \infty.
	\end{equation*}

The first step is to derive the $L^\infty(0,T;L^2(\mathbb{R}^n))$-estimate on the gradient of the velocity.
We would like to point out that the method adopted here is different from that in \cite{Feng-Zhu-Zi2017}.
More precisely, Feng, Zhu, and Zi \cite{Feng-Zhu-Zi2017} introduced an important quantity $w$ defined by $w=u-v$, where
$v\triangleq L^{-1}\mathcal{P}\diver(\mathbb{F}\mathbb{F}^\top)$
is the solution of the following Poisson equation (with suitable boundary condition)
    \begin{equation}
		\label{1.23}
        -\mu\De v=\mathcal{P}\diver(\mathbb{F}\mathbb{F}^\top).
	\end{equation}
It follows from \eqref{1.1} and \eqref{1.23} that $w$ satisfies
    \begin{align}
		\label{1.24}
        w_t-\mu\De w
        =&-\mathcal{P}(u\cdot\na u)
        -L^{-1}\mathcal{P}\diver(\mathbb{F}\mathbb{F}^\top)_t\nonumber\\
        =&-\mathcal{P}(u\cdot\na u)
        +L^{-1}\mathcal{P}\diver\diver((\mathbb{F}\mathbb{F}^\top)\otimes u)
        -L^{-1}\mathcal{P}\diver(\na u(\mathbb{F}\mathbb{F}^\top)
        +(\mathbb{F}\mathbb{F}^\top)(\na u)^\top).
	\end{align}
Following the ideas in \cite{Sun-Wang-Zhang2011} for the compressible Navier-Stokes equations, Feng et al. multiplied the above equation \eqref{1.24} by $w_t$ and derived an estimate of $w$.
Using this together with the a priori assumptions and the regularity estimates for $v$ from \eqref{1.23} yields the desired estimates of $\na u$.
However, the elastic effect introduces a tricky term $L^{-1}\mathcal{P}\diver(\na u(\mathbb{F}\mathbb{F}^\top)+(\mathbb{F}\mathbb{F}^\top)(\na u)^\top)$ in \eqref{1.24}, where $L^{-1}\mathcal{P}\diver$ is not a zeroth-order operator.
To overcome this, they employed the Poincar\'{e} inequality so that their method could not deal with the Cauchy problem.
Returning to our case, we multiply the momentum equation $\eqref{1.1}_1$ by $u_t$ and utilize the following Poisson equation for $w$ (rather than \eqref{1.24})
    \begin{equation}
		\label{1.25}
        \mu\De w=\mathcal{P}\dot{u}.
	\end{equation}
With the decomposition of the velocity,
we conduct a detailed analysis of the elastic term $\int_{\mathbb{R}^n}\diver(\mathbb{F}\mathbb{F}^\top-\mathbb{I})\cdot u_t\,dx$ and obtain the estimate on the $L^\infty(0,T;L^2(\mathbb{R}^n))$-norm of $\na u$ (see Lemma \ref{le3.3}).
In particular, we are able to avoid the Poincar\'{e} inequality by exploiting the structure of the system \eqref{1.1} and the regularity estimates for the Poisson equations \eqref{1.23} and \eqref{1.25}.
The second step is to improve the regularity of the component $w$.
For this, we first deduce the basic estimate on the material derivatives of the velocity $\dot{u}$.
Unlike the projection method used in \cite[Lemma 4.1]{Feng-Zhu-Zi2017}, we will utilize the relation \eqref{1.21} to bound the pressure term (see Lemma \ref{le3.4}).
Then, combining the energy estimates of $\dot{u}$ with the regularity estimates of \eqref{1.25}, we get the required higher-order regularity for $w$ (see Corollary \ref{co3.1}).
Subsequently, by virtue of the a priori estimates obtained and a logarithmic Sobolev inequality (see Lemma \ref{le2.6}), we arrive at the key estimates on both the $L^\infty(0,T;L^4(\mathbb{R}^n))$-norm of $\na \mathbb{F}$ and the $L^1(0,T;L^\infty(\mathbb{R}^n))$-norm of $\na u$.
Finally, based on the product and commutator estimates (see Lemma \ref{le2.4}), the higher-order energy estimates for the solutions are derived via the standard energy method.
One may check Section \ref{sec3} for further details.


{\it Comments on the proof of Theorem \ref{th1.2}}.
As noted in Remark \ref{re1.1}, the criteria \eqref{th1.1.2}-\eqref{th1.1.2-1} obtained in Theorem \ref{th1.1} demonstrate that the concentration of the deformation gradient poses a major obstacle to establishing the global classical solutions for system \eqref{1.1}.
Therefore, how to derive the time-independent upper bound of the deformation gradient is the key issue in proving Theorem \ref{th1.2}.
Observing that the transport structure of the equation $\eqref{1.1}_2$ for the deformation gradient $\mathbb{F}$ is analogous to that of the continuity equation for the density $\rho$ in the compressible Navier-Stokes system, i.e.,
    \begin{equation*}
        \rho_t+u\cdot\na\rho=-\rho\diver u,
    \end{equation*}
we would like to draw on the methods developed in \cite{Hoff1995,Hoff2002,Huang-Li-Xin2012,Li-Xin2019} for bounding the density from above to establish an upper bound for the deformation gradient.
However, the density is a non-negative scalar-valued function, while the deformation gradient is a matrix-valued function with a positive determinant.
The two exhibit fundamental differences in both dimension and sign, which makes it difficult to directly utilize the strategies in \cite{Hoff1995,Hoff2002,Huang-Li-Xin2012,Li-Xin2019} to bound the deformation gradient.
To overcome this difficulty,
we introduce two so-called effective viscous fluxes as follows:
    \begin{equation*}
        \begin{cases}
            \ \mathcal{G}_1
            =\mu\nabla u-
            (-\De)^{-1}\na{\mathcal{P}}\diver(\mathbb{F}\mathbb{F}^{\top}-\mathbb{I}),\\
            \ \mathcal{G}_2
            =\mu\cur u-
            (-\De)^{-1}\cur\diver(\mathbb{F}\mathbb{F}^{\top}-\mathbb{I}),
        \end{cases}
    \end{equation*}
which suitably combine the effects of the velocity $u$ and the deformation gradient  $\mathbb{F}$.
An important fact is that the effective viscous fluxes $\mathcal{G}_1$ and $\mathcal{G}_2$ possess more regularity information than $u$ and $\mathbb{F}$ (see \eqref{ck}-\eqref{ck1}).
Combining the effective viscous flux $\mathcal{G}_2$ with the Helmholtz decomposition
    \begin{equation*}
        \Delta =\na\diver+\diver\cur,
    \end{equation*}
one can rewrite the transport equation $\eqref{1.1}_2$ for $\mathbb{F}$ as follows
    \begin{equation}
        \label{1.26}
        (\mathbb F-\mathbb I)_t+u\cdot\na(\mathbb F-\mathbb I)
        +\frac{1}{\mu}(\mathbb F-\mathbb I)
        =\frac{1}{\mu}\mathcal{G}\mathbb F+\frac{1}{\mu}\mathcal{R}\mathbb F
        -\frac{1}{\mu}(\mathbb F-\mathbb I)^2,
    \end{equation}
where
$\mathcal{G}\sim\mathcal{G}_2$ and $\mathcal{R}\sim(\mathbb F-\mathbb I)^2$
(see \eqref{324-0311} for the definitions).
Note that the term $\frac{1}{\mu}(\mathbb F-\mathbb I)$ in the equation \eqref{1.26} reveals the damping mechanism of the deformation gradient.
Moreover, thanks to the regularity of $\mathcal{G}_i$ ($i=1,2$), we can derive the uniform-in-time upper bound for the deformation gradient without using the $L^1(0,T;L^\infty(\mathbb{R}^n))$-norm of the gradient of the velocity.
In fact, since the Riesz transform is not well-suited for direct use in the $L^\infty(\mathbb{R}^n)$ space, we will study the upper bound on the $\dot{B}_{p,1}^{\frac{n}{p}}(\mathbb{R}^n)(\subset L^\infty(\mathbb{R}^n))$-norm of the deformation gradient (see also \cite{Hu2018}).
In addition, the effective viscous fluxes are also important for understanding the rates of regularization near the initial layer and for closing the required a priori estimates.

The main steps in the proof of Theorem \ref{th1.2} are as follows.
To begin with, standard energy estimates (see Lemma \ref{le4.1}) show that the quantities on the left-hand side of \eqref{th1.2-A1} can be bounded by the initial energy plus two higher-order nonlinear terms
    \begin{equation*}
        \mathcal{N}_1\triangleq
        \left|\int_0^T\int_{\mathbb{R}^{n}}\sigma \na u:(\na u\na u)\,dxdt\right|
        \quad\text{and}\quad
        \mathcal{N}_2\triangleq
        \int_0^T\int_{\mathbb{R}^{n}}\si^n|\na u|^4\,dxdt
    \end{equation*}
arising from the convection term $u\cdot\na u$.
In order to control the above two higher-order terms, we first establish an improved $L^q(\mathbb{R}^n)$-estimate of the velocity for small time (see Lemma \ref{le4.2}).
Then, utilizing the properties of the effective viscous fluxes (see \eqref{ck}-\eqref{ck1}), the fact that $\na u$ is a linear combination of $\mathcal{G}_1$ and $(-\De)^{-1}\na{\mathcal{P}}\diver(\mathbb{F}\mathbb{F}^{\top}-\mathbb{I})$,
and the Riesz-Thorin interpolation argument (see \cite{B-L1976}), we derive some useful time-weighted estimates on both the gradient and the material derivative of the velocity; see Lemma \ref{le4.4}.
With such time-weighted estimates at hand, we can use the equation \eqref{1.26} to obtain an estimate of $\int_0^T\int_{\mathbb{R}^{n}}\sigma^n|\mathbb F-\mathbb I|^4\,dxdt$, which combined with the effective viscous fluxes $\mathcal{G}_1$ and $\mathcal{G}_2$ yields the required bound for
$\mathcal{N}_2$ (see \eqref{4.3.2}).
As for $\mathcal{N}_1$,
we point out that the estimate for $n=3$ is quite different from the treatment in \cite{Hu-Lin2016,Hu2018} for $n=2$.
Indeed, it follows Cauchy's inequality that
    \begin{equation*}
        \mathcal{N}_1
        \leq C\int_0^T\int_{\mathbb{R}^{n}}|\nabla u|^2\,dxdt
        +C\int_0^T\int_{\mathbb{R}^{n}}\sigma^2|\nabla u|^4\,dxdt.
    \end{equation*}
When $n=2$, the right-hand side of the above inequality can be bounded by the initial energy and the estimate of $\mathcal{N}_2$.
Unfortunately, when $n=3$ and $T\leq1$, one cannot use $\mathcal{N}_2\left(=\int_0^T\int_{\mathbb{R}^{3}}\si^3|\na u|^4\,dxdt\right)$ to control $\int_0^T\int_{\mathbb{R}^{3}}\sigma^2|\nabla u|^4\,dxdt$ directly.
To address this issue, we conduct a detailed estimate of $\mathcal{N}_1$ for small time (see \eqref{4.3.12}) by using the effective viscous flux $\mathcal{G}_1$ and time-weighted estimates obtained in Lemma \ref{le4.4}.
Collecting the bounds of the higher-order nonlinear terms $\mathcal{N}_1$ and $\mathcal{N}_2$, we derive the desired estimate \eqref{th1.2-A1} (see Lemma \ref{le4.3} for the details).
Furthermore, applying these key estimates and some well-known interpolation results (see Lemmas \ref{ri} and \ref{le2.3}), we can obtain the time-independent estimates on the
$L^{1}\left(0,\min\{1,T\};\dot{B}_{p,1}^{\frac{n}{p}}(\mathbb{R}^n)\right)$-norm
and $L^{\tilde{p}}\left(\min\{1,T\},T;\dot{B}_{p,1}^{\frac{n}{p}}(\mathbb{R}^n)\right)$-norm
of $\mathcal{G}_i$ ($i=1,2$), where $\tilde{p}=\frac{4p^2}{n^2p-2n^2}>1$.
In fact, since we do not require the $\dot{H}^\eta(\mathbb{R}^n)$-norm of the initial velocity to be small, we obtain only the boundedness (not smallness) of the time-weighted estimate for both the gradient and the material derivative of the velocity (see \eqref{4.4.1} and \eqref{4.4.1-0406}), and the lack of smallness makes
the estimate on the $L^{1}\left(0,\min\{1,T\};\dot{B}_{p,1}^{\frac{n}{p}}(\mathbb{R}^n)\right)$-norm
of the effective viscous fluxes (see \eqref{4.5.9})
more involved than that in \cite[Lemma 4.1]{Hu2018}.
Based on the equation \eqref{1.26}, the Littlewood-Paley theory, and some careful initial-layer analysis, we arrive at the uniform bound (with respect to time) on the $\dot{B}_{p,1}^{\frac{n}{p}}(\mathbb{R}^n)$-norm of the deformation gradient (see Lemma \ref{le4.5}).
Finally, combining these a priori estimates with the blow-up criteria obtained in Theorem \ref{th1.1}, we prove the global existence of classical solutions by a continuity argument.
See Section \ref{sec4} for more information.

The remaining part of this paper is organized as follows.
In Section \ref{sec2}, we first recall the Littlewood-Paley theory, and then present some important facts and inequalities that will be used in the later analysis.
Section \ref{sec3} is devoted to establishing some blow-up criteria of classical solutions to the problem \eqref{1.1}-\eqref{1.2}.
Finally, we investigate the global-in-time existence of classical solutions in Section \ref{sec4}.

\section{Preliminaries}\label{sec2}
\subsection{Littlewood-Paley theory}\

Let us recall some basic facts on the Littlewood-Paley theory; see \cite{Bahouri-Chemin-Danchin2011} for more details.
We first introduce the Littlewood-Paley decomposition.
Let $n=2$ or $3$.
The set of tempered distributions is denoted by $\mathcal S'(\mathbb R^n)$, which is the dual space of the Schwartz space $\mathcal{S}(\mathbb R^n)$.
Choose two smooth radial functions $\varphi$ and $\chi$, valued in the interval $[0,1]$, such that
    \begin{gather*}
        \sum_{j \in \mathbb{Z}} \varphi\Big(2^{-j} \xi\Big)=1,
        \quad\forall\,\xi\in \mathbb{R}^{n}\backslash\{0\},\\
        \chi(\xi)+\sum_{j\in\mathbb{N}}\varphi\Big(2^{-j}\xi\Big)=1,
        \quad\forall\,\xi\in\mathbb{R}^n,
    \end{gather*}
with
    \begin{equation*}
        \operatorname{supp} \varphi(\xi)\subset\bigg\{\xi\in{\mathbb R}^n\ \Big|\  \frac34\leq|\xi|\leq\frac83\bigg\}
        \quad\text{and}\quad
        \operatorname{supp} \chi(\xi)\subset\bigg\{\xi\in{\mathbb R}^n\ \Big|\  |\xi|\leq\frac43\bigg\}.
    \end{equation*}
For any $f\in\mathcal S'(\mathbb R^n)$, the homogeneous dyadic blocks $\dot\Delta_j$ and the homogeneous low-frequency cut-off operators $\dot S_j$ are defined for all $j\in\mathbb Z$ by
    \begin{gather*}
        \dot{\Delta}_{j}f
        =\varphi(2^{-j}D)f
        =2^{jn}\int_{\mathbb{R}^n}\check{\varphi}(2^jy)f(x-y) dy,\\
        \dot{S}_{j}f
        =\chi(2^{-j}D)f
        =2^{jn}\int_{\mathbb{R}^n}\check{\chi}(2^jy)f(x-y) dy,
    \end{gather*}
where $\check{\varphi}$ and $\check{\chi}$ are the inverse Fourier transform of $\varphi$ and $\chi$, respectively.
We denote by $\mathcal S'_h(\mathbb R^n)$ the space of temperate distributions $f\in\mathcal S'(\mathbb R^n)$ with
    \begin{equation*}
        \lim\limits_{j\to-\infty}\|\dot{S}_jf\|_{L^\infty(\mathbb{R}^n)}=0.
    \end{equation*}
Then it holds that
    \begin{equation*}
        f=\sum_{j\in\mathbb{Z}}\dot{\Delta}_jf,\quad \forall\, f\in \mathcal S'_h(\mathbb R^n).
    \end{equation*}
Moreover, by the definitions of $\varphi$ and $\chi$, one has
    \begin{equation}
        \label{0}
        \begin{cases}
            \ \dot{\Delta}_{j'}\dot{\Delta}_{j}f=0,
            &\text{if}\ |j'-j|\geq2,\\
            \ \dot{\Delta}_{j'}(\dot S_{j-1}f\dot{\Delta}_{j}f)=0,
            &\text{if}\ |j'-j|\geq5.
		\end{cases}
    \end{equation}

Now we can give the definition of homogeneous Besov spaces as follows.
\begin{Definition}\label{def2.1}
Let $s\in \mathbb R$ and $1\leq p,r\leq\infty$. The homogeneous Besov space $\dot{B}_{p,r}^{s}=\dot B_{p,r}^s(\mathbb R^n)$ is defined by
    \begin{equation*}
        \dot{B}_{p,r}^{s}(\mathbb{R}^n)
        =\{f\in \mathcal S'_h(\mathbb R^n)\mid
        \|f\|_{\dot{B}^{s}_{{p,r}}(\mathbb{R}^n)}<\infty\},
    \end{equation*}
where
    \begin{equation*}
        \|f\|_{\dot B_{p,r}^s(\mathbb{R}^n)}
        \triangleq\Bigg(\sum_{j\in\mathbb Z}2^{rjs}
        \|\dot\Delta_jf\|_{L^p(\mathbb{R}^n)}^r\Bigg)^{\frac1r}.
    \end{equation*}
\end{Definition}

In order to bound the nonlinearity in Besov spaces, we also introduce the homogeneous Bony's decomposition
    \begin{equation}
        \label{bony}
        fg
        =\sum_{j',j\in\mathbb{Z}}\dot{\Delta}_{j'}f \dot{\Delta}_jg
        =\dot{T}_fg+\dot{T}_gf+\dot{R}(f,g),
    \end{equation}
where the homogeneous paraproduct $\dot{T}_fg$ and remainder $\dot{R}(f,g)$ of $f$ and $g$ are defined as follows:
    \begin{equation}
        \label{bony-1}
        \dot{T}_fg=\sum_{j\in\mathbb{Z}}\dot{S}_{j-1}f \dot{\Delta}_jg
        \quad\text{and}\quad
        \dot{R}(f,g)=
        \sum_{j',j\in\mathbb{Z},\ |j'-j|\leq 1}
        \dot{\Delta}_{j'}f \dot{\Delta}_{j}g.
    \end{equation}

Below, we list some useful lemmas which will be used in the later analysis.
The detailed proofs of the following Lemmas \ref{ber}-\ref{jh} can be found in \cite[Chapter 2]{Bahouri-Chemin-Danchin2011}.

First, we show the Bernstein inequality.
\begin{Lemma}\label{ber}
Let $R_1,R_2,R_3>0$, $j\in\mathbb Z$, $k\in\mathbb N$, and $1 \leq p \leq q \leq \infty$.
If $f \in L^{p}(\mathbb{R}^{n})$, then
    \begin{gather*}
        \operatorname{supp}\mathcal F(f)\subset\{|\xi|\leq R_12^j\}
        \Rightarrow
        \sup\limits_{|\alpha|=k}\|\partial^\alpha f\|_{L^q(\mathbb{R}^n)}
        \leq C_12^{jk+nj(\frac1p-\frac1q)}\|f\|_{L^p(\mathbb{R}^n)},\\
        \operatorname{supp}\mathcal F(f)\subset\{R_22^{j}\leq|\xi|\leq R_32^{j}\}
        \Rightarrow
        C_{2}^{-1}2^{jk}\|f\|_{L^p(\mathbb{R}^n)}
        \leq \sup_{|\alpha|=k}\|\partial^{\alpha}f\|_{L^p(\mathbb{R}^n)(\mathbb{R}^n)}
        \leq C_22^{jk}\|f\|_{L^p(\mathbb{R}^n)},
    \end{gather*}
where the positive constants $C_1$ and $C_2$ are independent of $f$ and $j$.
\end{Lemma}

The following two lemmas are concerned with the continuity of the remainder homogeneous paraproduct and remainder operators.
\begin{Lemma}\label{TB0}
Let $s\in\mathbb{R}$ and $1\leq p,r\leq \infty$.
If $f\in L^\infty(\mathbb{R}^n)$ and $g\in B_{p,r}^s(\mathbb{R}^n)$,
then there exists a positive constant $C$ depending only on $n$ and $s$ such that
    \begin{equation*}
        \|\dot{T}_{f}g\|_{\dot B_{p,r}^s(\mathbb{R}^n)}
        \leq C\|f\|_{L^\infty}\|g\|_{\dot B_{p,r}^s(\mathbb{R}^n)}.
	\end{equation*}
\end{Lemma}

\begin{Lemma}\label{TB1}
Let $s_1,s_2>0$, and $1\leq p,p_1,p_2,r,r_1,r_2\leq \infty$ with
    \begin{equation*}
        \frac1{p}=\frac1{p_1}+\frac1{p_2}
        \quad\text{and}\quad
        \frac1{r}=\frac1{r_1}+\frac1{r_2}.
	\end{equation*}
If $f\in\dot B_{p_1,r_1}^{s_1}(\mathbb{R}^n)$ and $g\in\dot B_{p_2,r_2}^{s_2}(\mathbb{R}^n)$,
then there exists a positive constant $C$ depending only on $n$, $s_1$, and $s_2$ such that
    \begin{equation*}
        \|\dot{R}(f,g)\|_{\dot B_{p,r}^{s_1+s_2}(\mathbb{R}^n)}
        \leq C\|f\|_{\dot B_{p_1,r_1}^{s_1}(\mathbb{R}^n)}\|g\|_{\dot B_{p_2,r_2}^{s_2}(\mathbb{R}^n)}.
	\end{equation*}
\end{Lemma}

Using Bony's decomposition \eqref{bony} and Lemmas \ref{TB0}-\ref{TB1},
one has the following so-called tame estimates for the product of two functions in Besov spaces.
\begin{Lemma}\label{TB2}
Let $s>0$ and $1\leq p,r\leq \infty$.
If $(f,g)\in L^\infty(\mathbb{R}^n)\cap\dot B_{p,r}^s(\mathbb{R}^n)$, then there exists a positive constant $C$ depending only on $n$ and $s$ such that
    \begin{equation*}
        \|fg\|_{\dot B_{p,r}^s(\mathbb{R}^n)}
        \leq C\Big(\|f\|_{L^\infty(\mathbb{R}^n)}\|g\|_{\dot B_{p,r}^s(\mathbb{R}^n)}
        +\|g\|_{L^\infty(\mathbb{R}^n)}\|f\|_{\dot B_{p,r}^s(\mathbb{R}^n)}\Big).
	\end{equation*}
\end{Lemma}

The following lemma gives a classical commutator estimate.
\begin{Lemma}\label{jh}
Let $1\leq p,q\leq\infty$.
Assume that $f\in C^{1}(\mathbb{R}^{n})$ with $(1+|\cdot|)\mathcal{F}(f) \in L^{1}(\mathbb{R}^{n})$,
$g$ is a Lipschitz function with $\na g\in L^{p}(\mathbb{R}^{n})$, and $h\in L^{q}(\mathbb{R}^{n})$.
Define
    \begin{equation*}
        [f(\lambda^{-1} D), g] h
        \triangleq f(\lambda^{-1} D)(gh)-gf(\lambda^{-1} D)h
    \end{equation*}
with
    \begin{equation*}
        f(D)h=\mathcal F^{-1}(f\mathcal{F}(h)).
    \end{equation*}
For any $\lambda>0$, there exists a positive constant $C$ independent of $g$ and $h$ such that
    \begin{equation*}
        \|[f(\lambda^{-1} D), g] h\|_{L^{r}(\mathbb{R}^n)}
        \leq C \lambda^{-1}\|\nabla g\|_{L^{p}(\mathbb{R}^n)}\|h\|_{L^{q}(\mathbb{R}^n)}
         \quad\text{with}\quad
         \frac{1}{r}=\frac{1}{p}+\frac{1}{q}.
    \end{equation*}
\end{Lemma}

Finally, we present a well-known interpolation result, which relates the homogeneous Sobolev and Besov Spaces.
\begin{Lemma}\label{ri}
Let $s,s_1,s_2\in\mathbb{R}$ and $1\leq p,r\leq\infty$.
If $s=(1-\theta)s_1+\theta s_2$ for some $\theta\in(0,1)$, then
    \begin{equation*}
        (\dot H_p^{s_1}(\mathbb{R}^n),\dot H_p^{s_2}(\mathbb{R}^n))_{\theta,r}
        =\dot B_{p,r}^s(\mathbb{R}^n),
    \end{equation*}
where $(\dot H_p^{s_1}(\mathbb{R}^n),\dot H_p^{s_2}(\mathbb{R}^n))_{\theta,r}$ denotes the real interpolation space between $\dot H_p^{s_1}(\mathbb{R}^n)$ and $\dot H_p^{s_2}(\mathbb{R}^n)$.
\end{Lemma}

\begin{proof}
See \cite[Theorem 6.3.1]{B-L1976}.
\end{proof}

\subsection{Elementary facts and auxiliary inequalities}\

We first recall some intrinsic properties of the viscoelastic flows.
\begin{Lemma}
	\label{le2.2}
Assume that $(u,\mathbb{F})$ is the classical solution of the problem \eqref{1.1}-\eqref{1.2} in $\mathbb{R}^n$ ($n=2,3$).
Then the following identities
    \begin{equation}
		\label{2.2.1}
        \begin{cases}
            \ \det\mathbb{F}=1,\\
            \ \mathbb{F}^{lk}\na_l\mathbb{F}^{ij}
            =\mathbb{F}^{lj}\na_l\mathbb{F}^{ik},\quad
            1\leq i,j,k\leq n,\\
            \ \diver\mathbb{F}^{\top}=0
		\end{cases}
	\end{equation}
hold for all time $t>0$ if they are initially valid.
\end{Lemma}

\begin{proof}
Refer to \cite{Lei-Liu-Zhou2008,Qian-Zhang2010}.
Note that the three identities in \eqref{2.2.1} hold independently.
\end{proof}

In fact, the three identities in \eqref{2.2.1} are related as follows.
\begin{Lemma}
    \label{leH}
Assume that $(u,\mathbb{F})$ is the classical solution of the problem \eqref{1.1}-\eqref{1.2} in $\mathbb{R}^n$ ($n=2,3$).
It holds that for $1\leq i,j,k\leq n$,
    \begin{equation*}
        \left.
        \begin{aligned}
            \det\mathbb{F}=&\ 1\\
            \mathbb{F}^{lk}\na_l\mathbb{F}^{ij}
            =&\ \mathbb{F}^{lj}\na_l\mathbb{F}^{ik}\
		\end{aligned}
        \right\}
        \Rightarrow
        \diver\mathbb{F}^{\top}=0,
	\end{equation*}
and
    \begin{equation*}
        \left.
        \begin{aligned}
            \diver\mathbb{F}^{\top}=&\ 0\\
            \mathbb{F}^{lk}\na_l\mathbb{F}^{ij}
            =&\ \mathbb{F}^{lj}\na_l\mathbb{F}^{ik}\
		\end{aligned}
        \right\}
        \Rightarrow
        \det\mathbb{F}=1.
	\end{equation*}
\end{Lemma}

\begin{proof}
See \cite[Appendix A]{Hu-Zhao2020-1} and \cite[Appendix A]{Hu-Zhao2020}.
\end{proof}

We then present the classical Gagliardo-Nirenberg inequality, whose proof can be found in \cite{Nirenberg1959}.
\begin{Lemma}
    \label{le2.3}
Let $m>0$ be some positive integer and $1\leq r,q \leq\infty$.
Suppose $f\in L^r(\mathbb R^n)\cap D^{m,q}(\mathbb R^n)$.
Then we have for any integer $k\in[0,m)$,
    \begin{equation}
        \label{2.3.1}
		\|\na^k f\|_{L^{p}(\mathbb{R}^n)}
        \leq C\|\na^m f\|_{L^{q}(\mathbb{R}^n)}^{\theta}\|f\|_{L^{r}(\mathbb{R}^n)}^{1-\theta},
	\end{equation}
where
    \begin{equation*}
    \frac{1}{p}
    =\frac{k}{n}+\Big(\frac{1}{q}-\frac{m}{n}\Big)\theta+\frac{1}{r}(1-\theta)
	\end{equation*}
with
    \begin{equation*}
        \begin{cases}
            \ \theta\in[\frac{k}{m},1),&\text{if}\ 1<r<\infty\ \text{and}\ m-k-\frac{n}{q}\ \text{is a nonnegative integer},\\
            \ \theta\in[\frac{k}{m},1],&\text{otherwise}.
		\end{cases}
	\end{equation*}
Moreover, if $r=\infty$, $qm<n$ and $k=0$, then one additionally requires that either $f(x)\to 0$ as $|x|\to \infty$ or $|f|^{\tilde{r}}\in L^1(\mathbb{R}^n)$ for some $0<\tilde{r}< \infty$.
The constant $C>0$ in \eqref{2.3.1} depends only on $n$, $m$, $r$, $q$, $k$, and $\theta$.
\end{Lemma}


Below, we state some regularity results for the Poisson system.
\begin{Lemma}
	\label{le2.5}
Let $U=U(x)\in H^1(\mathbb R^n)$ be a weak solution to the following Poisson system
    \begin{equation*}
        \begin{cases}
			\ \De U=F,&x\in\mathbb{R}^n,\\
			\ U(x)\to0,&|x|\to\infty.
		\end{cases}
	\end{equation*}
Then, for $1<p<\infty$, there exist a positive constant $C$ depending only on $n$ and $p$ such that

\noindent (\romannumeral1) if $F\in L^p(\mathbb R^n)$, then
    \begin{equation}
        \label{2.5.1}
        \|\na^2 U\|_{L^p(\mathbb{R}^n)}\leq C\|F\|_{L^p(\mathbb{R}^n)};
	\end{equation}

\noindent (\romannumeral2) if $F=\diver f$ with $f=(f^{ij})_{n\times n}$ and $f^{ij}\in L^p(\mathbb R^n)$, then
    \begin{equation}
        \label{2.5.2}
        \|\na U\|_{L^p(\mathbb{R}^n)}\leq C\|f\|_{L^p(\mathbb{R}^n)};
	\end{equation}


\noindent (\romannumeral3) if $F=\diver f$ with $f=(f^{ij})_{n\times n}$ and $f_{ij}\in L^2(\mathbb R^n)\cap L^\infty(\mathbb R^n)$, then
    \begin{equation}
        \label{2.5.3}
        \|\na U\|_{BMO(\mathbb{R}^n)}
        \leq C(\|f\|_{L^2(\mathbb{R}^n)}+\|f\|_{L^\infty(\mathbb{R}^n)}).
	\end{equation}
\end{Lemma}

\begin{proof}
    The proof is similar to \cite[Propositions 2.1-2.2]{Sun-Wang-Zhang2011}.
\end{proof}

%
%
%
%
%

The following Sobolev inequality of logarithmic type, which can be found in \cite[Theorem 1]{Kozono-Taniuchi2000}, will be used later to estimate $\|\na \mathbb{F}\|_{L^\infty(0,T; L^4(\mathbb{R}^n))}$ and $\|\na u\|_{L^1(0,T; L^\infty(\mathbb{R}^n))}$.
\begin{Lemma}
	\label{le2.6}
Assume that $f\in W^{1,p}(\mathbb R^n)$ with $n<p<\infty$.
Then there exists a positive constant $C$ depending only on $n$ and $p$ such that
    \begin{equation}
        \label{2.6.1}
        \|f\|_{L^\infty(\mathbb{R}^n)}
        \leq C\big[1+\|f\|_{BMO(\mathbb{R}^n)}\ln(e+\|\na f\|_{L^p(\mathbb{R}^n)})\big].
	\end{equation}
\end{Lemma}

The following useful inequality for the Hardy space $\mathcal{H}^1({\mathbb{R}^n})$ is given in \cite{Coifman-Lions-Meyer-Semmes1993}.
\begin{Lemma}\label{RP}
Assume that $f\in L^2({\mathbb{R}^n})$ with $\diver f=0$, and $g\in H^1({\mathbb{R}^n})$.
Then there exists a positive constant $C$ depending only on $n$ such that
    \begin{equation*}
        \|f\cdot\nabla g\|_{\mathcal{H}^1(\mathbb{R}^n)}
        \leq C\|f\|_{L^2(\mathbb{R}^n)}\|\nabla g\|_{L^2(\mathbb{R}^n)}.
    \end{equation*}
\end{Lemma}


%


Finally, we recall the following product and commutator estimates, which will be used in the higher-order energy estimates.
\begin{Lemma}
    \label{le2.4}
Assume that $1\leq p,p_1,p_2,p_3,p_4\leq\infty$ and
    \begin{equation*}
		\frac{1}{p}
        =\frac{1}{p_{1}}+\frac{1}{p_{2}}=\frac{1}{p_{3}}+\frac{1}{p_{4}}.
	\end{equation*}
Then it holds that

\noindent (\romannumeral1) if $f\in L^{p_1}(\mathbb R^n)\cap D^{k,p_3}(\mathbb R^n)$ and $g\in L^{p_4}(\mathbb R^n)\cap D^{k,p_2}(\mathbb R^n)$ for some $k\geq0$, then
    \begin{equation*}
		\|\nabla^k(fg)\|_{L^{p}(\mathbb{R}^n)}
        \leq C\Big(\|f\|_{L^{p_1}(\mathbb{R}^n)}\|\nabla^k g\|_{L^{p_2}(\mathbb{R}^n)}
        +\|\nabla^k f\|_{L^{p_3}(\mathbb{R}^n)}\|g\|_{L^{p_4}(\mathbb{R}^n)}\Big);
	\end{equation*}

\noindent (\romannumeral2) if $f\in D^{1,p_1}(\mathbb R^n)\cap D^{l,p_3}(\mathbb R^n)$ and $g\in L^{p_4}(\mathbb R^n)\cap D^{l-1,p_2}(\mathbb R^n)$ for some $l\geq1$, then
    \begin{equation*}
        \|\nabla^l(fg)-f\nabla^l g\|_{L^{p}(\mathbb{R}^n)}
        \leq C\Big(\|\nabla f\|_{L^{p_1}(\mathbb{R}^n)}\|\nabla^{l-1} g\|_{L^{p_2}(\mathbb{R}^n)}
        +\|\nabla^l f\|_{L^{p_3}(\mathbb{R}^n)}\|g\|_{L^{p_4}(\mathbb{R}^n)}\Big),
	\end{equation*}
where the positive constant $C$ is independent of $f$ and $g$.
\end{Lemma}

\begin{proof}
Refer to \cite[Lemma 3.1]{Ju2004} and \cite[Lemma 3.4, p. 129]{Majda-Bertozzi2002}.
\end{proof}

\section{Blow-up Criteria}\label{sec3}
In this section, we study the blow-up criteria of classical solutions to the Cauchy problem \eqref{1.1}-\eqref{1.2}.
Assume that $(u,\mathbb{F})$ is the classical solution to \eqref{1.1}-\eqref{1.2} on $[0,T^*)$, where $T^*<\infty$ is the maximum existence time.

First of all, we have the following elementary energy estimate.
\begin{Lemma}\label{le3.1}
Under the conditions of Theorem \ref{th1.1}, it holds that for any $T\in[0,T^*)$,
    \begin{equation}
		\label{3.1.1}
        \sup\limits_{0\leq t\leq T}\int_{\mathbb{R}^{n}}\Big(|u|^2+|\mathbb F-\mathbb I|^2\Big)\,dx
        +2\mu\int_0^T\int_{\mathbb{R}^{n}}|\nabla u|^2\,dxdt
        \leq\int_{\mathbb{R}^{n}}\Big(|u_0|^2+|\mathbb F_0-\mathbb I|^2\Big)\,dx.
	\end{equation}
\end{Lemma}

\begin{proof}
Multiplying the first two equations in \eqref{1.1} by $u$ and $\mathbb{F}-\mathbb{I}$, respectively, summing up, and then integrating the resulting equation over $\mathbb{R}^n$ by parts, by $\diver u=0$ and $\diver \mathbb{F}^\top=0$, we deduce
    \begin{align*}
        &\ \frac{1}{2}\frac{d}{dt}\int_{\mathbb{R}^n}\Big(|u|^2
        +|\mathbb{F}-\mathbb{I}|^2\Big)\,dx
        +\mu\int_{\mathbb{R}^n}|\na u|^2\,dx\\
        =&\ \int_{\mathbb{R}^n}\diver(\mathbb{F}\mathbb{F}^\top)\cdot u\,dx
        +\int_{\mathbb{R}^n}(\na u\mathbb{F}):(\mathbb{F}-\mathbb{I})\,dx\\
        =&\ \int_{\mathbb{R}^n}\diver((\mathbb{F}-\mathbb{I})\mathbb{F}^\top)\cdot u\,dx
        +\int_{\mathbb{R}^n}\na u:(\mathbb{F}-\mathbb{I})\mathbb{F}^\top\,dx
        =0.
	\end{align*}
Integrating the above equation over $[0,t]$ with $0\leq t\leq T$ gives \eqref{3.1.1}.
\end{proof}

Next, we give an estimate on the $L^\infty(0,T;L^4(\mathbb{R}^n))$-norm of $u$ under the assumption that
    \begin{equation}
		\label{3.1-0406}
		\mathcal{K}\triangleq
        \begin{cases}
			\ \lim\limits_{T\to T^*}
            \|\mathbb{F}\|_{L^\infty(0,T;L^\infty(\mathbb{R}^2))}<\infty,&n=2,\\
            \ \lim\limits_{T\to T^*}
            \left(\|\mathbb{F}\|_{L^\infty(0,T;L^\infty(\mathbb{R}^3))}
            +\|P\|_{L^{s_2}(0,T;L^{r_2}(\mathbb{R}^3))}\right)<\infty,&n=3,
		\end{cases}
	\end{equation}
where $r_2$ and $s_2$ satisfy
    \begin{equation}
		\label{th1.1.3-0406}
        \frac{2}{s_2}+\frac{3}{r_2}\leq2,\quad
        \frac{3}{2}< r_2\leq \infty.
	\end{equation}
\begin{Lemma}\label{le3.2}
    Under the conditions of Theorem \ref{th1.1} and \eqref{3.1-0406}, it holds that for any $T\in[0,T^*)$,
    \begin{equation}
		\label{3.2.1}
        \sup\limits_{0\leq t\leq T}\|u\|_{L^4(\mathbb{R}^n)}
        +\int_0^{T}\int_{\mathbb{R}^n}\Big(|\na u|^2|u|^{2}+|\na|u|^2|^2\Big)\,dxdt
        \leq \widetilde{C},
	\end{equation}
where $\widetilde{C}>0$ depends only on $\mathcal{K}$, $T^*$, $n$, $\mu$,
$\|u_0\|_{L^2(\mathbb{R}^n)\cap L^4(\mathbb{R}^n)}$, and $\|\mathbb{F}_0-\mathbb{I}\|_{L^2(\mathbb{R}^n)}$.
\end{Lemma}

\begin{proof}
Multiplying $\eqref{1.1}_1$ by $4|u|^2u$, and integrating the resulting equation over $\mathbb{R}^n$, by integration by parts and $\diver u=0$, we have
    \begin{align}
		\label{3.2.2}
        & \frac{d}{dt}\int_{\mathbb{R}^n}|u|^4\,dx
        +\mu\int_{\mathbb{R}^n}\Big(4|\na u|^2|u|^{2}+2|\na|u|^2|^2\Big)\,dx\nonumber\\
        =& -4\int_{\mathbb{R}^n}(\mathbb{F}\mathbb{F}^\top):\na(|u|^2u)\,dx
        +4\int_{\mathbb{R}^n}P u\cdot\na |u|^2\,dx.
	\end{align}
For the first term on the right-hand side of \eqref{3.2.2}, since $|\nabla (|u|^2u)|\leq \widetilde{C}|\nabla u ||u|^2$, we deduce from H\"{o}lder's and Cauchy's inequalities that
    \begin{align}
		\label{3.2.3}
        -4\int_{\mathbb{R}^n}(\mathbb{F}\mathbb{F}^\top):\na(|u|^2u)\,dx
        \leq \widetilde{C}\|\mathbb{F}\|_{L^\infty(\mathbb{R}^n)}^2
        \|u\|_{L^2(\mathbb{R}^n)}\big\||\na u||u|\big\|_{L^2(\mathbb{R}^n)}
        \leq \mu\big\||\na u||u|\big\|_{L^2(\mathbb{R}^n)}^2+\widetilde{C}.
	\end{align}

We then estimate the second term on the right-hand side of \eqref{3.2.2}.
Applying the $\diver$ operator to $\eqref{1.1}_1$ yields
    \begin{equation}
		\label{3.4.4}
        \De P
        =-\diver(u\cdot\na u)+\diver\diver(\mathbb{F}\mathbb{F}^\top)
        =-\diver\diver(u\otimes u)
        +\diver\diver((\mathbb{F}-\mathbb{I})(\mathbb{F}^\top-\mathbb{I})),
	\end{equation}
where one has used $\diver u=0$ and the following identity
    \begin{align*}
        \diver\diver(\mathbb{F}\mathbb{F}^\top)
        &= \diver\diver((\mathbb{F}-\mathbb{I})(\mathbb{F}^\top-\mathbb{I}))
        +\diver\diver(\mathbb{F}-\mathbb{I})+\diver\diver(\mathbb{F}^\top-\mathbb{I})\nonumber\\
        &= \diver\diver((\mathbb{F}-\mathbb{I})(\mathbb{F}^\top-\mathbb{I}))
	\end{align*}
due to $\diver\mathbb{F}^\top=0$.
With the equation \eqref{3.4.4}, employing
$L^p$-boundedness of the Riesz transform,
the standard $L^p$ interpolation inequality, \eqref{3.1.1}, and \eqref{3.1-0406}, we get for any $k\in(1,\infty)$
    \begin{equation}
		\label{3.4.4-0336}
        \|P\|_{L^k(\mathbb{R}^n)}
        \leq \widetilde{C}\Big(\|u\|_{L^{2k}(\mathbb{R}^n)}^2
        +\|\mathbb{F}\|_{L^{2k}(\mathbb{R}^n)}^2\Big)
        \leq \widetilde{C}\Big(\|u\|_{L^{2k}(\mathbb{R}^n)}^2
        +\|\mathbb{F}\|_{L^2(\mathbb{R}^n) \cap L^\infty(\mathbb{R}^n)}^2\Big)
        \leq \widetilde{C}\Big(\|u\|_{L^{2k}(\mathbb{R}^n)}^2+1\Big).
	\end{equation}

For $n=2$, it follows from \eqref{3.4.4-0336}, H\"{o}lder's and Young's inequalities that
    \begin{align}
		\label{3.2.4-0336}
        4\int_{\mathbb{R}^2}P u\cdot\na |u|^2\,dx
        \leq&\ \widetilde{C}\|P\|_{L^3(\mathbb{R}^2)}\|u\|_{L^6(\mathbb{R}^2)}
        \big\|\na|u|^2\big\|_{L^2(\mathbb{R}^2)}\nonumber\\
        \leq&\ \frac{\mu}{2}\big\|\na|u|^2\big\|_{L^2(\mathbb{R}^2)}^2
        +\widetilde{C}\|P\|_{L^3(\mathbb{R}^2)}^2\|u\|_{L^6(\mathbb{R}^2)}^2\nonumber\\
        \leq&\ \frac{\mu}{2}\big\|\na|u|^2\big\|_{L^2(\mathbb{R}^2)}^2
        +\widetilde{C}\|u\|_{L^6(\mathbb{R}^2)}^6+\widetilde{C}.
	\end{align}
Moreover, Lemma \ref{le2.3} gives
    \begin{align*}
        \|u\|_{L^6(\mathbb{R}^2)}^{6}
        =&\ \big\||u|^2\big\|_{L^3(\mathbb{R}^2)}^{3}
        \leq \widetilde{C}\big\||u|^2\big\|_{L^2(\mathbb{R}^2)}^{2}
        \big\|\na|u|^2\big\|_{L^2(\mathbb{R}^2)}
        =\widetilde{C}\|u\|_{L^4(\mathbb{R}^2)}^{4}\big\|\na|u|^2\big\|_{L^2(\mathbb{R}^2)}\\
        \leq&\ \widetilde{C}\|u\|_{L^2(\mathbb{R}^2)}\|\na u\|_{L^2(\mathbb{R}^2)}
        \|u\|_{L^4(\mathbb{R}^2)}^{2}\big\|\na|u|^2\big\|_{L^2(\mathbb{R}^2)}
        \leq \widetilde{C}\|\na u\|_{L^2(\mathbb{R}^2)}\|u\|_{L^4(\mathbb{R}^2)}^{2}
        \big\|\na|u|^2\big\|_{L^2(\mathbb{R}^2)}.
	\end{align*}
Putting the above estimate into \eqref{3.2.4-0336}, and then applying Cauchy's inequality, we deduce
    \begin{equation}
		\label{3.2.5-0336}
        4\int_{\mathbb{R}^2}P u\cdot\na |u|^2\,dx
        \leq \mu\big\|\na|u|^2\big\|_{L^2(\mathbb{R}^2)}^2
        +\widetilde{C}\|\na u\|_{L^2}^2\|u\|_{L^4(\mathbb{R}^2)}^4+\widetilde{C}.
	\end{equation}

For $n=3$, utilizing H\"{o}lder's and Cauchy's inequalities, the standard $L^p$ interpolation inequality, \eqref{3.4.4-0336}, and Lemma \ref{le2.3}, we obtain for any $r_2>\frac{3}{2}$,
    \begin{align}
		\label{3.2.4}
        4\int_{\mathbb{R}^3}P u\cdot\na |u|^2\,dx
        \leq&\ \widetilde{C}\|P\|_{L^{\frac{4r_2}{r_2+1}}(\mathbb{R}^3)}
        \|u\|_{L^{\frac{4r_2}{r_2-1}}(\mathbb{R}^3)}
        \big\|\na|u|^2\big\|_{L^2(\mathbb{R}^3)}\nonumber\\
        \leq&\ \frac{\mu}{2}\big\|\na|u|^2\big\|_{L^2(\mathbb{R}^3)}^2
        +\widetilde{C}\|P\|_{L^{\frac{4r_2}{r_2+1}}(\mathbb{R}^3)}^2
        \|u\|_{L^{\frac{4r_2}{r_2-1}}(\mathbb{R}^3)}^2\nonumber\\
        \leq&\ \frac{\mu}{2}\big\|\na|u|^2\big\|_{L^2(\mathbb{R}^3)}^2
        +\widetilde{C}\|P\|_{L^{r_2}(\mathbb{R}^3)}\|P\|_{L^{\frac{2r_2}{r_2-1}}(\mathbb{R}^3)}
        \|u\|_{L^{\frac{4r_2}{r_2-1}}(\mathbb{R}^3)}^2\nonumber\\
        \leq&\ \frac{\mu}{2}\big\|\na|u|^2\big\|_{L^2(\mathbb{R}^3)}^2
        +\widetilde{C}\|P\|_{L^{r_2}(\mathbb{R}^3)}
        \Big(\|u\|_{L^{\frac{4r_2}{r_2-1}}(\mathbb{R}^3)}^4+1\Big)\nonumber\\
        \leq&\ \frac{\mu}{2}\big\|\na|u|^2\big\|_{L^2(\mathbb{R}^3)}^2
        +\widetilde{C}\|P\|_{L^{r_2}(\mathbb{R}^3)}
        \bigg(\|u\|_{L^4(\mathbb{R}^3)}^{\frac{4r_2-6}{r_2}}
        \|u\|_{L^{12}(\mathbb{R}^3)}^{\frac{6}{r_2}}+1\bigg).
	\end{align}
It follows from the embedding $D^{1,2}(\mathbb{R}^3)\hookrightarrow L^6(\mathbb{R}^3)$ that
    \begin{equation*}
        \|u\|_{L^{12}(\mathbb{R}^3)}^{\frac{6}{r_2}}
        =\big\||u|^2\big\|_{L^6(\mathbb{R}^3)}^{\frac{3}{r_2}}
        \leq\widetilde{C}\big\|\na|u|^2\big\|_{L^2(\mathbb{R}^3)}^{\frac{3}{r_2}},\quad n=3.
	\end{equation*}
Plugging the above inequality into \eqref{3.2.4} and using Young's inequality yield
    \begin{align}
		\label{3.2.5}
        4\int_{\mathbb{R}^3}P u\cdot\na |u|^2\,dx
        \leq&\ \mu\big\|\na|u|^2\big\|_{L^2(\mathbb{R}^3)}^2
        +\widetilde{C}\|P\|_{L^{r_2}(\mathbb{R}^3)}^{\frac{2r_2}{2r_2-3}}
        \|u\|_{L^4(\mathbb{R}^3)}^4+\widetilde{C}\|P\|_{L^{r_2}(\mathbb{R}^3)}\nonumber\\
        \leq&\ \mu\big\|\na|u|^2\big\|_{L^2(\mathbb{R}^3)}^2
        +\widetilde{C}\Big(1+\|P\|_{L^{r_2}(\mathbb{R}^3)}^{s_2}\Big)
        \Big(1+\|u\|_{L^4(\mathbb{R}^3)}^4\Big),
	\end{align}
where $r_2$ and $s_2$ satisfy \eqref{th1.1.3-0406}.
Combining \eqref{3.2.5-0336} and \eqref{3.2.5}, one derives
    \begin{equation}
		\label{3.2.6}
        4\int_{\mathbb{R}^n}P u\cdot\na |u|^2\,dx
        \leq
        \begin{cases}
			\ \mu\big\|\na|u|^2\big\|_{L^2(\mathbb{R}^2)}^2
            +\widetilde{C}\|\na u\|_{L^2(\mathbb{R}^2)}^2
            \|u\|_{L^4(\mathbb{R}^2)}^4+\widetilde{C},&n=2,\smallskip\\
            \ \mu\big\|\na|u|^2\big\|_{L^2(\mathbb{R}^3)}^2
            +\widetilde{C}\Big(1+\|P\|_{L^{r_2}(\mathbb{R}^3)}^{s_2}\Big)
            \Big(1+\|u\|_{L^4(\mathbb{R}^3)}^4\Big),&n=3.
		\end{cases}
	\end{equation}

Substituting \eqref{3.2.3} and \eqref{3.2.6} into \eqref{3.2.2} leads to
    \begin{equation*}
        \frac{d}{dt}\int_{\mathbb{R}^n}|u|^4\,dx
        +\mu\int_{\mathbb{R}^n}\Big(|\na u|^2|u|^{2}+|\na|u|^2|^2\Big)\,dx
        \leq
        \begin{cases}
			\ \widetilde{C}\|\na u\|_{L^2(\mathbb{R}^2)}^2\|u\|_{L^4(\mathbb{R}^2)}^4
            +\widetilde{C},
            &n=2,\smallskip\\
            \ \widetilde{C}\Big(1+\|P\|_{L^{r_2}(\mathbb{R}^3)}^{s_2}\Big)
            \Big(1+\|u\|_{L^4(\mathbb{R}^3)}^4\Big),
            &n=3,
		\end{cases}
	\end{equation*}
which implies \eqref{3.2.1} after using Gronwall's inequality with \eqref{3.1.1} and \eqref{3.1-0406}.
The proof of Lemma \ref{le3.2} is finished.
\end{proof}

For $n=3$, according to Lemma \ref{le3.2}, the boundedness of $\mathcal{K}$ (see \eqref{3.1-0406}) yields that $\|u\|_{L^\infty(0,T;L^4(\mathbb{R}^3))}$ is bounded, which shows that the velocity $u$ satisfies Serrin's regularity condition: $u\in L^{s_1}(0,T;L^{r_1}(\mathbb{R}^3))$ with $\frac{2}{s_1}+\frac{3}{r_1}\leq1$, $3< r_1\leq \infty$.
Therefore, to prove Theorem \ref{th1.1}, we just need to prove
    \begin{equation}
        \label{3.0-0406}
        \sup\limits_{0\leq t<T^*}
        (\|u\|_{H^m(\mathbb{R}^n)}+\|\mathbb{F}-\mathbb{I}\|_{H^m(\mathbb{R}^n)})
        \leq C,
	\end{equation}
under the assumption that
    \begin{equation}
		\label{3.1}
		\mathcal{U}\triangleq
        \begin{cases}
			\ \lim\limits_{T\to T^*}
            \|\mathbb{F}\|_{L^\infty(0,T;L^\infty(\mathbb{R}^2))}<\infty,&n=2,\\
            \ \lim\limits_{T\to T^*}\left(\|\mathbb{F}\|_{L^\infty(0,T;L^\infty(\mathbb{R}^3))}
            +\|u\|_{L^{s_1}(0,T;L^{r_1}(\mathbb{R}^3))}\right)<\infty,&n=3
		\end{cases}
	\end{equation}
for any $r_1$ and $s_1$ satisfying
    \begin{equation}
		\label{th1.1.3}
        \frac{2}{s_1}+\frac{3}{r_1}\leq1,\quad
        3< r_1\leq \infty.
	\end{equation}
In \eqref{3.0-0406}, as well as in the rest of this section, the positive constant $C$ depends only on $\mathcal{U}$, $T^*$, $n$, $m$, $\mu$, and the initial data.

Let us decompose the velocity field into two parts,
i.e., $u=v+w$, where $v$ and $w$ satisfy the following two Poisson problems, respectively,
    \begin{equation}
        \label{v-eq}
        \begin{cases}
        \ \mu\De v
        =-\mathcal{P}\diver(\mathbb{F}\mathbb{F}^\top-\mathbb{I}),
        &x\in\mathbb{R}^n,\ t>0,\\
        \ v(x,t)\to0,&|x|\to\infty,\ t>0,
        \end{cases}
	\end{equation}
and
    \begin{equation*}
        \begin{cases}
        \ \mu\De w=\mathcal{P}\dot{u},
        &x\in\mathbb{R}^n,\ t>0,\\
        \ w(x,t)\to0,&|x|\to\infty,\ t>0,
        \end{cases}
	\end{equation*}
where $\dot{u}=u_t+u\cdot\na u$ is the material derivative of the velocity and $\mathcal{P}$ is the Helmholtz projection operator \eqref{n-p}.
Note that
    \begin{equation}
        \label{F-0127}
        \mathbb{F}\mathbb{F}^\top-\mathbb{I}
        =(\mathbb{F}-\mathbb{I})(\mathbb{F}-\mathbb{I})^\top
        +(\mathbb{F}-\mathbb{I})+(\mathbb{F}^\top-\mathbb{I}).
	\end{equation}
Using \eqref{2.5.1}-\eqref{2.5.3}, \eqref{3.1.1}, \eqref{3.1}, \eqref{F-0127}, and the $L^p$-boundedness of the operator $\mathcal{P}$, we have for any $k\in[2,\infty)$,
    \begin{equation}
        \label{v-es}
        \begin{cases}
        \ \|\na^2 v\|_{L^k(\mathbb{R}^n)}
        \leq C\|\na\mathbb{F}\|_{L^k(\mathbb{R}^n)},\\
        \ \|\na v\|_{L^k(\mathbb{R}^n)}
        \leq C
        \|\mathbb{F}-\mathbb{I}\|_{L^k(\mathbb{R}^n)}
        \leq C,\\
        \ \|\na v\|_{BMO(\mathbb{R}^n)}
        \leq C
        \|\mathbb{F}-\mathbb{I}\|_{L^2(\mathbb{R}^n)\cap L^\infty(\mathbb{R}^n)}
        \leq C,
        \end{cases}
	\end{equation}
and
    \begin{equation}
        \label{w-es}
        \|\na^2 w\|_{L^k(\mathbb{R}^n)}
        \leq C\|\dot{u}\|_{L^k(\mathbb{R}^n)}.
	\end{equation}

Below, we derive the crucial estimate on the $L^\infty(0,T;L^2(\mathbb{R}^n))$-norm of $\na u$.
\begin{Lemma}\label{le3.3}
Under the conditions of Theorem \ref{th1.1} and \eqref{3.1}, it holds that for any $T\in[0,T^*)$,
    \begin{equation}
		\label{3.3.1}
        \sup\limits_{0\leq t\leq T}\|\na u\|_{L^2(\mathbb{R}^n)}
        +\int_0^T\|\dot{u}\|_{L^2(\mathbb{R}^n)}^2\,dt
        \leq C.
	\end{equation}
\end{Lemma}

\begin{proof}
Taking $L^2(\mathbb{R}^n)$ inner product of $\eqref{1.1}_1$ with $u_t$, and then using integration by parts and $\diver u=0$, we have
    \begin{equation}
		\label{3.3.2}
        \frac{1}{2}\frac{d}{dt}\int_{\mathbb{R}^n}\mu|\na u|^2\,dx
        +\int_{\mathbb{R}^n}|\dot{u}|^2\,dx
        =\int_{\mathbb{R}^n}\dot{u}\cdot(u\cdot\na u)\,dx
        +\int_{\mathbb{R}^n}\diver(\mathbb{F}\mathbb{F}^\top-\mathbb{I})\cdot u_t\,dx
        \triangleq I_1+I_2.
	\end{equation}
For $I_1$, by H\"{o}lder's and Young's inequalities, \eqref{3.2.1} (for $n=2$), \eqref{v-es}, \eqref{w-es}, and Lemma \ref{le2.3}, one deduces
    \begin{align}
		\label{3.3.3-0326}
        I_1
        \leq&\ C\|\dot{u}\|_{L^2(\mathbb{R}^2)}\|u\|_{L^4(\mathbb{R}^2)}\|\na u\|_{L^4(\mathbb{R}^2)}\nonumber\\
        \leq&\ \frac{1}{8}\|\dot{u}\|_{L^2(\mathbb{R}^2)}^2+C\Big(\|\na v\|_{L^4(\mathbb{R}^2)}^2+\|\na w\|_{L^4(\mathbb{R}^2)}^2\Big)\nonumber\\
        \leq&\ \frac{1}{8}\|\dot{u}\|_{L^2(\mathbb{R}^2)}^2+C\Big(1+\|\na w\|_{L^2(\mathbb{R}^2)}\|\na^2 w\|_{L^2(\mathbb{R}^2)}\Big)\nonumber\\
        \leq&\ \frac{1}{8}\|\dot{u}\|_{L^2(\mathbb{R}^2)}^2+C\Big[1+\Big(\|\na u\|_{L^2(\mathbb{R}^2)}+\|\na v\|_{L^2(\mathbb{R}^2)}\Big)\|\dot{u}\|_{L^2(\mathbb{R}^2)}\Big]\nonumber\\
        \leq&\ \frac{1}{4}\|\dot{u}\|_{L^2(\mathbb{R}^2)}^2+C\Big(1+\|\na u\|_{L^2(\mathbb{R}^2)}^2\Big),
        \quad \text{for}\ n=2,
	\end{align}
and
    \begin{align}
		\label{3.3.3-0326-0406}
        I_1
        \leq&\ C\|\dot{u}\|_{L^2(\mathbb{R}^3)}\|u\|_{L^{r_1}(\mathbb{R}^3)}\|\na u\|_{L^{\frac{2r_1}{r_1-2}}(\mathbb{R}^3)}\nonumber\\
        \leq&\ \frac{1}{8}\|\dot{u}\|_{L^2(\mathbb{R}^3)}^2+C\|u\|_{L^{r_1}(\mathbb{R}^3)}^2
        \bigg(\|\na v\|_{L^{\frac{2r_1}{r_1-2}}(\mathbb{R}^3)}^2
        +\|\na w\|_{L^{\frac{2r_1}{r_1-2}}(\mathbb{R}^3)}^2\bigg)\nonumber\\
        \leq&\ \frac{1}{8}\|\dot{u}\|_{L^2(\mathbb{R}^3)}^2+C\|u\|_{L^{r_1}(\mathbb{R}^3)}^2
        \bigg(1+\|\na w\|_{L^2(\mathbb{R}^3)}^{\frac{2(r_1-3)}{r_1}}
        \|\na^2 w\|_{L^2(\mathbb{R}^3)}^{\frac{6}{r_1}}\bigg)\nonumber\\
        \leq&\ \frac{1}{8}\|\dot{u}\|_{L^2(\mathbb{R}^3)}^2+C\|u\|_{L^{r_1}(\mathbb{R}^3)}^2
        \bigg[1+\bigg(\|\na u\|_{L^2(\mathbb{R}^3)}^{\frac{2(r_1-3)}{r_1}}
        +\|\na v\|_{L^2(\mathbb{R}^3)}^{\frac{2(r_1-3)}{r_1}}\bigg)
        \|\dot{u}\|_{L^2(\mathbb{R}^3)}^{\frac{6}{r_1}}\bigg]\nonumber\\
        \leq&\ \frac{1}{4}\|\dot{u}\|_{L^2(\mathbb{R}^3)}
        +C\bigg(1+\|u\|_{L^{r_1}(\mathbb{R}^3)}^{\frac{2r_1}{r_1-3}}\bigg)
        \Big(1+\|\na u\|_{L^2(\mathbb{R}^3)}^2\Big)\nonumber\\
        \leq&\ \frac{1}{4}\|\dot{u}\|_{L^2(\mathbb{R}^3)}
        +C\Big(1+\|u\|_{L^{r_1}(\mathbb{R}^3)}^{s_1}\Big)
        \Big(1+\|\na u\|_{L^2(\mathbb{R}^3)}^2\Big),
        \quad \text{for}\ n=3,
	\end{align}
where $r_1$ and $s_1$ satisfy \eqref{th1.1.3}.
For $I_2$, employing integration by parts, $\mathcal{P}u=u$, and the orthogonality between $\mathcal{P}$ and $\mathcal{Q}$ (see \eqref{n-p} and \eqref{n-q}), we get
    \begin{align}
		\label{3.3.4}
        I_2
        =& -\frac{d}{dt}\int_{\mathbb{R}^n}(\mathbb{F}\mathbb{F}^\top-\mathbb{I}):\na u\,dx
        -\int_{\mathbb{R}^n}\diver(\mathbb{F}\mathbb{F}^\top-\mathbb{I})_t\cdot\mathcal{P} u\,dx\nonumber\\
        =& -\frac{d}{dt}\int_{\mathbb{R}^n}(\mathbb{F}\mathbb{F}^\top-\mathbb{I}):\na u\,dx
        -\int_{\mathbb{R}^n}\mathcal{P}\diver(\mathbb{F}\mathbb{F}^\top-\mathbb{I})_t\cdot v\,dx
        -\int_{\mathbb{R}^n}\mathcal{P}\diver(\mathbb{F}\mathbb{F}^\top-\mathbb{I})_t\cdot w\,dx\nonumber\\
        =& -\frac{1}{2}\frac{d}{dt}\int_{\mathbb{R}^n}
        \Big[2(\mathbb{F}\mathbb{F}^\top-\mathbb{I}):\na u+\mu|\na v|^2\Big]\,dx
        +\int_{\mathbb{R}^n}(\mathbb{F}\mathbb{F}^\top)_t:\na w\,dx
        +\int_{\mathbb{R}^n}
        \mathcal{Q}\diver(\mathbb{F}\mathbb{F}^\top)_t\cdot w\,dx,
	\end{align}
where in the third identity we have used
the following fact
    \begin{equation*}
        -\int_{\mathbb{R}^n}\mathcal{P}
        \diver(\mathbb{F}\mathbb{F}^\top-\mathbb{I})_t: v\,dx
        =\int_{\mathbb{R}^n}\mu\De v_t: v\,dx
        =-\frac{1}{2}\frac{d}{dt}\int_{\mathbb{R}^n}\mu|\na v|^2\,dx
	\end{equation*}
due to \eqref{v-eq}.
By the equation $\eqref{1.1}_2$, a direct calculation gives
    \begin{equation}
		\label{3.3.5}
        (\mathbb{F}\mathbb{F}^\top)_t
        =-u\cdot\na(\mathbb{F}\mathbb{F}^\top)
        +\na u(\mathbb{F}\mathbb{F}^\top)
        +(\mathbb{F}\mathbb{F}^\top)(\na u)^\top.
	\end{equation}
For the second term on the right-hand side of \eqref{3.3.4}, applying the equation \eqref{3.3.5}, integration by parts, $\diver u=0$, \eqref{3.1}, H\"{o}lder's and Cauchy's inequalities, we derive
    \begin{align}
		\label{3.3.6}
        \int_{\mathbb{R}^n}(\mathbb{F}\mathbb{F}^\top)_t:\na w\,dx
        =&\ \int_{\mathbb{R}^n}\mathbb{F}\mathbb{F}^\top:[u\cdot\na(\na w)]\,dx
        +\int_{\mathbb{R}^n}\big[\na u(\mathbb{F}\mathbb{F}^\top)
        +(\mathbb{F}\mathbb{F}^\top)(\na u)^\top\big]:\na w\,dx\nonumber\\
        \leq&\ C\Big(\|\mathbb{F}\|_{L^\infty(\mathbb{R}^n)}^2\|u\|_{L^2(\mathbb{R}^n)}\|\na^2 w\|_{L^2(\mathbb{R}^n)}
        +\|\mathbb{F}\|_{L^\infty(\mathbb{R}^n)}^2\|\na u\|_{L^2(\mathbb{R}^n)}\|\na w\|_{L^2(\mathbb{R}^n)}\Big)\nonumber\\
        \leq&\ C\big[\|\dot{u}\|_{L^2(\mathbb{R}^n)}+\|\na u\|_{L^2(\mathbb{R}^n)}(\|\na u\|_{L^2(\mathbb{R}^n)}+\|\na v\|_{L^2(\mathbb{R}^n)})\big]\nonumber\\
        \leq&\ \frac{1}{8}\|\dot{u}\|_{L^2(\mathbb{R}^n)}^2+C\Big(1+\|\na u\|_{L^2(\mathbb{R}^n)}^2\Big),
	\end{align}
where in the last inequality we have used \eqref{3.1.1}, \eqref{v-es}, and \eqref{w-es}.
For the third term on the right-hand side of \eqref{3.3.4}, we first use Plancherel's Theorem to obtain
    \begin{equation}
		\label{3.3.7-1215}
        \|\na^2 \mathcal{Q}w\|_{L^2(\mathbb{R}^n)}^2
        =\frac{1}{(2\pi)^{\frac{n}{2}}}\int_{\mathbb{R}^n}
        |\xi|^4\Big|\frac{\xi\xi^\top}{|\xi|^2}\hat{w}(\xi)\Big|^2\,d\xi
        \leq\frac{1}{(2\pi)^{\frac{n}{2}}}\int_{\mathbb{R}^n}
        |\xi|^4|\hat{w}(\xi)|^2\,d\xi
        =\|\na^2 w\|_{L^2(\mathbb{R}^n)}^2,
	\end{equation}
where $\hat{w}$ is the Fourier transform of $w$.
Utilizing the definition of the operator $\mathcal{Q}$ (see \eqref{n-q}), the estimate \eqref{3.3.7-1215}, and the $L^p$-boundedness of the Riesz transform, we get from the arguments similar to the ones used in \eqref{3.3.6} that
    \begin{align}
		\label{3.3.7}
        &\ \int_{\mathbb{R}^n}\mathcal{Q}\diver(\mathbb{F}\mathbb{F}^\top)_t\cdot w\,dx\nonumber\\
        =&\ \int_{\mathbb{R}^n}\mathcal{Q}\diver\big[
        -u\cdot\na(\mathbb{F}\mathbb{F}^\top)
        +\na u(\mathbb{F}\mathbb{F}^\top)
        +(\mathbb{F}\mathbb{F}^\top)(\na u)^\top\big]\cdot w\,dx\nonumber\\
        =& -\int_{\mathbb{R}^n}(\mathbb{F}\mathbb{F}^\top)
        :[u\cdot\na(\na\mathcal{Q}w)]\,dx
        +\int_{\mathbb{R}^n}(-\De)^{-1}\diver\diver
        \big[\na u(\mathbb{F}\mathbb{F}^\top)
        +(\mathbb{F}\mathbb{F}^\top)(\na u)^\top\big]:\diver w\,dx\nonumber\\
        \leq&\ C\Big(\|\mathbb{F}\|_{L^\infty(\mathbb{R}^n)}^2\|u\|_{L^2(\mathbb{R}^n)}\|\na^2 \mathcal{Q}w\|_{L^2(\mathbb{R}^n)}
        +\|\na u(\mathbb{F}\mathbb{F}^\top)
        +(\mathbb{F}\mathbb{F}^\top)(\na u)^\top\|_{L^2(\mathbb{R}^n)}\|\na w\|_{L^2(\mathbb{R}^n)}\Big)\nonumber\\
        \leq&\ \frac{1}{8}\|\dot{u}\|_{L^2(\mathbb{R}^n)}^2+C\Big(1+\|\na u\|_{L^2(\mathbb{R}^n)}^2\Big).
	\end{align}
Plugging \eqref{3.3.6} and \eqref{3.3.7} into \eqref{3.3.4} yields
    \begin{equation}
		\label{3.3.9}
        I_2
        \leq -\frac{1}{2}\frac{d}{dt}\int_{\mathbb{R}^n}
        \Big[2(\mathbb{F}\mathbb{F}^\top-\mathbb{I}):\na u+\mu|\na v|^2\Big]\,dx
        +\frac{1}{4}\|\dot{u}\|_{L^2(\mathbb{R}^n)}^2
        +C\Big(1+\|\na u\|_{L^2(\mathbb{R}^n)}^2\Big).
	\end{equation}
Substituting \eqref{3.3.3-0326} (for $n=2$), \eqref{3.3.3-0326-0406} (for $n=3$), and \eqref{3.3.9} into \eqref{3.3.2}, we deduce
    \begin{equation}
		\label{3.3.101}
        \Phi'(t)+\|\dot{u}\|_{L^2(\mathbb{R}^n)}^2
        \leq
        \begin{cases}
			\ C\Big(1+\|\na u\|_{L^2(\mathbb{R}^2)}^2\Big),&n=2,\smallskip\\
            \ C\Big(1+\|u\|_{L^{r_1}(\mathbb{R}^3)}^{s_1}\Big)\Big(1+\|\na u\|_{L^2(\mathbb{R}^3)}^2\Big),&n=3,
		\end{cases}
	\end{equation}
where
    \begin{equation*}
        \Phi(t)
        \triangleq \int_{\mathbb{R}^n}\Big[\mu|\na u|^2-2(\mathbb{F}\mathbb{F}^\top-\mathbb{I}):\na u+\mu|\na v|^2\Big]\,dx.
	\end{equation*}
By Cauchy's inequality and \eqref{F-0127}, one derives
    \begin{equation*}
        \left|\int_{\mathbb{R}^n}2(\mathbb{F}\mathbb{F}^\top-\mathbb{I}):\na u\,dx\right|
        \leq C(\|\mathbb{F}\|_{L^\infty(\mathbb{R}^n)}+1)
        \|\mathbb{F}-\mathbb{I}\|_{L^2(\mathbb{R}^n)}\|\na u\|_{L^2(\mathbb{R}^n)}
        \leq \frac{\mu}{2}\|\na u\|_{L^2(\mathbb{R}^n)}^2+C,
	\end{equation*}
which combined with \eqref{v-es} gives
    \begin{equation}
        \label{3.3.102}
        \frac{\mu}{2}\|\na u\|_{L^2(\mathbb{R}^n)}^2-C
        \leq \Phi(t)
        \leq \frac{3\mu}{2}\|\na u\|_{L^2(\mathbb{R}^n)}^2+C.
	\end{equation}
Integrating \eqref{3.3.101} over $[0,t]$ with $0\leq t\leq T$, and employing \eqref{3.3.102}, we have
    \begin{align*}
        \label{3.3.105-0336}
        \frac{\mu}{2}\|\na u\|_{L^2}^2
        +\int_0^t\|\dot{u}\|_{L^2}^2\,d\tau
        &\leq C+\Phi(t)
        +\int_0^t\|\dot{u}\|_{L^2}^2\,d\tau\\
        &\leq
        \begin{cases}
			\ \displaystyle
            C+C\int_0^t\Big(1+\|\na u\|_{L^2}^2\Big)\,d\tau,&n=2,\smallskip\\
            \ \displaystyle C+C\int_0^t\Big(1+\|u\|_{L^{r_1}(\mathbb{R}^3)}^{s_1}\Big)\Big(1+\|\na u\|_{L^2(\mathbb{R}^3)}^2\Big)\,d\tau,&n=3,
		\end{cases}
	\end{align*}
which yields \eqref{3.3.1} after using Gronwall's inequality and \eqref{3.1}.
We complete of the proof of Lemma \ref{le3.3}.
\end{proof}

In order to improve the regularity of the component $w$, we then establish the $L^\infty(0,T;L^2(\mathbb{R}^n))$-estimate for the material derivative of the velocity $\dot{u}$.
\begin{Lemma}\label{le3.4}
Under the conditions of Theorem \ref{th1.1} and \eqref{3.1}, it holds that for any $T\in[0,T^*)$,
    \begin{equation}
		\label{3.4.1}
        \sup\limits_{0\leq t\leq T}\|\dot{u}\|_{L^2(\mathbb{R}^n)}
        +\int_0^T\|\na\dot{u}\|_{L^2(\mathbb{R}^n)}^2\,dt
        \leq C.
	\end{equation}
\end{Lemma}

\begin{proof}
Applying the material derivative (see \eqref{n-m}) to the momentum equation $\eqref{1.1}_1$, we have
    \begin{align*}
        &\ \dot{u}_{t}+u\cdot\nabla\dot{u}+\nabla P_{t}+u\cdot\na(\na P)\\
        =&\ \mu[\De\dot{u}-\De(u\cdot\na u)+u\cdot\na(\De u)]
        +\diver(\mathbb{F}\mathbb{F}^\top)_t
        +u\cdot\na(\diver(\mathbb{F}\mathbb{F}^\top)).
	\end{align*}
Taking $L^2(\mathbb{R}^n)$ inner product of the above equation with $\dot{u}$, and using integration by parts, one gets
    \begin{align}
        \label{3.4.2}
        &\ \frac{1}{2}\frac{d}{dt}\int_{\mathbb{R}^n}|\dot{u}|^2\,dx
        +\mu\int_{\mathbb{R}^n}|\na\dot{u}|^2\,dx\nonumber\\
        =&\ -\int_{\mathbb{R}^n}[\nabla P_{t}+u\cdot\na(\na P)]\cdot\dot{u}\,dx
        +\mu\int_{\mathbb{R}^n}[-\De(u\cdot\na u)+u\cdot\na(\De u)]\cdot\dot{u}\,dx\nonumber\\
        &\ +\int_{\mathbb{R}^n}[\diver(\mathbb{F}\mathbb{F}^\top)_t
        +u\cdot\na(\diver(\mathbb{F}\mathbb{F}^\top))]\cdot\dot{u}\,dx\nonumber\\
        \triangleq&\ \sum_{i=1}^{3}J_i.
	\end{align}
For $J_1$, employing integration by parts and $\diver u=0$, we obtain
    \begin{align}
        \label{3.4.3}
        J_1
        =&\ \int_{\mathbb{R}^n}[P_{t}\diver\dot{u}
        +\na P\cdot(u\cdot\na\dot{u})]\,dx\nonumber\\
        =&\ \frac{d}{dt}\int_{\mathbb{R}^n}P\diver\dot{u}\,dx
        -\int_{\mathbb{R}^n}P[\diver\dot{u}_t+u\cdot\na (\diver\dot{u})]\,dx
        -\int_{\mathbb{R}^n}P(\na\dot{u})^\top:\na u\,dx\nonumber\\
        =&\ \frac{d}{dt}\int_{\mathbb{R}^n}P(\na u)^\top:\na u\,dx
        -3\int_{\mathbb{R}^n}P(\na\dot{u})^\top:\na u\,dx
        +2\int_{\mathbb{R}^n}P(\na u)^\top:(\na u\na u)\,dx,
	\end{align}
where we have used the following two identities
    \begin{equation}
        \label{3.4.3-0311}
        \diver\dot{u}=(\na u)^\top:\na u
	\end{equation}
and
    \begin{equation*}
        \diver\dot{u}_t+u\cdot\na (\diver\dot{u})
        =2[(\na u_t)^\top+u\cdot\na(\na u)^\top]:\na u
        =2(\na\dot{u})^\top:\na u-2(\na u)^\top:(\na u\na u).
	\end{equation*}
For the second term on the right-hand side of \eqref{3.4.3}, utilizing \eqref{3.4.4-0336}, \eqref{3.3.1}, H\"{o}lder's and Cauchy's inequalities, we derive
    \begin{align}
        \label{3.4.5}
        -3\int_{\mathbb{R}^n}P(\na\dot{u})^\top:\na u\,dx
        \leq&\ C\|P\|_{L^3(\mathbb{R}^n)}\|\na \dot{u}\|_{L^2(\mathbb{R}^n)}\|\na u\|_{L^6(\mathbb{R}^n)}\nonumber\\
        \leq&\ C\Big(1+\|u\|_{L^6(\mathbb{R}^n)}^2\Big)\|\na \dot{u}\|_{L^2(\mathbb{R}^n)}\|\na u\|_{L^6(\mathbb{R}^n)}\nonumber\\
        \leq&\ C\Big(1+\|u\|_{H^1(\mathbb{R}^n)}^2\Big)\|\na\dot{u}\|_{L^2(\mathbb{R}^n)}
        (1+\|\dot{u}\|_{L^2(\mathbb{R}^n)})\nonumber\\
        \leq&\ \frac{\mu}{6}\|\na \dot{u}\|_{L^2(\mathbb{R}^n)}^2
        +C\Big(1+\|\dot{u}\|_{L^2(\mathbb{R}^n)}^2\Big),
	\end{align}
where in the third inequality we have used Lemma \ref{le2.3} and the following estimate
    \begin{align}
		\label{3.4.5-1}
        \|\na u\|_{L^6(\mathbb{R}^n)}
        \leq&\ C(\|\na v\|_{L^6(\mathbb{R}^n)}+\|\na w\|_{L^6(\mathbb{R}^n)})\nonumber\\
        \leq&\ C(1+\|\na w\|_{L^2(\mathbb{R}^n)}+\|\na^2 w\|_{L^2(\mathbb{R}^n)})\nonumber\\
        \leq&\ C(1+\|\na u\|_{L^2(\mathbb{R}^n)}+\|\na v\|_{L^2(\mathbb{R}^n)}+\|\dot{u}\|_{L^2(\mathbb{R}^n)})\nonumber\\
        \leq&\ C(1+\|\dot{u}\|_{L^2(\mathbb{R}^n)})
	\end{align}
due to \eqref{v-es}-\eqref{3.3.1}.
For the third term on the right-hand side of \eqref{3.4.3}, using H\"{o}lder's inequality, \eqref{3.4.4-0336}, \eqref{3.3.1}, \eqref{3.4.5-1}, and Lemma \ref{le2.3}, we have
    \begin{align}
        \label{3.4.6}
        2\int_{\mathbb{R}^n}P(\na u)^\top:(\na u\na u)\,dx
        \leq&\ C\|P\|_{L^6(\mathbb{R}^n)}\|\na u\|_{L^2(\mathbb{R}^n)}\|\na u\|_{L^6(\mathbb{R}^n)}^2\nonumber\\
        \leq&\ C\Big(1+\|u\|_{L^{12}(\mathbb{R}^n)}^2\Big)\|\na u\|_{L^2(\mathbb{R}^n)}\|\na u\|_{L^6(\mathbb{R}^n)}^2\nonumber\\
        \leq&\ C\Big(1+\|u\|_{L^2(\mathbb{R}^n)}^2+\|\na u\|_{L^6(\mathbb{R}^n)}^2\Big)\|\na u\|_{L^2(\mathbb{R}^n)}\|\na u\|_{L^6(\mathbb{R}^n)}^2\nonumber\\
        \leq&\ C\Big(1+\|\dot{u}\|_{L^2(\mathbb{R}^n)}^4\Big).
	\end{align}
Putting \eqref{3.4.5} and \eqref{3.4.6} into \eqref{3.4.3} gives
    \begin{align}
        \label{3.4.8}
        J_1
        \leq \frac{d}{dt}\int_{\mathbb{R}^n}P(\na u)^\top:\na u\,dx
        +\frac{\mu}{6}\|\na \dot{u}\|_{L^2(\mathbb{R}^n)}^2
        +C\Big(1+\|\dot{u}\|_{L^2(\mathbb{R}^n)}^4\Big).
	\end{align}
For $J_2$, it follows from integration by parts, $\diver u=0$, and H\"{o}lder's inequality that
    \begin{align}
        \label{3.4.9}
        J_2
        =&\ \mu\int_{\mathbb{R}^n}[\na(u\cdot\na u):\na \dot{u}-( u\cdot\na \dot{u})\cdot\De u]\,dx\nonumber\\
        =&\ \mu\int_{\mathbb{R}^n}\Big\{[\na u\na u+u\cdot\na(\na u)]:\na \dot{u}+\na(u\cdot\na\dot{u}):\na u\Big\}\,dx\nonumber\\
        =&\ \mu\int_{\mathbb{R}^n}\Big\{(\na u\na u):\na\dot{u}+[-u\cdot\na(\na \dot{u})+\na(u\cdot\na\dot{u})]:\na u\Big\}\,dx\nonumber\\
        =&\ \mu\int_{\mathbb{R}^n}[(\na u\na u):\na\dot{u}+(\na\dot{u}\na u):\na u]\,dx\nonumber\\
        \leq&\ C\|\na \dot{u}\|_{L^2(\mathbb{R}^n)}\|\na u\|_{L^4(\mathbb{R}^n)}^2\nonumber\\
        \leq&\ \frac{\mu}{6}\|\na \dot{u}\|_{L^2(\mathbb{R}^n)}^2
        +C\Big(1+\|\dot{u}\|_{L^2(\mathbb{R}^n)}^4\Big),
	\end{align}
where in the last inequality one has used Cauchy's inequality and the following estimate
    \begin{equation*}
        \|\na u\|_{L^4(\mathbb{R}^n)}
        \leq C(\|\na u\|_{L^2(\mathbb{R}^n)}+\|\na u\|_{L^6(\mathbb{R}^n)})
        \leq C(1+\|\dot{u}\|_{L^2(\mathbb{R}^n)})
	\end{equation*}
by \eqref{3.3.1} and \eqref{3.4.5-1}.
For $J_3$, applying integration by parts, $\diver u=0$, the equation \eqref{3.3.5}, H\"{o}lder's and Cauchy's inequalities, \eqref{3.1}, and \eqref{3.3.1}, we obtain
    \begin{align}
        \label{3.4.10}
        J_3
        =&\ -\int_{\mathbb{R}^n}\big[(\mathbb{F}\mathbb{F}^\top)_{t}:\na\dot{u}
        +(u\cdot\na\dot{u})\cdot\diver(\mathbb{F}\mathbb{F}^\top)\big]\,dx\nonumber\\
        =&\ \int_{\mathbb{R}^n}\Big\{\big[u\cdot\na(\mathbb{F}\mathbb{F}^\top)
        -\na u(\mathbb{F}\mathbb{F}^\top)
        -(\mathbb{F}\mathbb{F}^\top)(\na u)^\top\big]:\na\dot{u}
        +\na(u\cdot\na\dot{u}):(\mathbb{F}\mathbb{F}^\top)\Big\}\,dx\nonumber\\
        =&\ \int_{\mathbb{R}^n}\big[-\na u(\mathbb{F}\mathbb{F}^\top)
        -(\mathbb{F}\mathbb{F}^\top)(\na u)^\top\big]:\na\dot{u}\,dx
        +\int_{\mathbb{R}^n}[-u\cdot\na(\na \dot{u})+\na(u\cdot\na\dot{u})]:(\mathbb{F}\mathbb{F}^\top)\,dx\nonumber\\
        =&\ \int_{\mathbb{R}^n}\Big[\big(-\na u(\mathbb{F}\mathbb{F}^\top)
        -(\mathbb{F}\mathbb{F}^\top)(\na u)^\top\big):\na\dot{u}+(\na\dot{u}\na u):(\mathbb{F}\mathbb{F}^\top)\Big]\,dx\nonumber\\
        \leq&\ C\|\mathbb{F}\|_{L^\infty(\mathbb{R}^n)}^2\|\na u\|_{L^2(\mathbb{R}^n)}\|\na \dot{u}\|_{L^2(\mathbb{R}^n)}\nonumber\\
        \leq&\ \frac{\mu}{6}\|\na \dot{u}\|_{L^2(\mathbb{R}^n)}^2+C.
	\end{align}
Substituting \eqref{3.4.8}-\eqref{3.4.10} into \eqref{3.4.2} gives
    \begin{equation}
        \label{3.4.13}
        \Psi'(t)
        +\mu\int_{\mathbb{R}^n}|\na\dot{u}|^2\,dx
        \leq C\Big(1+\|\dot{u}\|_{L^2(\mathbb{R}^n)}^4\Big),
	\end{equation}
where
    \begin{equation*}
        \Psi(t)
        \triangleq \int_{\mathbb{R}^n}\Big(|\dot{u}|^2-2P(\na u)^\top:\na u\Big)\,dx.
	\end{equation*}
By the arguments similar to the ones used in \eqref{3.4.5}, we derive
    \begin{align*}
        \left|-2\int_{\mathbb{R}^n}P(\na u)^\top:\na u\,dx\right|
        \leq&\ C\|P\|_{L^3(\mathbb{R}^n)}\|\na u\|_{L^2(\mathbb{R}^n)}\|\na u\|_{L^6(\mathbb{R}^n)}\nonumber\\
        \leq&\ C\Big(1+\|u\|_{H^1(\mathbb{R}^n)}^2\Big)\|\na u\|_{L^2(\mathbb{R}^n)}(1+\|\dot{u}\|_{L^2(\mathbb{R}^n)})\nonumber\\
        \leq&\ \frac{1}{2}\|\dot{u}\|_{L^2(\mathbb{R}^n)}^2+C,
	\end{align*}
which implies
    \begin{equation}
        \label{3.4.17}
        \frac{1}{2}\|\dot{u}\|_{L^2(\mathbb{R}^n)}^2-C
        \leq \Psi(t)
        \leq \frac{3}{2}\|\dot{u}\|_{L^2(\mathbb{R}^n)}^2+C.
	\end{equation}
Integrating \eqref{3.4.13} over $[0,t]$ with $0\leq t\leq T$, and utilizing \eqref{3.4.17}, one deduces
    \begin{equation*}
        \frac{1}{2}\|\dot{u}\|_{L^2(\mathbb{R}^n)}^2
        +\mu\int_0^t\|\na\dot{u}\|_{L^2(\mathbb{R}^n)}^2\,d\tau
        \leq C+\Psi(t)+\mu\int_0^t\|\na\dot{u}\|_{L^2(\mathbb{R}^n)}^2\,d\tau
        \leq C+C\int_0^t\|\dot{u}\|_{L^2(\mathbb{R}^n)}^4\,d\tau,
	\end{equation*}
which yields \eqref{3.4.1} after using Gronwall's inequality and \eqref{3.3.1}.
This completes the proof of Lemma \ref{le3.4}.
\end{proof}

As a corollary of Lemma \ref{le3.4}, we obtain the higher-order regularity of the quantity $w$ as follows.
\begin{Corollary}
	\label{co3.1}
Under the conditions of Theorem \ref{th1.1} and \eqref{3.1}, it holds that for any $T\in[0,T^*)$,
    \begin{equation}
        \label{co3.1.1}
        \sup\limits_{0\leq t\leq T}\|\na^2 w\|_{L^2(\mathbb{R}^n)}
        +\int_0^T\Big(\|\na w\|_{L^\infty(\mathbb{R}^n)}^2+\|\na^2 w\|_{L^4(\mathbb{R}^n)}^2\Big)\,dt
        \leq C.
	\end{equation}
\end{Corollary}

\begin{proof}
It follows from \eqref{v-es}-\eqref{3.3.1}, \eqref{3.4.1} and Lemma \ref{le2.3} that
    \begin{equation*}
        \sup\limits_{0\leq t\leq T}\|\na^2 w\|_{L^2(\mathbb{R}^n)}
        \leq C\sup\limits_{0\leq t\leq T}\|\dot{u}\|_{L^2(\mathbb{R}^n)}
        \leq C,
	\end{equation*}
and
    \begin{align*}
        \int_0^T\Big(\|\na w\|_{L^\infty(\mathbb{R}^n)}^2+\|\na^2 w\|_{L^4(\mathbb{R}^n)}^2\Big)\,dt
        \leq&\ C\int_0^T\Big(\|\na w\|_{L^2(\mathbb{R}^n)}^2+\|\na^2 w\|_{L^4(\mathbb{R}^n)}^2\Big)\,dt\\
        \leq&\ C\int_0^T\Big(\|\na u\|_{L^2(\mathbb{R}^n)}^2+\|\na v\|_{L^2(\mathbb{R}^n)}^2+\|\dot{u}\|_{L^4(\mathbb{R}^n)}^2\Big)\,dt\\
        \leq&\ C\int_0^T\Big(1+\|\dot{u}\|_{H^1(\mathbb{R}^n)}^2\Big)\,dt\\
        \leq&\ C.
	\end{align*}
The proof of Corollary \ref{co3.1} is finished.
\end{proof}

In the following lemma, we give the estimates on both the $L^\infty(0,T;L^4(\mathbb{R}^n))$-norm of $\na \mathbb{F}$ and the $L^1(0,T;L^\infty(\mathbb{R}^n))$-norm of $\na u$, which are key to establishing high-order energy estimates for the solutions.
\begin{Lemma}\label{le3.5}
Under the conditions of Theorem \ref{th1.1} and \eqref{3.1}, it holds that for any $T\in[0,T^*)$,
    \begin{equation}
        \label{3.5.1}
        \sup\limits_{0\leq t\leq T}\|\na \mathbb{F}\|_{L^2(\mathbb{R}^n)\cap L^4(\mathbb{R}^n)}
        +\int_0^T\Big(\|\na u\|_{L^\infty(\mathbb{R}^n)}^2+\|\na^2 u\|_{L^4(\mathbb{R}^n)}^2\Big)\,dt
        \leq C.
	\end{equation}
\end{Lemma}

\begin{proof}
For $2\leq k\leq 4$, applying $\na$ to $\eqref{1.1}_2$, multiplying the resulting equation by $k|\na\mathbb{F}|^{k-2}\na\mathbb{F}$, and then integrating it over $\mathbb{R}^n$ by parts, by H\"{o}lder's inequality, \eqref{3.1}, and \eqref{v-es}, we obtain
    \begin{align}
        \label{3.5.2}
        \frac{d}{dt}\|\na\mathbb{F}\|_{L^k(\mathbb{R}^n)}^k
        =&\ \int_{\mathbb{R}^n}\pa_i(-u\cdot\na\mathbb{F}+\na u\mathbb{F})
        :(k|\na\mathbb{F}|^{k-2}\pa_i\mathbb{F})\,dx\nonumber\\
        =&\ \int_{\mathbb{R}^n} [-(\pa_i u)\cdot\na\mathbb{F}
        +\na(\pa_i u)\mathbb{F}
        +\na u(\pa_i \mathbb{F})]
        :(k|\na\mathbb{F}|^{k-2}\pa_i\mathbb{F})\,dx\nonumber\\
        \leq&\ C\Big(\|\na u\|_{L^\infty(\mathbb{R}^n)}\|\na\mathbb{F}\|_{L^k(\mathbb{R}^n)}^k
        +\|\na^2 u\|_{L^k(\mathbb{R}^n)}\|\na\mathbb{F}\|_{L^k(\mathbb{R}^n)}^{k-1}\Big)\nonumber\\
        \leq&\ C(\|\na v\|_{L^\infty(\mathbb{R}^n)}+\|\na w\|_{L^\infty(\mathbb{R}^n)})\|\na\mathbb{F}\|_{L^k(\mathbb{R}^n)}^k
        +C\Big(\|\na^2 v\|_{L^k(\mathbb{R}^n)}+\|\na^2 w\|_{L^k(\mathbb{R}^n)}\Big)\|\na\mathbb{F}\|_{L^k(\mathbb{R}^n)}^{k-1}\nonumber\\
        \leq&\ C(1+\|\na v\|_{L^\infty(\mathbb{R}^n)}+\|\na w\|_{L^\infty(\mathbb{R}^n)})\|\na\mathbb{F}\|_{L^k(\mathbb{R}^n)}^k
        +C\|\na^2 w\|_{L^k(\mathbb{R}^n)}\|\na\mathbb{F}\|_{L^k(\mathbb{R}^n)}^{k-1},
	\end{align}
where in the second identity we have used the following fact
    \begin{equation*}
        -\int_{\mathbb{R}^n}u\cdot\na(\pa_i\mathbb{F}): (p|\na\mathbb{F}|^{p-2}\pa_i\mathbb{F})\,dx
        =\int_{\mathbb{R}^n}\diver u|\na\mathbb{F}|^p\,dx
        =0
	\end{equation*}
due to $\diver u=0$.
It follows from \eqref{3.5.2} that
    \begin{equation}
        \label{3.5.3}
        \frac{d}{dt}\ln(e+\|\na \mathbb{F}\|_{L^k(\mathbb{R}^n)})
        \leq C(1+\|\na v\|_{L^\infty(\mathbb{R}^n)}+\|\na w\|_{L^\infty(\mathbb{R}^n)}+\|\na^2 w\|_{L^k(\mathbb{R}^n)}).
	\end{equation}
Moreover, by \eqref{2.6.1} and \eqref{v-es}, we have
    \begin{equation}
        \label{3.5.4}
        \|\na v\|_{L^\infty(\mathbb{R}^n)}
        \leq C\big[1+\|\na v\|_{BMO(\mathbb{R}^n)}\ln(e+\|\na v\|_{W^{1,4}(\mathbb{R}^n)})\big]
        \leq C\big[1+\ln(e+\|\na \mathbb{F}\|_{L^4(\mathbb{R}^n)})\big].
	\end{equation}
Taking $p=4$ in \eqref{3.5.3} and using \eqref{3.5.4}, one deduces
    \begin{equation*}
        \frac{d}{dt}\ln(e+\|\na \mathbb{F}\|_{L^4(\mathbb{R}^n)})
        \leq C\big[1+\|\na w\|_{L^\infty(\mathbb{R}^n)}+\|\na^2 w\|_{L^4(\mathbb{R}^n)}+
        \ln(e+\|\na \mathbb{F}\|_{L^4(\mathbb{R}^n)})\big],
	\end{equation*}
which combined with Gronwall's inequality and \eqref{co3.1.1} yields
    \begin{equation}
        \label{3.5.5}
        \sup\limits_{0\leq t\leq T}\|\na \mathbb{F}\|_{L^4(\mathbb{R}^n)}
        \leq C.
	\end{equation}
It follows from \eqref{v-es}, \eqref{3.5.5} and Lemma \ref{le2.3} that
    \begin{equation}
        \label{3.5.6}
        \sup\limits_{0\leq t\leq T}\Big(\|\na v\|_{L^\infty(\mathbb{R}^n)}+\|\na^2 v\|_{L^4(\mathbb{R}^n)}\Big)
        \leq C\sup\limits_{0\leq t\leq T}\Big(\|\na v\|_{L^2(\mathbb{R}^n)}+\|\na^2 v\|_{L^4(\mathbb{R}^n)}\Big)
        \leq C,
	\end{equation}
which together with \eqref{co3.1.1} gives
    \begin{equation*}
        \int_0^T\Big(\|\na u\|_{L^\infty(\mathbb{R}^n)}^2+\|\na^2 u\|_{L^4(\mathbb{R}^n)}^2\Big)\,dt
        \leq C.
	\end{equation*}

In addition, taking $p=2$ in \eqref{3.5.3} and using \eqref{3.5.6}, one gets
    \begin{equation*}
        \frac{d}{dt}\ln(e+\|\na \mathbb{F}\|_{L^2(\mathbb{R}^n)})
        \leq C(1+\|\na w\|_{L^\infty(\mathbb{R}^n)}+\|\na^2 w\|_{L^2(\mathbb{R}^n)}),
	\end{equation*}
which, along with Gronwall's inequality and \eqref{co3.1.1}, implies
    \begin{equation*}
        \sup\limits_{0\leq t\leq T}\|\na\mathbb{F}\|_{L^2(\mathbb{R}^n)}
        \leq C.
	\end{equation*}
The proof of Lemma \ref{le3.5} is completed.
\end{proof}

Finally, we derive the following high-order regularity estimates for $(u,\mathbb{F})$.
\begin{Lemma}\label{le3.6}
Under the conditions of Theorem \ref{th1.1} and \eqref{3.1}, it holds that for any $T\in[0,T^*)$ and $2\leq l\leq m$,
    \begin{equation}
        \label{3.6.1}
        \sup\limits_{0\leq t\leq T}\Big(\|\na^l u\|_{L^2(\mathbb{R}^n)}+\|\na^l \mathbb{F}\|_{L^2(\mathbb{R}^n)}\Big)
        +\int_0^T\|\na^{l+1} u\|_{L^2(\mathbb{R}^n)}^2\,dt
        \leq C.
	\end{equation}
\end{Lemma}

\begin{proof}
For $2\leq l\leq m$, applying $\na^l$ to the first two equations of \eqref{1.1}, taking $L^2(\mathbb{R}^n)$ inner product of the resulting equations with $\na^l u$ and $\na^l \mathbb{F}$, respectively, and then summing up, by integration by parts and $\diver u=0$, we get
	\begin{align}
		\label{3.6.2}
		&\ \frac{1}{2}\frac{d}{dt}\int_{\mathbb{R}^n}\Big(|\na^l u|^2+|\na^l\mathbb{F}|^2\Big)\,dx
        +\mu\int_{\mathbb{R}^n}|\na^{l+1} u|^2\,dx\nonumber\\
		=&\ -\int_{\mathbb{R}^n}\Big[\na^l(u\cdot\nabla u)-u\cdot\nabla(\na^l u)\Big]: \na^l u\,dx
        -\int_{\mathbb{R}^n}\na^l(\mathbb{F}\mathbb{F}^{\top}): \na^l(\nabla u)\,dx\nonumber\\
		&\ -\int_{\mathbb{R}^n}\Big[\na^l(u\cdot\nabla \mathbb{F})-u\cdot\nabla(\na^l \mathbb{F})\Big]: \na^l \mathbb{F}\,dx
        +\int_{\mathbb{R}^n}\na^l(\nabla u\mathbb{F}): \na^l\mathbb{F}\,dx\nonumber\\
        \triangleq& \sum_{i=1}^{4}W_i.
	\end{align}
Using H\"{o}lder's and Cauchy's inequalities, Lemmas \ref{le2.3} and \ref{le2.4}, \eqref{3.1}, \eqref{3.3.1}, and \eqref{3.5.1}, we can estimate $W_1$-$W_4$ as follows:
    \begin{align*}
		W_1
        \leq&\ C\|\na^l(u\cdot\nabla u)-u\cdot\nabla(\na^l u)\|_{L^2(\mathbb{R}^n)}\|\na^l u\|_{L^2(\mathbb{R}^n)}
        \leq C\|\nabla u\|_{L^\infty(\mathbb{R}^n)}\|\na^l u\|_{L^2(\mathbb{R}^n)}^2;\\
        W_2
        \leq&\ C\|\na^l(\mathbb{F}\mathbb{F}^{\top})\|_{L^2(\mathbb{R}^n)}\|\na^{l+1} u\|_{L^2(\mathbb{R}^n)}
		\leq C\|\mathbb{F}\|_{L^\infty(\mathbb{R}^n)}
        \|\na^l\mathbb{F}\|_{L^2(\mathbb{R}^n)}\|\na^{l+1} u\|_{L^2(\mathbb{R}^n)}\\
        \leq&\ \frac{\mu}{6}\|\na^l\nabla u\|_{L^2(\mathbb{R}^n)}^2+C\|\na^l \mathbb{F}\|_{L^2(\mathbb{R}^n)}^2;\\
        W_3
		\leq&\ C\|\na^l \mathbb{F}\|_{L^2(\mathbb{R}^n)}\|\na^l(u\cdot\nabla \mathbb{F})-u\cdot\nabla(\na^l \mathbb{F})\|_{L^2(\mathbb{R}^n)}\\
		\leq&\ C\|\na^l \mathbb{F}\|_{L^2(\mathbb{R}^n)}\Big(\|\nabla u\|_{L^\infty(\mathbb{R}^n)}\|\na^l \mathbb{F}\|_{L^2(\mathbb{R}^n)}+\|\nabla \mathbb{F}\|_{L^4(\mathbb{R}^n)}\|\na^l u\|_{L^4(\mathbb{R}^n)}\Big)\\
        \leq&\ C\|\na^l \mathbb{F}\|_{L^2(\mathbb{R}^n)}\Big(\|\nabla u\|_{L^\infty(\mathbb{R}^n)}\|\na^l \mathbb{F}\|_{L^2(\mathbb{R}^n)}+\|\na u\|_{L^2(\mathbb{R}^n)}+\|\na^{l+1} u\|_{L^2(\mathbb{R}^n)}\Big)\\
        \leq&\ \frac{\mu}{6}\|\na^{l+1} u\|_{L^2(\mathbb{R}^n)}^2
        +C(1+\|\nabla u\|_{L^\infty(\mathbb{R}^n)})\|\na^l \mathbb{F}\|_{L^2(\mathbb{R}^n)}^2+C;\\
        W_4
        \leq&\ C\|\na^l(\nabla u\mathbb{F})\|_{L^2(\mathbb{R}^n)}\|\na^l \mathbb{F}\|_{L^2(\mathbb{R}^n)}\\
        \leq&\ C\Big(\|\na^{l+1} u\|_{L^2(\mathbb{R}^n)}\|\mathbb{F}\|_{L^\infty(\mathbb{R}^n)}
        +\|\na u\|_{L^\infty(\mathbb{R}^n)}\|\na^l \mathbb{F}\|_{L^2(\mathbb{R}^n)}\Big)
        \|\na^l \mathbb{F}\|_{L^2(\mathbb{R}^n)}\\
        \leq&\ \frac{\mu}{6}\|\na^{l+1} u\|_{L^2(\mathbb{R}^n)}^2
        +C(1+\|\nabla u\|_{L^\infty(\mathbb{R}^n)})
        \|\na^l \mathbb{F}\|_{L^2(\mathbb{R}^n)}^2.
	\end{align*}
Substituting the above estimates into \eqref{3.6.2} yields
    \begin{equation*}
		\frac{d}{dt}\Big(\|\na^l u\|_{L^2(\mathbb{R}^n)}^2+\|\na^l \mathbb{F}\|_{L^2(\mathbb{R}^n)}^2\Big)
        +\mu\|\na^{l+1} u\|_{L^2(\mathbb{R}^n)}^2
		\leq C(1+\|\nabla u\|_{L^\infty(\mathbb{R}^n)})
        \Big(\|\na^l u\|_{L^2(\mathbb{R}^n)}^2+\|\na^l \mathbb{F}\|_{L^2(\mathbb{R}^n)}^2\Big)+C,
	\end{equation*}
which combined with Gronwall's inequality and \eqref{3.5.1} leads to \eqref{3.6.1}.
The proof of Lemma \ref{le3.6} is finished.
\end{proof}

Combining all the estimates obtained in Lemmas \ref{le3.1}-\ref{le3.6} with the local existence result in Proposition \ref{prop1.1}, we are now ready to prove Theorem \ref{th1.1}.
\begin{proof}[\textbf{Proof of Theorem \ref{th1.1}}]
From the identities \eqref{2.2.1} in Lemma \ref{le2.2} and the a priori estimates \eqref{3.1.1}, \eqref{3.2.1}, \eqref{3.3.1}, \eqref{3.4.1}, \eqref{3.5.1}, \eqref{3.6.1} in Lemmas \ref{le3.1}-\ref{le3.6}, it follows that the function
    \begin{equation*}
        (\tilde{u},\widetilde{\mathbb{F}})(x)
        =(u,\mathbb{F})(x,T^*)
        \triangleq\lim\limits_{t\to T^*}(u,\mathbb{F})(x,t)
	\end{equation*}
satisfies $(\tilde{u},\widetilde{\mathbb{F}}-\mathbb{I})\in H^m(\mathbb{R}^n)$ and
    \begin{equation*}
		\begin{cases}
            \ \diver \tilde{u}=0,\\
            \ \det\widetilde{\mathbb{F}}=1,\\
			\ \widetilde{\mathbb{F}}^{lk}\na_l\widetilde{\mathbb{F}}^{ij}
            =\widetilde{\mathbb{F}}^{lj}\na_l\widetilde{\mathbb{F}}^{ik},\quad
            1\leq i,j,k\leq n.
		\end{cases}
	\end{equation*}
Therefore, one can take $(\tilde{u},\widetilde{\mathbb{F}})$ as the initial data and utilize Proposition \ref{prop1.1} to extend the local classical solution beyond $T^*$.
This contradicts the assumption that $T^*$ is the maximum existence time of the solution.
Hence, we complete the proof of Theorem \ref{th1.1}.
\end{proof}

\section{Global Well-posedness}\label{sec4}
In this section, we investigate the global-in-time existence of classical solutions to the Cauchy problem \eqref{1.1}-\eqref{1.2}.
For some $T>0$, we assume in the rest of this section that $(u,\mathbb{F})$ is the classical solution to \eqref{1.1}-\eqref{1.2} in $\mathbb{R}^n\times[0,T]$.
Let $\si=\si(t)\triangleq\min\{1,t\}$ and define the following energy functionals:
    \begin{equation}\label{A1}
    A_1(T)\triangleq
    \sup\limits_{0\leq t\leq T}\int_{\mathbb{R}^{n}}\Big(\si|\nabla u|^2
    +\si^n|\mathcal{P}\dot{u}|^2\Big)\,dx
    +\int_0^T\int_{\mathbb{R}^n}\Big(\si|\mathcal{P}\dot{u}|^2
    +\si^n|\nabla\mathcal{P}\dot{u}|^2\Big)\,dxdt,
    \end{equation}
    \begin{equation}\label{A2}
    A_2(T)\triangleq
    \sup\limits_{0\leq t\leq T}\int_{\mathbb{R}^{n}}|u|^q\,dx
        +\int_0^{T}\int_{\mathbb{R}^n}\Big(|\na u|^2|u|^{q-2}
        +\big|\na|u|^{\frac{q}{2}}\big|^2\Big)\,dxdt,
    \end{equation}
and
    \begin{equation}\label{B}
    B(T)\triangleq\sup\limits_{0\leq t\leq T}
    \|\mathbb{F}-\mathbb I\|_{\dot B_{p,1}^{\frac{n}{p}}(\mathbb{R}^n)}^2,
    \end{equation}
where $\eta$ and $p$ satisfy \eqref{th1.2.0}-\eqref{th1.2.0-0413}, and
    \begin{equation}
        \label{q-0404}
        q\triangleq
        \begin{cases}
            \ \frac{4}{2-\eta},&\text{if}\ n=2,\\
            \ 3,&\text{if}\ n=3.
		\end{cases}
    \end{equation}

We have the following key a priori estimates on the solution $(u,\mathbb{F})$.
\begin{Proposition}\label{prop4.1}
Under the conditions of Theorem \ref{th1.2}, there exists a positive constant $\de$ depending only on $n$, $m$, $M$, $\eta$, $p$, and $\mu$, such that if $(u,\mathbb{F})$ is a classical solution of the Cauchy problem \eqref{1.1}-\eqref{1.2} in $\mathbb{R}^n\times[0,T]$ satisfying $\mathcal{E}_0\leq \de$ and
    \begin{equation}
        \label{prop4.1.1}
        A_1(T)\leq 2\mathcal{E}_0^{\frac{1}{4}},\quad
        A_2(\sigma(T))\leq 2\mathcal{E}_0^{\al},\quad
        B(T)\leq 2\mathcal{E}_0^{\be},
    \end{equation}
then the following estimates hold
    \begin{equation*}
        A_1(T)\leq \mathcal{E}_0^{\frac{1}{4}},\quad
        A_2(\sigma(T))\leq \mathcal{E}_0^{\al},\quad
        B(T)\leq \mathcal{E}_0^{\be},
    \end{equation*}
where
    \begin{align*}
        &n=2:\al=\frac{1}{2(2-\eta)}\in\Big(\frac{1}{4},\frac{1}{2}\Big],\
        \be=\frac{1}{16};\\
        &n=3:\al=\frac{3(2\eta-1)}{8\eta}\in\Big(\frac{3}{64},\frac{3}{8}\Big],\
        \be=\frac{(15\eta-8)}{192\eta}\in\Big(0,\frac{7}{192}\Big].
    \end{align*}

\end{Proposition}

\begin{proof}
Proposition \ref{prop4.1} follows from Lemmas \ref{le4.2}, \ref{le4.3}, and \ref{le4.5} below.
\end{proof}

Throughout this section, we assume $\mathcal{E}_0\leq1$ for simplicity.
Moreover, we denote by $C$ or $C_i$ $(1\leq i\leq 7)$ the generic positive constant which may depend on $n$, $m$, $M$, $\eta$, $p$, and $\mu$, but is independent of time $T$.
In addition, we use $\de_i$ $(1\leq i\leq 4)$ to denote a small positive constant, which may depend on $n$, $m$, $M$, $\eta$, $p$, and $\mu$, but is independent of time $T$.

To begin with, we derive the elementary energy estimate for $(u,\mathbb F)$ and the preliminary bound for $A_1(T)$ in the following lemma.
\begin{Lemma}\label{le4.1}
Under the conditions of Proposition \ref{prop4.1}, it holds that
    \begin{equation}
        \label{4.1.1}
        \sup\limits_{0\leq t\leq T}\int_{\mathbb{R}^{n}}\Big(|u|^2+|\mathbb F-\mathbb I|^2\Big)\,dx
        +2\mu\int_0^T\int_{\mathbb{R}^{n}}|\nabla u|^2\,dxdt\leq C\mathcal{E}_0,
    \end{equation}
and
    \begin{equation}
        \label{4.1.2}
        A_1(T)\leq C\left(\mathcal{E}_0
        +\left|\int_0^T\int_{\mathbb{R}^{n}}\sigma \na u:(\na u\na u)\,dxdt\right|
        +\int_0^T\int_{\mathbb{R}^{n}}\si^n|\na u|^4\,dxdt\right).
    \end{equation}
\end{Lemma}

\begin{proof}
The elementary energy estimate \eqref{4.1.1} follows from \eqref{3.1.1} directly.
Let
    \begin{equation*}
    A_{11}(T)\triangleq
    \sup\limits_{0\leq t\leq T}\int_{\mathbb{R}^{n}}\si|\nabla u|^2\,dx
    +\int_0^T\int_{\mathbb{R}^n}\si|\mathcal{P}\dot{u}|^2\,dxdt
    \end{equation*}
and
    \begin{equation*}
    A_{12}(T)\triangleq
    \sup\limits_{0\leq t\leq T}\int_{\mathbb{R}^{n}}\si^n|\mathcal{P}\dot{u}|^2\,dx
    +\int_0^T\int_{\mathbb{R}^n}\si^n|\nabla\mathcal{P}\dot{u}|^2\,dxdt.
    \end{equation*}


We first give an estimate of $A_{11}(T)$.
Applying the Helmholtz projection operator $\mathcal{P}$ to $\eqref{1.1}_1$ gives
    \begin{equation}
        \label{4.1.3}
        \mathcal{P}\dot{u}
        =\mu\Delta u+\mathcal{P}\operatorname{div}(\mathbb F\mathbb F^{\top}).
    \end{equation}
Multiplying \eqref{4.1.3} by $\sigma\mathcal{P}\dot{u}$, and integrating the resulting equation over $\mathbb{R}^n\times [0, t]$ with $0\leq t\leq T$, we have
    \begin{align}
        \label{4.1.4}
        &\ \int_0^t\int_{\mathbb{R}^{n}}\sigma |\mathcal{P}\dot {u}|^2\,dxd\tau\nonumber\\
        =&\ \mu\int_0^t\int_{\mathbb{R}^{n}}\Delta u\cdot \sigma\mathcal{P} \dot{u}\,dxd\tau
        +\int_0^t\int_{\mathbb{R}^{n}}\mathcal{P}\operatorname{div}(\mathbb F\mathbb F^{\top})\cdot \sigma\mathcal{P} \dot{u}\,dxd\tau\nonumber\\
        =&\ \mu\int_0^t\int_{\mathbb{R}^{n}}\sigma \Delta u\cdot[u_{\tau}
        +(I_d-\mathcal{Q})(u\cdot \nabla u)]\,dxd\tau
        +\int_0^t\int_{\mathbb{R}^{n}}\sigma\mathcal{P}\operatorname{div}(\mathbb F\mathbb F^{\top})\cdot (I_d-\mathcal{Q}) \dot{u}\,dxd\tau\nonumber\\
        =&\ \mu\int_0^t\int_{\mathbb{R}^{n}}\sigma \Delta u\cdot u_{\tau}\,dxd\tau
        +\mu\int_0^t\int_{\mathbb{R}^{n}}\sigma \Delta u\cdot(u\cdot\nabla u)\,dxd\tau\nonumber\\
        &\ +\int_0^t\int_{\mathbb{R}^{n}}\sigma \operatorname{div}(\mathbb F\mathbb F^{\top})\cdot \dot{u}\,dxd\tau
        -\int_0^t\int_{\mathbb{R}^{n}}\sigma\mathcal{Q} \operatorname{div}(\mathbb F\mathbb F^{\top})\cdot\dot u\,dxd\tau\nonumber\\
        \triangleq&\ \sum_{i=1}^4 K_i,
    \end{align}
where in the third inequality we have used the fact
    \begin{equation}
        \label{4.1.4-0319}
        \mu\int_0^t\int_{\mathbb{R}^{n}}\sigma \Delta u\cdot\mathcal{Q}(u\cdot\nabla u)\,dxd\tau
        =\mu\int_0^t\int_{\mathbb{R}^{n}}\sigma \mathcal{P}\Delta u\cdot\mathcal{Q}(u\cdot\nabla u)\,dxd\tau
        =0
    \end{equation}
due to
the orthogonality between $\mathcal{P}$ and $\mathcal{Q}$ (see \eqref{n-p} and \eqref{n-q}).
For $K_1$, employing integration by parts and \eqref{4.1.1}, one gets
    \begin{align}
        \label{4.1.5}
        K_1
        & =-\mu\int_0^t
        \Big(\frac{d}{d\tau}\int_{\mathbb{R}^{n}}\frac{\sigma}{2}|\nabla u|^2\,dx
        -\frac{1}{2}\int_{\mathbb{R}^{n}}\sigma'|\nabla u|^2\,dx\Big)\,d\tau\nonumber\\
        & =-\frac{\mu}{2}\int_{\mathbb{R}^{n}}\sigma|\nabla u|^2\,dx+\frac{1}{2}\int_0^{\sigma}\int_{\mathbb{R}^{n}}|\nabla u|^2\,dxd\tau\nonumber\\
        & \leq-\frac{\mu}{2}\int_{\mathbb{R}^{n}}\sigma|\nabla u|^2\,dx+C\mathcal{E}_0.
    \end{align}
For $K_2$, it follows from integration by parts that
    \begin{align}
        \label{4.1.6}
        K_2
        =&-\mu\int_0^t\int_{\mathbb{R}^{n}}\sigma \na u:\na(u\cdot\nabla u)\,dxd\tau
        =-\mu\int_0^t\int_{\mathbb{R}^{n}}\sigma \na u:(\na u\na u)\,dxd\tau\nonumber\\
        \leq&\ C\left|\int_0^t\int_{\mathbb{R}^{n}}\sigma \na u:(\na u\na u)\,dxd\tau\right|,
    \end{align}
where one has used the identity
    \begin{equation*}
        -\mu\int_0^t\int_{\mathbb{R}^{n}}\si \na u:[u\cdot\na (\na u)]\,dxd\tau
        =\frac{\mu}{2}\int_0^t\int_{\mathbb{R}^{n}}\si \diver u |\na u|^2\,dxd\tau
        =0
    \end{equation*}
due to $\diver u=0$.
For $K_3$, we first deduce from $\diver\mathbb{F}^{\top}=0$ that
    \begin{equation}
        \label{4.1.7-0311}
        \diver(\mathbb{F}\mathbb{F}^{\top})
        =\diver((\mathbb{F}-\mathbb I)(\mathbb{F}-\mathbb I)^\top)+\diver(\mathbb{F}-\mathbb I).
    \end{equation}
By the above identity and integration by parts, one obtains
    \begin{align}
        \label{4.1.7}
        K_3
        =& \int_0^t\int_{\mathbb{R}^{n}}\sigma (\operatorname{div}(\mathbb F-\mathbb I)+\operatorname{div}((\mathbb F-\mathbb I)(\mathbb F-\mathbb I)^\top))\cdot \dot{u}\,dxd\tau\nonumber\\
        =& -\int_0^t\int_{\mathbb{R}^{n}}\sigma (\mathbb F-\mathbb I):\nabla u_{\tau}\,dxd\tau
        -\int_0^t\int_{\mathbb{R}^{n}}\sigma (\mathbb F-\mathbb I):\nabla(u\cdot \nabla u)\,dxd\tau\nonumber\\
        & -\int_0^t\int_{\mathbb{R}^{n}}\sigma ((\mathbb F-\mathbb I)(\mathbb F-\mathbb I)^\top):\nabla u_{\tau}\,dxd\tau
        -\int_0^t\int_{\mathbb{R}^{n}}\sigma ((\mathbb F-\mathbb I)(\mathbb F-\mathbb I)^\top):\nabla(u\cdot \nabla u)\,dxd\tau\nonumber\\
        =& -\int_{\mathbb{R}^{n}}\sigma (\mathbb F-\mathbb I):\nabla u\,dx
        +\int_0^t\int_{\mathbb{R}^{n}}\si'(\mathbb F-\mathbb I):\nabla u\,dxd\tau\nonumber\\
        & +\int_0^t\int_{\mathbb{R}^{n}}\sigma (\mathbb F-\mathbb I)_{\tau}:\nabla u\,dxd\tau
        -\int_0^t\int_{\mathbb{R}^{n}}\sigma (\mathbb F-\mathbb I):\nabla(u\cdot \nabla u)\,dxd\tau\nonumber\\
        & -\int_{\mathbb{R}^{n}}\sigma ((\mathbb F-\mathbb I)(\mathbb F-\mathbb I)^\top):\nabla u\,dx
        +\int_0^t \int_{\mathbb{R}^{n}}\si'((\mathbb F-\mathbb I)(\mathbb F-\mathbb I)^\top):\nabla u\,dxd\tau\nonumber\\
        & +\int_0^t\int_{\mathbb{R}^{n}}\sigma ((\mathbb F-\mathbb I)(\mathbb F-\mathbb I)^\top)_{\tau}:\nabla u\,dxd\tau
        -\int_0^t\int_{\mathbb{R}^{n}}\sigma ((\mathbb F-\mathbb I)(\mathbb F-\mathbb I)^\top):\nabla(u\cdot \nabla u)\,dxd\tau\nonumber\\
        \triangleq& \sum_{i=1}^8 K_{3i}.
    \end{align}
For $K_{31}$ and $K_{32}$, using Cauchy's inequality and \eqref{4.1.1}, we derive
    \begin{align}
        \label{4.1.8}
        K_{31}+K_{32}
        \leq&\ C\Big(\sigma\int_{\mathbb{R}^{n}}|\nabla u||\mathbb F-\mathbb I|\,dx
        +\int_0^\sigma\int_{\mathbb{R}^{n}}|\nabla u||\mathbb F-\mathbb I|\,dxd\tau\Big)\nonumber\\
        \leq&\ \frac{\mu}{8}\int_{\mathbb{R}^{n}}\sigma|\nabla u|^2\,dx
        +C\Big(\int_{\mathbb{R}^{n}}|\mathbb F-\mathbb I|^2\,dx
        +\int_0^{\si(T)}\int_{\mathbb{R}^{n}}|\nabla u|^2\,dxdt
        +\int_0^{\si(T)}\int_{\mathbb{R}^{n}}|\mathbb F-\mathbb I|^2\,dxdt\Big)\nonumber\\
        \leq&\ \frac{\mu}{8}\int_{\mathbb{R}^{n}}\sigma|\nabla u|^2\,dx+C\mathcal{E}_0.
    \end{align}
For $K_{33}$ and $K_{34}$, by the equation $\eqref{1.1}_2$, integration by parts, \eqref{prop4.1.1}, and \eqref{4.1.1}, we have
    \begin{align}
        \label{4.1.9}
        K_{33}+K_{34}
        =&\ \int_0^t\int_{\mathbb{R}^{n}}\sigma (-u\cdot \nabla(\mathbb F-\mathbb I)+\nabla u\mathbb F):\nabla u\,dxd\tau-\int_0^t\int_{\mathbb{R}^{n}}\sigma (\mathbb F-\mathbb I):\nabla(u\cdot \nabla u)\,dxd\tau \nonumber\\
        =&\ \int_0^t\int_{\mathbb{R}^{n}}\sigma \operatorname{div}u(\mathbb F-\mathbb I):\nabla u\,dxd\tau+\int_0^t\int_{\mathbb{R}^{n}}\sigma (\mathbb F-\mathbb I):(u\cdot \nabla (\nabla u))\,dxd\tau\nonumber\\
        &\ +\int_0^t\int_{\mathbb{R}^{n}}\sigma (\nabla u\mathbb F):\nabla u\,dxd\tau-\int_0^t\int_{\mathbb{R}^{n}}\sigma (\mathbb F-\mathbb I):(u\cdot \nabla (\nabla u)+\nabla u\nabla u)\,dxd\tau \nonumber\\
        =&\ \int_0^t\int_{\mathbb{R}^{n}}\sigma (\nabla u\mathbb F):\nabla u\,dxd\tau-\int_0^t\int_{\mathbb{R}^{n}}\sigma (\mathbb F-\mathbb I):(\nabla u\nabla u)\,dxd\tau \nonumber\\
        \leq&\ C\Big(\sup\limits_{0\leq t\leq T}\|\mathbb{F}-\mathbb{I}\|_{L^\infty(\mathbb{R}^n)}+1\Big)
        \int_0^T\int_{\mathbb{R}^{n}}|\nabla u|^2\,dxdt\nonumber\\
        \leq&\ C\mathcal{E}_0.
    \end{align}
For $K_{35}$ and $K_{36}$, it follows from Cauchy's inequality and \eqref{4.1.1} that
    \begin{align}
        \label{4.1.10}
        K_{35}+K_{36}
        \leq&\ C\Big(\sigma\int_{\mathbb{R}^{n}}|\nabla u||\mathbb F-\mathbb I|^2\,dx
        +\int_0^\sigma\int_{\mathbb{R}^{n}}|\nabla u||\mathbb F-\mathbb I|^2\,dxd\tau\Big)\nonumber\\
        \leq&\ \frac{\mu}{8}\int_{\mathbb{R}^{n}}\sigma|\nabla u|^2\,dx
        +C\Big(\sup\limits_{0\leq t\leq \si(T)}\|\mathbb{F}-\mathbb{I}\|_{L^\infty(\mathbb{R}^n)}^2+1\Big)\nonumber\\
        &\ \times\Big(\int_{\mathbb{R}^{n}}|\mathbb F-\mathbb I|^2\,dx
        +\int_0^{\si(T)}\int_{\mathbb{R}^{n}}|\nabla u|^2\,dxdt
        +\int_0^{\si(T)}\int_{\mathbb{R}^{n}}|\mathbb F-\mathbb I|^2\,dxdt\Big)\nonumber\\
        \leq&\ \frac{\mu}{8}\int_{\mathbb{R}^{n}}\sigma|\nabla u|^2\,dx
        +C\mathcal{E}_0.
    \end{align}
For $K_{37}$ and $K_{38}$, employing the arguments similar to the ones used in \eqref{4.1.9}, one gets
    \begin{align}
        \label{4.1.11}
        K_{37}+K_{38}
        =&\ \int_0^t\int_{\mathbb{R}^{n}}\sigma
        [\nabla u\mathbb F(\mathbb F-\mathbb I)^{\top}+(\mathbb F-\mathbb I)(\nabla u\mathbb F)^{\top}]:\nabla u\,dxd\tau\nonumber\\
        &\ -\int_0^t\int_{\mathbb{R}^{n}}\sigma ((\mathbb F-\mathbb I)(\mathbb F-\mathbb I)^{\top}):(\nabla u\nabla u)\,dxd\tau\nonumber\\
        \leq&\ C\sup\limits_{0\leq t\leq T}\Big(\|\mathbb{F}-\mathbb{I}\|_{L^\infty(\mathbb{R}^n)}^2
        +\|\mathbb{F}-\mathbb{I}\|_{L^\infty(\mathbb{R}^n)}\Big)
        \int_0^T\int_{\mathbb{R}^{n}}|\nabla u|^2\,dxdt\nonumber\\
        \leq&\ C\mathcal{E}_0.
    \end{align}
Putting \eqref{4.1.8}-\eqref{4.1.11} into \eqref{4.1.7} gives
    \begin{equation}
        \label{4.1.12}
        K_3
        \leq \frac{\mu}{4}\int_{\mathbb{R}^{n}}\sigma|\nabla u|^2\,dx+C\mathcal{E}_0.
    \end{equation}
For $K_{4}$, we deduce from integration by parts, \eqref{F-0127}, \eqref{3.4.3-0311}, \eqref{prop4.1.1}, and \eqref{4.1.1} that
    \begin{align}
        \label{4.1.13}
        K_{4}
        =&\ \int_0^t\int_{\mathbb{R}^{n}}\sigma\nabla(-\Delta)^{-1}\diver\diver(\mathbb F\mathbb F^{\top}-\mathbb{I})\cdot\dot u\,dxd\tau\nonumber\\
        =&\ -\int_0^t\int_{\mathbb{R}^{3}}\sigma(-\Delta)^{-1}\diver\diver(\mathbb F\mathbb F^{\top}-\mathbb{I})((\nabla u)^\top:\nabla u)\,dxd\tau\nonumber\\
        \leq&\ C \int_{0}^{t}\|(-\Delta)^{-1}\diver\diver(\mathbb F\mathbb F^{\top}-\mathbb{I})\|_{BMO(\mathbb{R}^n)}
        \|(\nabla u)^\top:\nabla u\|_{\mathcal{H}^{1}(\mathbb{R}^n)}\,d\tau\nonumber\\
        \leq&\ C\int_{0}^{t}\|\mathbb F\mathbb F^\top-\mathbb{I}\|_{L^\infty(\mathbb{R}^n)}
        \|\nabla u\|_{L^2(\mathbb{R}^n)}^2\,d\tau\nonumber\\
        \leq&\ C\sup\limits_{0\leq t\leq T}\Big(\|\mathbb{F}-\mathbb{I}\|_{L^\infty(\mathbb{R}^n)}^2
        +\|\mathbb{F}-\mathbb{I}\|_{L^\infty(\mathbb{R}^n)}\Big)
        \int_{0}^{T}\|\nabla u\|_{L^2(\mathbb{R}^n)}^2\,dt\nonumber\\
        \leq&\ C\mathcal{E}_0,
    \end{align}
where in the second inequality we have used the boundedness of the Riesz transform on the $BMO(\mathbb{R}^n)$ space, the embedding $L^\infty(\mathbb{R}^n)\hookrightarrow BMO(\mathbb{R}^n)$, and the following estimate
    \begin{equation*}
        \|(\nabla u)^\top:\nabla u\|_{\mathcal{H}^{1}(\mathbb{R}^n)}
        =\|\pa_iu^j\pa_ju^i\|_{\mathcal{H}^{1}(\mathbb{R}^n)}
        \leq C\|\nabla u\|_{L^2(\mathbb{R}^n)}^2
    \end{equation*}
by $\diver u=0$ and Lemma \ref{RP}.
Substituting \eqref{4.1.5}, \eqref{4.1.6}, \eqref{4.1.12}, and \eqref{4.1.13} into \eqref{4.1.4} implies
    \begin{equation}
        \label{4.1.14}
        A_{11}(T)
        \leq C\mathcal{E}_0+C\left|\int_0^T\int_{\mathbb{R}^{n}}\sigma \na u:(\na u\na u)\,dxdt\right|.
    \end{equation}

Next, we focus on $A_{12}(T)$.
Applying the material derivative (see \eqref{n-m}) to the equation $\eqref{4.1.3}$ yields
    \begin{equation}
        \label{4.1.15}
        (\mathcal{P}\dot{u})_t+u\cdot\nabla(\mathcal{P}\dot{u})
        =\mu[\Delta\dot{u}-\De(u\cdot\na u)+u\cdot\nabla(\Delta u)]
        +(\mathcal{P}\diver(\mathbb F\mathbb F^{\top}))_t
        +u\cdot\nabla(\mathcal{P}\diver(\mathbb F\mathbb F^{\top})).
    \end{equation}
Multiplying (\ref{4.1.15}) by $\sigma^n\mathcal P \dot u$, and integrating the resulting equation over $\mathbb{R}^n\times [0, t]$ with $0\leq t\leq T$, we obtain
    \begin{align}
        \label{4.1.16}
        &\ \frac{1}{2}\int_{\mathbb{R}^{n}}\sigma^n |\mathcal{P}\dot{u}|^2\,dxd\tau
        -\frac{n}{2}\int_{0}^{t}\int_{\mathbb{R}^{n}}
        \si^{n-1}\si'|\mathcal{P}\dot{u}|^2\,dxd\tau\nonumber\\
        =&\ \mu\int_{0}^{t}\int_{\mathbb{R}^{n}}\sigma^n[\Delta\dot{u}-\De(u\cdot\na u)+u\cdot\nabla(\Delta u)]\cdot\mathcal P \dot u\,dxd\tau\nonumber\\
        &+\int_{0}^{t}\int_{\mathbb{R}^{n}}\sigma^n[\mathcal{P}\diver(\mathbb F\mathbb F^{\top})_{\tau}+u\cdot\nabla(\mathcal{P}\diver(\mathbb F\mathbb F^{\top}))]\cdot\mathcal P \dot u\,dxd\tau\nonumber\\
        =&\ \mu\int_{0}^{t}\int_{\mathbb{R}^{n}}\sigma^n\Delta\dot{u}\cdot\mathcal P \dot u\,dxd\tau
        +\mu\int_{0}^{t}\int_{\mathbb{R}^{n}}\sigma^n[-\De(u\cdot\na u)+u\cdot\nabla(\Delta u)]\cdot\mathcal P \dot u\,dxd\tau\nonumber\\
        &\ +\int_{0}^{t}\int_{\mathbb{R}^{n}}\sigma^n[\diver(\mathbb F\mathbb F^{\top})_{\tau}+u\cdot\nabla(\diver(\mathbb F\mathbb F^{\top}))]\cdot\mathcal P \dot u\,dxd\tau\nonumber\\
        &\ -\int_{0}^{t}\int_{\mathbb{R}^{n}}
        \sigma^n(u\cdot\nabla(\mathcal{Q}\diver(\mathbb F\mathbb F^{\top})))\cdot\mathcal P \dot u\,dxd\tau\nonumber\\
        \triangleq&\ \sum_{i=1}^4 L_i.
    \end{align}
For $L_1$, employing integration by parts yields
    \begin{equation}
        \label{4.1.17}
        L_1
        =-\mu\int_{0}^{t} \int_{\mathbb{R}^{n}} \sigma^{n}(\nabla \mathcal{P}\dot{u}+\nabla \mathcal{Q}\dot{u}):\nabla \mathcal{P} \dot{u}\,dxd\tau
        =-\mu\int_{0}^{t} \int_{\mathbb{R}^{n}} \sigma^{n}|\nabla \mathcal{P} \dot{u}|^2\,dxd\tau,
    \end{equation}
where one has used the following identity
    \begin{equation*}
        -\mu\int_{0}^{t} \int_{\mathbb{R}^{n}} \sigma^{n}\nabla \mathcal{Q}\dot{u}:\nabla \mathcal{P} \dot{u}\,dxd\tau
        =\mu\int_{0}^{t} \int_{\mathbb{R}^{n}} \sigma^{n}\mathcal{Q}\dot{u}\cdot \mathcal{P} \Delta\dot{u}\,dxd\tau
        =0
    \end{equation*}
by
the orthogonality between $\mathcal{P}$ and $\mathcal{Q}$ (see \eqref{n-p} and \eqref{n-q}).
For $L_2$, it follows from integration by parts, $\diver u=0$, and Cauchy's inequality that
    \begin{align}
        \label{4.1.18}
        L_2
        =&\ -\mu\int_{0}^{t} \int_{\mathbb{R}^{n}} \sigma^{n}[\na(u \cdot \nabla u):\nabla \mathcal{P} \dot{u}+(u\cdot\na(\mathcal{P} \dot{u}))\cdot\De u]\,dxd\tau\nonumber\\
        =&\ \mu\int_{0}^{t} \int_{\mathbb{R}^{n}} \sigma^{n}\Big\{
        [\na u\na u+u\cdot\na (\na u)]:\na \mathcal{P} \dot{u}
        +\na(u\cdot\na(\mathcal{P} \dot{u})):\na u\Big\}\,dxd\tau\nonumber\\
        =&\ \mu\int_{0}^{t} \int_{\mathbb{R}^{n}} \sigma^{n}\Big\{(\na u\na u):\na \mathcal{P} \dot{u}+[-u\cdot\na(\na\mathcal{P} \dot{u})+\na(u\cdot\na(\mathcal{P} \dot{u}))]:\na u\Big\}\,dxd\tau\nonumber\\
        =&\ \mu\int_{0}^{t} \int_{\mathbb{R}^{n}} \sigma^{n}[(\na u\na u):\na \mathcal{P} \dot{u}+(\nabla\mathcal{P} \dot{u}\nabla u) :\nabla u]\,dxd\tau\nonumber\\
        \leq&\ \frac{\mu}{6}\int_{0}^{t} \int_{\mathbb{R}^{n}} \sigma^{n}|\nabla \mathcal{P} \dot{u}|^2\,dxd\tau+C\int_{0}^{t} \int_{\mathbb{R}^{n}} \sigma^{n}|\nabla u|^4\,dxd\tau.
    \end{align}
For $L_3$, utilizing integration by parts, $\diver u=0$, the equation \eqref{3.3.5}, Cauchy's inequality, \eqref{prop4.1.1}, and \eqref{4.1.1}, we deduce
    \begin{align}
        \label{4.1.19}
        L_3
        =&\ -\int_{0}^{t}\int_{\mathbb{R}^{n}}\sigma^n[
        (\mathbb F\mathbb F^{\top})_{\tau}:\na\mathcal P \dot u
        +(u\cdot\nabla(\mathcal P \dot u))\cdot\diver(\mathbb F\mathbb F^{\top})]\,dxd\tau\nonumber\\
        =&\ \int_{0}^{t}\int_{\mathbb{R}^{n}}\sigma^n\Big\{
        \big[u\cdot\na(\mathbb{F}\mathbb{F}^\top)
        -\na u(\mathbb{F}\mathbb{F}^\top)
        -(\mathbb{F}\mathbb{F}^\top)(\na u)^\top\big]:\na\mathcal P \dot u
        +\na(u\cdot\nabla(\mathcal P \dot u)):(\mathbb F\mathbb F^{\top})\Big\}\,dxd\tau\nonumber\\
        =&\ \int_{0}^{t}\int_{\mathbb{R}^{n}}\sigma^n\big(-\na u(\mathbb{F}\mathbb{F}^\top)
        -(\mathbb{F}\mathbb{F}^\top)(\na u)^\top\big):\na(\mathcal P \dot u)\,dxd\tau\nonumber\\
        &\ +\int_{0}^{t}\int_{\mathbb{R}^{n}}\sigma^n\big(-u\cdot\na(\na\mathcal P \dot u)+\na(u\cdot\nabla(\mathcal P \dot u))\big):(\mathbb F\mathbb F^{\top})\,dxd\tau\nonumber\\
        =&\ \int_{0}^{t}\int_{\mathbb{R}^{n}}\sigma^n\Big[\big(-\na u(\mathbb{F}\mathbb{F}^\top)
        -(\mathbb{F}\mathbb{F}^\top)(\na u)^\top\big):\na(\mathcal P \dot u)
        +(\nabla\mathcal{P} \dot{u}\nabla u) :(\mathbb{F}\mathbb{F}^\top)\Big]\,dxd\tau\nonumber\\
        \leq &\ \frac{\mu}{6}\int_{0}^{t} \int_{\mathbb{R}^{n}} \sigma^{n}|\nabla \mathcal{P} \dot{u}|^2\,dxd\tau+C\Big(\sup\limits_{0\leq t\leq T}\|\mathbb{F}-\mathbb{I}\|_{L^\infty(\mathbb{R}^n)}^4+1\Big)
        \int_{0}^{T} \int_{\mathbb{R}^{n}} \sigma^{n}|\nabla u|^2\,dxdt\nonumber\\
        \leq &\ \frac{\mu}{6}\int_{0}^{t} \int_{\mathbb{R}^{n}} \sigma^{n}|\nabla \mathcal{P} \dot{u}|^2\,dxd\tau
        +C\mathcal{E}_0.
	\end{align}
For $L_4$, by integration by parts, $\diver u=0$, \eqref{F-0127}, \eqref{prop4.1.1}, and \eqref{4.1.1}, we have
    \begin{align}
        \label{4.1.20}
        L_4
        =&\ \int_{0}^{t}\int_{\mathbb{R}^{n}}\sigma^n(u\cdot\nabla)\nabla(-\Delta)^{-1}\diver\diver(\mathbb F\mathbb F^{\top}-\mathbb{I}) \cdot\mathcal P \dot u\,dxd\tau\nonumber\\
        =&\ \int_{0}^{t}\int_{\mathbb{R}^{n}}\sigma^n
        (-\Delta)^{-1}\diver\diver(\mathbb F\mathbb F^{\top}-\mathbb{I})
        \diver(u\cdot\nabla(\mathcal P \dot u))\,dxd\tau\nonumber\\
        \leq&\ C\int_{0}^{t}\sigma^n\|(-\Delta)^{-1}\diver\diver(\mathbb F\mathbb F^{\top}-\mathbb{I})\|_{BMO(\mathbb{R}^n)}
        \|\diver(u\cdot\nabla(\mathcal P \dot u))\|_{\mathcal{H}^1(\mathbb{R}^n)}\,d\tau\nonumber\\
        \leq&\ C\int_{0}^{t}\sigma^n\|\mathbb F\mathbb F^{\top}-\mathbb{I}\|_{L^\infty(\mathbb{R}^n)}
        \|\nabla u\|_{L^2(\mathbb{R}^n)}\|\nabla\mathcal P \dot u\|_{L^2(\mathbb{R}^n)}\,d\tau\nonumber\\
        \leq&\ \frac{\mu}{6}\int_{0}^{t} \int_{\mathbb{R}^{n}} \sigma^{n}|\nabla \mathcal{P} \dot{u}|^2\,dxd\tau
        +C\sup\limits_{0\leq t\leq T}\Big(\|\mathbb{F}-\mathbb{I}\|_{L^\infty(\mathbb{R}^n)}^4
        +\|\mathbb{F}-\mathbb{I}\|_{L^\infty(\mathbb{R}^n)}^2\Big)
        \int_{0}^{T} \int_{\mathbb{R}^{n}} \sigma^{n}|\nabla u|^2\,dxdt\nonumber\\
        \leq&\ \frac{\mu}{6}\int_{0}^{t} \int_{\mathbb{R}^{n}} \sigma^{n}|\nabla \mathcal{P} \dot{u}|^2\,dxd\tau
        +C\mathcal{E}_0,
    \end{align}
where in the second inequality we have used the boundedness of the Riesz transform on the $BMO({\mathbb{R}^n})$ space, the embedding $L^\infty(\mathbb{R}^n)\hookrightarrow BMO(\mathbb{R}^n)$, and the following estimate
    \begin{equation*}
        \|\diver(u\cdot\nabla(\mathcal P \dot u))\|_{\mathcal{H}^1(\mathbb{R}^n)}
        =\|(\nabla u)^\top:\nabla(\mathcal P \dot u)\|_{\mathcal{H}^{1}(\mathbb{R}^n)}
        =\|\pa_iu^j\pa_j(\mathcal P \dot u)^i\|_{\mathcal{H}^{1}(\mathbb{R}^n)}
        \leq C\|\nabla u\|_{L^2(\mathbb{R}^n)}\|\nabla\mathcal P \dot u\|_{L^2(\mathbb{R}^n)}
    \end{equation*}
due to $\diver(\mathcal{P}\dot{u})=0$ and Lemma \ref{RP}.
Plugging \eqref{4.1.17}-\eqref{4.1.20} into \eqref{4.1.16}, and noticing that
    \begin{equation*}
        \frac{n}{2}\int_{0}^{t}\int_{\mathbb{R}^{n}}
        \si^{n-1}\si'|\mathcal{P}\dot{u}|^2\,dxd\tau
        \leq C\int_{0}^{t}\int_{\mathbb{R}^{n}}\si|\mathcal{P}\dot{u}|^2\,dxd\tau
        \leq CA_{11}(T),
    \end{equation*}
one thus derives
    \begin{equation}
        \label{4.1.21}
        A_{12}(T)\leq C\Big(\mathcal{E}_0+A_{11}(T)+\int_0^T\int_{\mathbb{R}^{n}}\si^n|\na u|^4\,dxdt\Big).
    \end{equation}
Combining \eqref{4.1.14} with \eqref{4.1.21}, we conclude \eqref{4.1.2} and hence complete the proof of Lemma \ref{le4.1}.
\end{proof}

The following lemma shows the boundedness of $A_2(\si(T))$, which will be used to establish some short-time estimates on the gradient and the material derivative of the velocity.
\begin{Lemma}\label{le4.2}
Under the conditions of Proposition \ref{prop4.1}, there exists a positive constant $\de_1$ such that if $\mathcal{E}_0\leq\de_1$, then
    \begin{equation}
        \label{4.2.1}
        A_2(\si(T))
        =\sup\limits_{0\leq t\leq \si(T)}\int_{\mathbb{R}^{n}}|u|^q\,dx
        +\int_0^{\si(T)}\int_{\mathbb{R}^n}\Big(|\na u|^2|u|^{q-2}
        +\big|\na|u|^{\frac{q}{2}}\big|^2\Big)\,dxdt
        \leq \mathcal{E}_0^\al,
    \end{equation}
where $q$ is defined in \eqref{q-0404}.
\end{Lemma}

\begin{proof}
By the identity \eqref{4.1.7-0311}, we can rewrite the equation \eqref{4.1.3} as follows:
    \begin{equation}
        \label{4.2.2}
        u_t+u\cdot\na u-\mu\Delta u
        =\mathcal{Q}(u\cdot\na u)+\mathcal{P}\diver((\mathbb F-\mathbb I)(\mathbb F-\mathbb I)^{\top})+\diver(\mathbb F-\mathbb I),
    \end{equation}
where we have used $\mathcal{P}\diver(\mathbb F-\mathbb I)=\diver(\mathbb F-\mathbb I)$ since $\diver\diver(\mathbb F-\mathbb I)=0$.
Multiplying \eqref{4.2.2} by $q|u|^{q-2}u$, and integrating the resulting equation over $\mathbb{R}^n\times[0, t]$ with $0\leq t\leq\sigma(T)$, by integration by parts and $\diver u=0$, one deduces
    \begin{align}
        \label{4.2.3}
        &\ \int_{\mathbb{R}^{n}}|u|^q\,dx\Big|_0^t
        +q\mu\int_0^t\int_{\mathbb{R}^{n}}|\nabla u|^2|u|^{q-2}\,dxd\tau
        +\frac{4(q-2)}{q}\mu\int_0^t\int_{\mathbb{R}^{n}}
        \big|\nabla|u|^{\frac{q}{2}}\big|^2\,dxd\tau
        \nonumber\\
        =&\ q\int_0^t\int_{\mathbb{R}^{n}}\Big[\mathcal{Q}(u\cdot\nabla u)+\mathcal{P}\diver((\mathbb F-\mathbb I)(\mathbb F-\mathbb I)^{\top})+\diver(\mathbb F-\mathbb I)\Big]\cdot |u|^{q-2}u\,dxd\tau\nonumber\\
        =&\ q\int_0^t\int_{\mathbb{R}^{n}}\Big[-\na(-\Delta)^{-1}\diver\diver(u\otimes u)
        +\diver((\mathbb F-\mathbb I)(\mathbb F-\mathbb I)^{\top})\nonumber\\
        &\ +\nabla(-\Delta)^{-1}\diver\diver((\mathbb F-\mathbb I)(\mathbb F-\mathbb I)^{\top})
        +\diver(\mathbb F-\mathbb I)\Big]\cdot |u|^{q-2}u\,dxd\tau\nonumber\\
        =&\
        q\int_0^t\int_{\mathbb{R}^{n}}(-\Delta)^{-1}\diver\diver(u\otimes u) u\cdot\na|u|^{q-2}\,dxd\tau
        -q\int_0^t\int_{\mathbb{R}^{n}}(\mathbb F-\mathbb I)(\mathbb F-\mathbb I)^{\top}:\na(|u|^{q-2} u)\,dxd\tau\nonumber\\
        &\
        -q\int_0^t\int_{\mathbb{R}^{n}}(-\Delta)^{-1}\diver\diver((\mathbb F-\mathbb I)(\mathbb F-\mathbb I)^{\top})u\cdot\na|u|^{q-2}\,dxd\tau
        -q\int_0^t\int_{\mathbb{R}^{n}}(\mathbb F-\mathbb I):\na(|u|^{q-2} u)\,dxd\tau\nonumber\\
        \triangleq&\ \sum_{i=1}^4 R_i.
    \end{align}
For $R_1$, by H\"{o}lder's and Young's inequalities, the $L^p$-boundedness of the Riesz transform, Lemma \ref{le2.3}, and \eqref{prop4.1.1}, we derive
    \begin{align}
        \label{4.2.4}
        R_1
        \leq&\ C\Big(\int_0^t\int_{\mathbb{R}^{n}}\big|\nabla|u|^{\frac{q}{2}}\big|^2\,dxd\tau\Big)^\frac{1}{2}
        \Big(\int_0^t\int_{\mathbb{R}^{n}}|u|^{q+2}\,dxd\tau\Big)^\frac{q-2}{2(q+2)}
        \Big(\int_0^t\int_{\mathbb{R}^{n}}|(-\Delta)^{-1}\diver\diver(u\otimes u)|^{\frac{q+2}{2}}\,dxd\tau\Big)^\frac{2}{q+2}\nonumber\\
        \leq&\ \frac{\mu}{4}\int_0^t\int_{\mathbb{R}^{n}}\big|\nabla|u|^{\frac{q}{2}}\big|^2\,dxd\tau
        +C\int_0^{t}\int_{\mathbb{R}^{n}}|u|^{q+2}\,dxd\tau\nonumber\\
        \leq&\ \frac{\mu}{4}\int_0^t\int_{\mathbb{R}^{n}}\big|\nabla|u|^{\frac{q}{2}}\big|^2\,dxd\tau
        +C\int_0^{t}\Big(\int_{\mathbb{R}^{n}}|u|^q\,dx\Big)^{\frac{q+2-n}{q}}
        \Big(\int_{\mathbb{R}^{n}}\big|\nabla|u|^{\frac{q}{2}}\big|^2\,dx\Big)^{\frac{n}{q}}\,d\tau\nonumber\\
        \leq&\ \frac{\mu}{4}\int_0^t\int_{\mathbb{R}^{n}}\big|\nabla|u|^{\frac{q}{2}}\big|^2\,dxd\tau
        +C\Big(\sup\limits_{0\leq t\leq \si(T)}\int_{\mathbb{R}^{n}}|u|^q\,dx\Big)^{\frac{q+2-n}{q}}
        \Big(\int_0^{\si(T)}\int_{\mathbb{R}^{n}}
        \big|\nabla|u|^{\frac{q}{2}}\big|^2\,dxd\tau\Big)^\frac{n}{q}
        \nonumber\\
        \leq&\ \frac{\mu}{4}\int_0^t\int_{\mathbb{R}^{n}}\big|\nabla|u|^{\frac{q}{2}}\big|^2\,dxd\tau
        +CA_2(\si(T))^{\frac{q+2}{q}}\nonumber\\
        \leq&\ \frac{\mu}{4}\int_0^t\int_{\mathbb{R}^{n}}\big|\nabla|u|^{\frac{q}{2}}\big|^2\,dxd\tau
        +C\mathcal{E}_0^{\frac{(q+2)\al}{q}}.
    \end{align}
For $R_2$, since $|\nabla (|u|u)|\leq C|\nabla u ||u|$, utilizing Cauchy's and H\"{o}lder's inequalities with \eqref{prop4.1.1}, one gets
    \begin{align}
        \label{4.2.5}
        R_2
        \leq&\ \frac{\mu}{4}\int_0^t\int_{\mathbb{R}^{n}}|\nabla u|^2|u|^{q-2}\,dxd\tau
        +C\int_0^{\si(T)}\int_{\mathbb{R}^{n}}
        |u|^{q-2}|\mathbb{F}-\mathbb{I}|^4\,dxdt\nonumber\\
        \leq&\ \frac{\mu}{4}\int_0^t\int_{\mathbb{R}^{n}}|\nabla u|^2|u|^{q-2}\,dxd\tau
        +C\sup\limits_{0\leq t\leq \si(T)}\Big(\|u\|_{L^q(\mathbb{R}^n)}^{q-2}
        \|\mathbb{F}-\mathbb{I}\|_{L^\infty(\mathbb{R}^n)}^{\frac{4(q-1)}{q}}
        \|\mathbb{F}-\mathbb{I}\|_{L^2(\mathbb{R}^n)}^{\frac{4}{q}}\Big)
        \nonumber\\
        \leq&\ \frac{\mu}{4}\int_0^t\int_{\mathbb{R}^{n}}|\nabla u|^2|u|^{q-2}\,dxd\tau
        +CA_2(\si(T))^{\frac{q-2}{q}}
        B(T)^{\frac{2(q-1)}{q}}
        \mathcal{E}_0^{\frac{2}{q}}\nonumber\\
        \leq&\ \frac{\mu}{4}\int_0^t\int_{\mathbb{R}^{n}}|\nabla u|^2|u|^{q-2}\,dxd\tau
        +C\mathcal{E}_0^{\frac{2+(q-2)\al+2(q-1)\be}{q}}.
    \end{align}
For $R_3$ and $R_4$, it follows from the arguments analogous to those used in \eqref{4.2.4} and \eqref{4.2.5} that
    \begin{align}
        \label{4.2.5-0327}
        R_3
        \leq&\ C\Big(\int_0^t\int_{\mathbb{R}^{n}}\big|\nabla|u|^{\frac{q}{2}}\big|^2\,dxd\tau\Big)^\frac{1}{2}
        \Big(\int_0^t\int_{\mathbb{R}^{n}}|u|^{q}\,dxd\tau\Big)^\frac{q-2}{2q}
        \Big(\int_0^t\int_{\mathbb{R}^{n}}|(-\Delta)^{-1}\diver\diver((\mathbb F-\mathbb I)(\mathbb F-\mathbb I)^{\top})|^q\,dxd\tau\Big)^\frac{1}{q}
        \nonumber\\
        \leq&\ \frac{\mu}{4}\int_0^t\int_{\mathbb{R}^{n}}\big|\nabla|u|^{\frac{q}{2}}\big|^2\,dxd\tau
        +C\Big(\int_0^{\si(T)}\int_{\mathbb{R}^{n}}|u|^q\,dxdt\Big)^{\frac{q-2}{q}}
        \Big(\int_0^{\si(T)}\int_{\mathbb{R}^{n}}
        |\mathbb{F}-\mathbb{I}|^{2q}\,dxdt\Big)^{\frac{2}{q}}\nonumber\\
        \leq&\ \frac{\mu}{4}\int_0^t\int_{\mathbb{R}^{n}}\big|\nabla|u|^{\frac{q}{2}}\big|^2\,dxd\tau
        +C\sup\limits_{0\leq t\leq \si(T)}\Big(\|u\|_{L^q(\mathbb{R}^n)}^{q-2}
        \|\mathbb{F}-\mathbb{I}\|_{L^\infty(\mathbb{R}^n)}^{\frac{4(q-1)}{q}}
        \|\mathbb{F}-\mathbb{I}\|_{L^2(\mathbb{R}^n)}^{\frac{4}{q}}\Big)\nonumber\\
        \leq&\ \frac{\mu}{4}\int_0^t\int_{\mathbb{R}^{n}}\big|\nabla|u|^{\frac{q}{2}}\big|^2\,dxd\tau
        +C\mathcal{E}_0^{\frac{2+(q-2)\al+2(q-1)\be}{q}},
    \end{align}
and
    \begin{align}
        \label{4.2.6}
        R_4
        \leq&\ \frac{\mu}{4}\int_0^t\int_{\mathbb{R}^{n}}|\nabla u|^2|u|^{q-2}\,dxd\tau
        +C\int_0^{\si(T)}\int_{\mathbb{R}^{n}}|u|^{q-2}|\mathbb{F}-\mathbb{I}|^2\,dxdt\nonumber\\
        \leq&\ \frac{\mu}{4}\int_0^t\int_{\mathbb{R}^{n}}|\nabla u|^2|u|^{q-2}\,dxd\tau
        +C\sup\limits_{0\leq t\leq \si(T)}\Big(\|u\|_{L^q(\mathbb{R}^n)}^{q-2}
        \|\mathbb{F}-\mathbb{I}\|_{L^\infty(\mathbb{R}^n)}^{\frac{2(q-2)}{q}}
        \|\mathbb{F}-\mathbb{I}\|_{L^2(\mathbb{R}^n)}^{\frac{4}{q}}\Big)
        \nonumber\\
        \leq&\ \frac{\mu}{4}\int_0^t\int_{\mathbb{R}^{n}}|\nabla u|^2|u|^{q-2}\,dxd\tau
        +C\mathcal{E}_0^{\frac{2+(q-2)\al+(q-2)\be}{q}}.
    \end{align}
In addition, by the standard $L^p$ interpolation inequality and the Sobolev embedding theorem, one has for $n=2$,
    \begin{equation}
        \label{4.2.7}
        \|u_0\|_{L^{\frac{4}{2-\eta}}(\mathbb{R}^2)}^{\frac{4}{2-\eta}}\leq C\|u_0\|_{L^2(\mathbb{R}^2)}^{\frac{2}{2-\eta}}
        \|u_0\|_{\dot{H}^\eta(\mathbb{R}^2)}^{\frac{2}{2-\eta}}
        \leq C\mathcal{E}_0^{2\al},
    \end{equation}
where $\eta\in(0,1]$ and $\al=\frac{1}{2(2-\eta)}$.
Similarly, for $n=3$,
    \begin{equation}
        \label{4.2.7-0403}
        \|u_0\|_{L^3(\mathbb{R}^3)}^3\leq C\|u_0\|_{L^2(\mathbb{R}^3)}^\frac{3(2\eta-1)}{2\eta}
        \|u_0\|_{\dot{H}^\eta(\mathbb{R}^3)}^\frac{3}{2\eta}
        \leq C\mathcal{E}_0^{2\al},
    \end{equation}
where $\eta\in(\frac{8}{15},1]$ and $\al=\frac{3(2\eta-1)}{8\eta}$.
Substituting \eqref{4.2.4}-\eqref{4.2.7-0403} into \eqref{4.2.3}, we obtain
    \begin{equation}
        \label{4.2.8}
        \sup\limits_{0\leq t\leq \si(T)}\int_{\mathbb{R}^{n}}|u|^q\,dx
        +\int_0^{\si(T)}\int_{\mathbb{R}^n}\Big(|\na u|^2|u|^{q-2}
        +\big|\nabla|u|^{\frac{q}{2}}\big|^2\Big)\,dxdt
        \leq C_1\mathcal{E}_0^{\frac{5\al}{3}}
    \end{equation}
since $0<\al<\frac{1}{2}$.
Therefore, it follows from \eqref{4.2.8} that \eqref{4.2.1} holds if
    \begin{equation*}
        \mathcal{E}_0
        \leq\de_1
        \triangleq
        \min\Big\{1,C_1^{-\frac{3}{2\al}}\Big\}.
    \end{equation*}
The proof of Lemma \ref{le4.2} is finished.
\end{proof}

Define the effective viscous fluxes by
    \begin{equation}
        \label{ck0}
        \begin{cases}
            \ \mathcal{G}_1
            \triangleq\mu\nabla u-
            (-\De)^{-1}\na{\mathcal{P}}\diver(\mathbb{F}\mathbb{F}^{\top}-\mathbb{I}),\\
            \ \mathcal{G}_2
            \triangleq\mu\cur u-
            (-\De)^{-1}\cur\diver(\mathbb{F}\mathbb{F}^{\top}-\mathbb{I}),
        \end{cases}
    \end{equation}
where $\cur u=\na u-(\na u)^\top$.
It is easily check from $\eqref{1.1}_1$ that
    \begin{equation*}
        \Delta \mathcal{G}_1=\na{\mathcal{P}}\dot{u}\quad \text{and} \quad
        \Delta \mathcal{G}_2=\cur(\mathcal P\dot u).
    \end{equation*}
By \eqref{2.5.1}-\eqref{2.5.2} in Lemma \ref{le2.5}, one has for $i=1,2$,
    \begin{equation}
        \label{ck}
        \|\na^2\mathcal{G}_i\|_{L^2(\mathbb{R}^n)}
        \leq C\|\nabla \mathcal{P}\dot{u}\|_{L^2(\mathbb{R}^n)},
    \end{equation}
and
    \begin{equation}
        \label{ck1}
        \|\nabla\mathcal{G}_i\|_{L^2(\mathbb{R}^n)}
        \leq C\|\mathcal{P}\dot{u}\|_{L^2(\mathbb{R}^n)}.
    \end{equation}

Let $\mathbb{E}\triangleq\mathbb{F}-\mathbb{I}$.
By $\diver u=0$, $\diver{\mathbb F}^{\top}=0$, and the Helmholtz decomposition
    \begin{equation*}
        \Delta =\na\diver+\diver\cur,
    \end{equation*}
we deduce
    \begin{align*}
        \diver{\mathcal G}_2
        =&\ \mu\diver\cur u-(-\Delta)^{-1}\diver\cur\diver({\mathbb{E}}{\mathbb{E}}^{\top})
        -(-\Delta)^{-1}\diver\cur\diver{\mathbb{E}}\nonumber\\
        =&\ \mu\Delta u-(-\Delta)^{-1}\Delta\diver\mathbb E+(-\Delta)^{-1}\na\diver\diver\mathbb E-(-\Delta)^{-1}\diver\cur\diver({\mathbb{E}}{\mathbb{E}}^{\top})\nonumber\\
        =&\ \mu\Delta u+\diver \mathbb{E}
        -(-\Delta)^{-1}\diver\cur\diver({\mathbb{E}}{\mathbb{E}}^{\top}).
    \end{align*}
Applying $\na$ to the above equation yields
    \begin{align*}
        \na\diver{\mathcal G}_2
        =&\ \mu\na\Delta u+\na\diver \mathbb{E}-\na(-\Delta)^{-1}\diver\cur\diver({\mathbb{E}}{\mathbb{E}}^{\top})\\
        =&\ \Delta(\mu\na u+\mathbb E)-\diver\cur\mathbb E
        -\na(-\Delta)^{-1}\diver\cur\diver({\mathbb{E}}{\mathbb{E}}^{\top}),
    \end{align*}
i.e.,
    \begin{equation}
        \label{329}
        \mu\na u+\mathbb E
        =-(-\Delta)^{-1}\na\diver{\mathcal G}_2
        -(-\Delta)^{-1}\diver\cur\mathbb E
        -(-\Delta)^{-1}\na(-\Delta)^{-1}\diver\cur\diver({\mathbb{E}}{\mathbb{E}}^{\top}).
    \end{equation}
It follows from the `$\cur$' structural condition $\eqref{1.4-1224}_2$ that
    \begin{equation*}
        {\mathbb E}^{lk}\pa_l{\mathbb E}^{ij}+\pa_k{\mathbb E}^{ij}
        ={\mathbb E}^{lj}\pa_l{\mathbb E}^{ik}+\pa_j{\mathbb E}^{ik},
    \end{equation*}
which combined with the `$\diver$' identity \eqref{1.4.1-1224} implies
    \begin{equation*}
        (\diver\cur\mathbb E)^{ij}
        =\pa_k(\pa_k{\mathbb E}^{ij}-\pa_j{\mathbb E}^{ik})
        =\pa_k({\mathbb E}^{lj}\pa_l{\mathbb E}^{ik}-{\mathbb E}^{lk}\pa_l{\mathbb E}^{ij})
        =\pa_k\pa_l({\mathbb E}^{lj}{\mathbb E}^{ik}-{\mathbb E}^{lk}{\mathbb E}^{ij}).
    \end{equation*}
Putting the above equality into \eqref{329} leads to
    \begin{equation}
        \label{324}
        \mu\na u+\mathbb E=\mathcal G+\mathcal R
    \end{equation}
with $\mathbb{E}=\mathbb{F}-\mathbb{I}$ and
    \begin{equation}
        \label{324-0311}
        \begin{cases}
            \
            \mathcal{G}\triangleq-(-\Delta)^{-1}\nabla\operatorname{div}\mathcal{G}_2,\\
            \ \mathcal R^{ij}
            \triangleq-(-\Delta)^{-1}\pa_k\pa_l({\mathbb{E}}^{lj}{\mathbb{E}}^{ik}-
            {\mathbb{E}}^{lk}{\mathbb{E}}^{ij})
            -((-\Delta)^{-1}\na\diver
            (-\Delta)^{-1}\cur\diver({\mathbb{E}}{\mathbb{E}}^{\top}))^{ij}.
        \end{cases}
    \end{equation}

In the following lemma, we give some time-weighted estimates on both the gradient and the material derivative of the velocity for small time.
These estimates are crucial for dealing with the initial layer.
\begin{Lemma}\label{le4.4}
Under the conditions of Proposition \ref{prop4.1}, there exists a positive constant $\de_2\leq\de_1$ such that if $\mathcal{E}_0\leq\de_2$, then
    \begin{equation}
        \label{4.4.1}
        \sup\limits_{0\leq t\leq \sigma(T)}\int_{\mathbb{R}^n}
        t^{1-\eta}|\nabla u|^2\,dx
        +\int_0^{\sigma(T)}\int_{\mathbb{R}^n}
        t^{1-\eta}|\mathcal{P}\dot{u}|^2\,dxdt
        \leq C,
    \end{equation}
and
    \begin{equation}
        \label{4.4.1-0406}
        \sup\limits_{0\leq t\leq \sigma(T)}\int_{\mathbb{R}^n}
        t^{2-\eta}|\mathcal P\dot u|^2\,dx
        +\int_0^{\sigma(T)}\int_{\mathbb{R}^n}
        t^{2-\eta}|\na\mathcal{P}\dot{u}|^2\,dxdt\leq C.
    \end{equation}
\end{Lemma}

\begin{proof}
%
Define a differential operator $\mathcal L$ acting on functions $\om:\mathbb R^n\times\mathbb R^+\to\mathbb R^n$ with $\diver \om=0$
by
    \begin{equation*}
        \mathcal L \om
        =\om_t+ \mathcal P(u\cdot \na \om)-\mu\Delta \om
        =\mathcal P\dot{\om}-\mu\Delta \om.
    \end{equation*}
For some given $\om_{10}$, let $\om_1$ and $\om_2$ be the solutions to the following two initial value problems, respectively,
    \begin{equation}
        \label{w1}
        \begin{cases}
        \ \mathcal{L} \om_1=0,\quad x\in\mathbb{R}^n,\ t>0,\\
        \ \om_1(x,0)=\om_{10}(x),\quad \diver \om_{10}=0,\\
        \ \om_1(x,t)\to0,\quad |x|\to\infty,\ t>0,
        \end{cases}
	\end{equation}
and
    \begin{equation}
        \label{w2}
        \begin{cases}
        \ \mathcal{L} \om_2=\mathcal P\diver(\mathbb F\mathbb F^\top-\mathbb I),
        \quad x\in\mathbb{R}^n,\ t>0,\\
        \ \om_2(x,0)=0,\\
        \ \om_2(x,t)\to0,\quad|x|\to\infty,\ t>0.
        \end{cases}
	\end{equation}
Note that if $\om_{10}=u_0$, then $\om_1+\om_2=u$.

We first show the estimate \eqref{4.4.1}.
Multiplying \eqref{w1} and \eqref{w2} by $\om_1$ and $\om_2$, respectively, and then integrating the resulting equations over $\mathbb{R}^n\times[0,t]$ with $0\leq t\leq \si(T)$, by integration by parts and $\diver \om_1=\diver \om_2=0$, we can deduce
    \begin{equation}
        \label{4.4.3}
        \sup\limits_{0\leq t\leq \si(T)}
        \int_{\mathbb{R}^n}|\om_1|^2\,dx
        +\int_0^{\si(T)}\int_{\mathbb{R}^n}|\na \om_1|^2\,dxd\tau
        \leq C\int_{\mathbb{R}^n}| \om_{10}|^2\,dx,
    \end{equation}
and
    \begin{align}
        \label{4.4.4}
        &\ \sup\limits_{0\leq t\leq \si(T)}
        \int_{\mathbb{R}^n}|\om_2|^2\,dx
        +\int_0^{\si(T)}\int_{\mathbb{R}^n}|\na \om_2|^2\,dxd\tau
        \leq C\int_0^{\si(T)}\int_{\mathbb{R}^n}|\mathbb F \mathbb F^\top-\mathbb I|^2\,dxdt\nonumber\\
        \leq&\ C\Big(\sup\limits_{0\leq t\leq \si(T)}\|\mathbb F-\mathbb I\|_{L^\infty(\mathbb{R}^n)}^2+1\Big)
        \Big(\sup\limits_{0\leq t\leq \si(T)}\|\mathbb{F}-\mathbb{I}\|_{L^2(\mathbb{R}^n)}^2\Big)
        \leq  C\mathcal{E}_0.
    \end{align}
where we have used \eqref{prop4.1.1} and \eqref{4.1.1}.
Moreover, it follows from the equation \eqref{w1} and the inequality \eqref{2.5.1} in Lemma \ref{le2.5} that
    \begin{equation}
        \label{4.4.6-0404}
        \|\na^2 \om_1\|_{L^2(\mathbb{R}^n)}
        \leq C\|\mathcal P\dot \om_1\|_{L^2(\mathbb{R}^n)}.
    \end{equation}

Multiplying \eqref{w1} by $ \om_{1t}$, integrating the resulting equation over $\mathbb R^n$ by parts, and then employing \eqref{4.2.1}, H\"{o}lder's and Young's inequalities, we derive for $q$ defined in \eqref{q-0404},
    \begin{align}
        \label{4.4.5}
        &\ \frac{\mu}{2}\frac{d}{dt}\int_{\mathbb{R}^n}|\na \om_1|^2\,dx
        +\int_{\mathbb{R}^n}|\mathcal P\dot \om_1|^2\,dx\nonumber\\
        =&\ \int_{\mathbb{R}^n}\mathcal P\dot \om_1\cdot\mathcal P(u\cdot \na \om_1)\,dxd\tau\nonumber\\
        \leq&\ C\|\mathcal P\dot \om_1\|_{L^2(\mathbb{R}^n)}\|u\|_{L^q(\mathbb{R}^n)}\|\na \om_1\|_{L^{\frac{2q}{q-2}}(\mathbb{R}^n)}\nonumber\\
        \leq&\ C\mathcal{E}_0^{\frac{\al}{q}}
        \|\mathcal P\dot \om_1\|_{L^2(\mathbb{R}^n)}^{\frac{q+n}{q}}
        \|\na \om_1\|_{L^2(\mathbb{R}^n)}^{\frac{q-n}{q}}\nonumber\\
        \leq&\ C_2\mathcal{E}_0^{\frac{\al}{q}}\|\mathcal P\dot \om_1\|_{L^2(\mathbb{R}^n)}^2
        +C\|\na \om_1\|_{L^2(\mathbb{R}^n)}^2,
    \end{align}
where in the second inequality we have used the following estimate
    \begin{equation*}
        \|\na \om_1\|_{L^q(\mathbb{R}^n)}
        \leq C\|\na \om_1\|_{L^2(\mathbb{R}^n)}^{\frac{q-n}{q}}\|\na^2 \om_1\|_{L^2(\mathbb{R}^n)}^{\frac{n}{q}}
        \leq C\|\na \om_1\|_{L^2(\mathbb{R}^n)}^{\frac{q-n}{q}}\|\mathcal P\dot \om_1\|_{L^2(\mathbb{R}^n)}^{\frac{n}{q}}
    \end{equation*}
due to Lemma \ref{le2.3} and the inequality \eqref{4.4.6-0404}.
Then applying Gronwall's inequality to \eqref{4.4.5} gives
    \begin{equation}
        \label{4.4.7}
        \sup\limits_{0\leq t\leq \si(T)}\int_{\mathbb{R}^n}|\nabla \om_1|^2dx
        +\int_{0}^{\si(T)}\int_{\mathbb{R}^n}|\mathcal P\dot \om_1|^2dxdt
        \leq C\int_{\mathbb{R}^n}|\na\om_{10}|^2\,dx,
    \end{equation}
provided that
    \begin{equation*}
        \mathcal{E}_0
        \leq\de_{2,1}
        \triangleq \min\Big\{\de_1,(2C_2)^{-\frac{q}{\al}}\Big\}.
    \end{equation*}
Furthermore, multiplying \eqref{4.4.5} by $t$, and then
using Gronwall's inequality and \eqref{4.4.3}, we obtain
    \begin{equation}
        \label{4.4.8}
        \sup\limits_{0\leq t\leq \si(T)}\Big(t\int_{\mathbb{R}^n}|\na \om_1|^2dx\Big)
        +\int_{0}^{\si(T)}\int_{\mathbb{R}^n}t|\mathcal P\dot \om_1|^2dxdt
        \leq C\int_0^{\si(T)}\int_{\mathbb{R}^n}|\na \om_1|^2\,dxd\tau
        \leq C\int_{\mathbb{R}^n}| \om_{10}|^2\,dx.
    \end{equation}
Since the solution operator $\om_{10}\mapsto\om_1(\cdot,t)$ is linear,
by \eqref{4.4.7}, \eqref{4.4.8}, and the standard Riesz-Thorin interpolation argument
(cf. \cite{B-L1976}),
one derives
    \begin{equation}
        \label{4.4.9}
        \sup\limits_{0\leq t\leq \si(T)}\Big(t^{1-\eta}\int_{\mathbb{R}^n}|\nabla \om_1|^2dx\Big)
        +\int_{0}^{\si(T)}\int_{\mathbb{R}^n}t^{1-\eta}|\mathcal P\dot \om_1|^2dxdt
        \leq C \|\om_{10}\|_{\dot{H}^\eta(\mathbb{R}^n)}^2.
    \end{equation}

In order to obtain the estimate of $\om_2$,
we define
    \begin{equation}
        \label{4.4.10-0404}
        \mathcal{W}
        \triangleq \mu\nabla \om_2-(-\Delta)^{-1}\nabla{\mathcal{P}}
        \diver(\mathbb{F}\mathbb{F}^{\top}-\mathbb I).
    \end{equation}
Similar to the estimate of $\mathcal {G}_1$ in \eqref{ck1}, one has
    \begin{equation}
        \label{4.4.10}
        \|\nabla\mathcal{W}\|_{L^2(\mathbb{R}^n)}
        \leq C\|\mathcal{P}\dot{\om}_2\|_{L^2(\mathbb{R}^n)}.
    \end{equation}
Multiplying \eqref{w2} by $\mathcal P\dot \om_2$, and integrating the resulting equation over $\mathbb{R}^n\times[0,t]$ with $0\leq t\leq \si(T)$, by \eqref{4.4.4}, the definition of $\mathcal{W}$ \eqref{4.4.10-0404}, and the arguments analogous to those used for \eqref{4.1.14}, we obtain
    \begin{align}
        \label{4.4.11}
        &\ Y(\si(T))
        \triangleq\sup_{0\leq t\leq \si(T)}\int_{\mathbb{R}^n}|\nabla \om_2|^2\,dx
        +\int_{0}^{\si(T)}\int_{\mathbb{R}^n}|\mathcal{P}\dot{\om_2}|^2\,dxdt\nonumber\\
        \leq&\ C\mathcal{E}_0
        +C\left|\int_0^{\si(T)}\int_{\mathbb{R}^{n}}\na \om_2:(\na\om_2\na u)\,dxdt\right|\nonumber\\
        \leq&\ C\mathcal{E}_0
        +C\left|\int_0^{\si(T)}\int_{\mathbb{R}^{n}}\Big(\mathcal{W}+(-\Delta)^{-1}\nabla{\mathcal{P}}
        \diver(\mathbb{F}\mathbb{F}^{\top}-\mathbb I)\Big)
        :\Big[\Big(\mathcal{W}+(-\Delta)^{-1}\nabla{\mathcal{P}}
        \diver(\mathbb{F}\mathbb{F}^{\top}-\mathbb I)\Big)\na u\Big]\,dxdt\right|\nonumber\\
        \leq&\ C\mathcal{E}_0
        +C\int_0^{\si(T)}\int_{\mathbb{R}^{n}}|u||\mathcal W||\na\mathcal W|\,dxdt
        +C\int_0^{\si(T)}\int_{\mathbb{R}^{n}}|\nabla u||(-\Delta)^{-1}\nabla{\mathcal{P}}\diver(\mathbb{F}\mathbb{F}^{\top}-\mathbb I)|^2\,dxdt\nonumber\\
        &\ +C\int_0^{\si(T)}\int_{\mathbb{R}^{n}}|\nabla u||\mathcal W||(-\Delta)^{-1}\nabla{\mathcal{P}}\diver(\mathbb{F}\mathbb{F}^{\top}-\mathbb I)|\,dxdt\nonumber\\
        \triangleq&\ C\mathcal{E}_0+\sum_{i=1}^3 Y_i.
    \end{align}
Indeed, it follows from H\"{o}lder's inequality, \eqref{prop4.1.1}, and \eqref{4.1.1} that for any $k\geq2$,
    \begin{align}
        \label{4.3.13.1}
        \sup\limits_{0\leq t\leq T}
        \|(\mathbb{F}\mathbb{F}^{\top}-\mathbb I)\|_{L^k(\mathbb{R}^n)}^k
        \leq C\sup\limits_{0\leq t\leq T}
        \Big[\Big(\|\mathbb{F}-\mathbb{I}\|_{L^\infty(\mathbb{R}^n)}^{2k-2}
        +\|\mathbb{F}-\mathbb{I}\|_{L^\infty(\mathbb{R}^n)}^{k-2}\Big)
        \|\mathbb{F}-\mathbb{I}\|_{L^2(\mathbb{R}^n)}^2\Big]
        \leq C\mathcal{E}_0.
    \end{align}
For $Y_1$, utilizing H\"{o}lder's and Young's inequalities, the $L^p$-boundedness of the Riesz transform, \eqref{4.1.1}, \eqref{4.2.1}, and \eqref{4.3.13.1}, one deduces for $q$ defined in \eqref{q-0404},
    \begin{align}
        \label{4.4.11-1}
        Y_1
        \leq&\ C\int_0^{\si(T)}\|u\|_{L^q(\mathbb{R}^n)}\|\mathcal W\|_{L^{\frac{2q}{q-2}}(\mathbb{R}^n)}\|\na \mathcal W\|_{L^2(\mathbb{R}^n)}\,dt\nonumber\\
        \leq&\ C\mathcal{E}_0^{\frac{\al}{q}}
        \int_0^{\si(T)}(\|\na\om_2\|_{L^2(\mathbb{R}^n)}
        +\|(\mathbb{F}\mathbb{F}^{\top}-\mathbb I)\|_{L^2(\mathbb{R}^n)})^{\frac{q-n}{q}}
        \|\mathcal P\dot \om_2\|_{L^2(\mathbb{R}^n)}^{\frac{q+n}{q}}\,dt\nonumber\\
        \leq&\ C\mathcal{E}_0
        +C\mathcal{E}_0^{\frac{\al}{q}}Y(\si(T)),
    \end{align}
where we have used the following estimate
    \begin{equation*}
        \|\mathcal W\|_{L^{\frac{2q}{q-2}}(\mathbb{R}^n)}
        \leq C\|\mathcal W\|_{L^2(\mathbb{R}^n)}^{\frac{q-n}{q}}
        \|\na\mathcal W\|_{L^2(\mathbb{R}^n)}^{\frac{n}{q}}
        \leq C(\|\na\om_2\|_{L^2(\mathbb{R}^n)}
        +\|(\mathbb{F}\mathbb{F}^{\top}-\mathbb I)\|_{L^2(\mathbb{R}^n)})^{\frac{q-n}{q}}
        \|\mathcal P\dot \om_2\|_{L^2(\mathbb{R}^n)}^{\frac{n}{q}}
    \end{equation*}
by Lemma \ref{le2.3} and \eqref{4.4.10-0404}-\eqref{4.4.10}.
For $Y_2$, by H\"{o}lder's inequality, \eqref{4.1.1}, and \eqref{4.3.13.1}, we get
    \begin{equation}
        \label{4.4.11-2}
        Y_2
        \leq C\Big(\int_0^{\si(T)}\|\na u\|_{L^2(\mathbb{R}^n)}^2\,dt\Big)^{\frac{1}{2}}
        \Big(\int_0^{\si(T)}\|\mathbb{F}\mathbb{F}^{\top}-\mathbb I\|_{L^4(\mathbb{R}^n)}^4\,dt\Big)^{\frac{1}{2}}
        \leq C\mathcal{E}_0.
    \end{equation}
For $Y_3$, employing the arguments similar to the ones used in \eqref{4.4.11-1} yields for $q$ defined in \eqref{q-0404},
    \begin{align}
        \label{4.4.11-3}
        Y_3
        \leq&\ C\int_0^{\si(T)}\|\na u\|_{L^2(\mathbb{R}^n)}\|\mathcal W\|_{L^{\frac{2q}{q-2}}(\mathbb{R}^n)}\|\mathbb{F}\mathbb{F}^{\top}-\mathbb I\|_{L^q(\mathbb{R}^n)}\,dt\nonumber\\
        \leq&\ C\int_0^{\si(T)}\|\na u\|_{L^2(\mathbb{R}^n)}^2\,dt
        +\sup\limits_{0\leq t\leq \si(T)}
        \|(\mathbb{F}\mathbb{F}^{\top}-\mathbb I)\|_{L^q(\mathbb{R}^n)}^2\nonumber\\
        &\ \times\int_0^{\si(T)}(\|\na\om_2\|_{L^2(\mathbb{R}^n)}
        +\|(\mathbb{F}\mathbb{F}^{\top}-\mathbb I)\|_{L^2(\mathbb{R}^n)})^{\frac{2(q-n)}{q}}
        \|\mathcal P\dot \om_2\|_{L^2(\mathbb{R}^n)}^{\frac{2n}{q}}\,dt\nonumber\\
        \leq&\ C\mathcal{E}_0
        +C\mathcal{E}_0^{\frac{2}{q}}Y(\si(T)).
    \end{align}
Plugging \eqref{4.4.11-1}-\eqref{4.4.11-3} into \eqref{4.4.11} leads to
    \begin{equation*}
        Y(\si(T))\leq C\mathcal{E}_0+C_3\mathcal{E}_0^{\frac{\al}{q}}Y(\si(T)),
    \end{equation*}
which implies
    \begin{equation}
        \label{4.4.12}
        \sup\limits_{0\leq t\leq T}\int_{\mathbb{R}^n}|\nabla w_2|^2dx
        +\int_{0}^T\int_{\mathbb{R}^n}|\mathcal P\dot w_2|^2dxdt
        \leq C\mathcal{E}_0
    \end{equation}
as long as
    \begin{equation*}
        \mathcal{E}_0
        \leq\de_{2}
        \triangleq \min\Big\{\de_{2,1},(2C_3)^{-\frac{q}{\al}}\Big\}.
    \end{equation*}
Choosing $\om_{10}=u_0$, we have $\om_1+\om_2=u$ and hence obtain \eqref{4.4.1} from \eqref{4.4.9} and \eqref{4.4.12}.
Furthermore, it follows from \eqref{ck0}, \eqref{ck1}, and \eqref{4.4.1} that
    \begin{align}
        \label{4.4.13}
        \sup\limits_{0\leq t\leq \sigma(T)}t^{1-\eta}
        \int_{\mathbb{R}^{n}}|\mathcal{G}_i|^2\,dx
        \leq&\ C\sup\limits_{0\leq t\leq \sigma(T)}t^{1-\eta}\int_{\mathbb{R}^{n}}\Big(|\na u|^2+|\mathbb F\mathbb F^\top-\mathbb I|^2\Big)\,dx\nonumber\\
        \leq&\ C\sup\limits_{0\leq t\leq \sigma(T)}\Big[t^{1-\eta}\int_{\mathbb{R}^{n}}|\nabla u|^2\,dx
        +\Big(\|\mathbb F-\mathbb I\|_{L^\infty(\mathbb{R}^n)}^2+1\Big)
        \int_{\mathbb{R}^{n}}|\mathbb F-\mathbb I|^2\,dx\Big]\nonumber\\
        \leq&\ C,
    \end{align}
and
    \begin{equation}
        \label{4.4.14}
        \int_0^{\sigma(T)}\int_{\mathbb{R}^{n}}t^{1-\eta}|\nabla\mathcal{G}_i|^2\,dxdt
        \leq C\int_0^{\sigma(T)}\int_{\mathbb{R}^{n}}t^{1-\eta}|\mathcal{P}\dot{u}|^2\,dxdt
        \leq C.
    \end{equation}

Next, we prove \eqref{4.4.1-0406}.
Multiplying \eqref{4.1.15} by $t^{2-\eta}\mathcal P \dot u$, and integrating the resulting equation over $\mathbb{R}^n\times [0, t]$ with $0\leq t\leq T$, one can deduce from the arguments similar to those used for \eqref{4.1.21} that
    \begin{align}
        \label{4.4.14-0312}
        &\ \sup\limits_{0\leq t\leq \si(T)}\int_{\mathbb{R}^{n}}t^{2-\eta}|\mathcal{P}\dot {u}|^2\,dx
        +\int_0^{\si(T)}\int_{\mathbb{R}^{n}}t^{2-\eta}
        |\nabla\mathcal{P}\dot{u}|^2\,dxdt\nonumber\\
        \leq&\ C\mathcal{E}_0
        +C\int_0^{\si(T)}\int_{\mathbb{R}^{n}}t^{2-\eta}|\nabla u|^4\,dxdt\nonumber\\
        \leq&\ C\mathcal{E}_0
        +C\int_0^{\si(T)}\int_{\mathbb{R}^{n}}|\mathbb F\mathbb F^\top-\mathbb I|^4\,dxdt
        +C\int_0^{\si(T)}\int_{\mathbb{R}^{n}}t^{2-\eta}|\mathcal{G}_1|^4\,dxdt\nonumber\\
        \leq&\ C\mathcal{E}_0
        +C\int_0^{\si(T)}\int_{\mathbb{R}^{n}}t^{2-\eta}|\mathcal{G}_1|^4\,dxdt,
    \end{align}
where we have used the definition of $\mathcal{G}_1$ (see \eqref{ck0}), the $L^p$-boundedness of the Riesz transform, and \eqref{4.3.13.1}.
When $n=2$, it follows from Lemma \ref{le2.3}, \eqref{4.4.13}, and \eqref{4.4.14} that for $\eta\in(0,1]$,
    \begin{align}
        \label{4.4.15-0312}
        \int_0^{\si(T)}\int_{\mathbb{R}^{2}}t^{2-\eta}|\mathcal{G}_1|^4\,dxdt
        \leq&\ C\int_0^{\si(T)}t^{\eta}
        \Big(t^{1-\eta}\int_{\mathbb{R}^{2}}|\mathcal{G}_1|^2\,dx\Big)
        \Big(t^{1-\eta}\int_{\mathbb{R}^{2}}|\nabla\mathcal{G}_1|^2\,dx\Big)\,dt\nonumber\\
        \leq&\ C\Big(\sup_{0\leq t\leq \si(T)}
        t^{1-\eta}\int_{\mathbb{R}^{2}}|\mathcal{G}_1|^2\,dx\Big)
        \int_0^{\si(T)}\int_{\mathbb{R}^{2}}t^{1-\eta}|\nabla\mathcal{G}_1|^2\,dxdt\nonumber\\
        \leq&\ C.
    \end{align}
When $n=3$, repeating the arguments analogous to the ones used in \eqref{4.4.15-0312}, we have for $\eta\in(\frac{8}{15},1]$,
    \begin{align}
        \label{4.4.16-0312}
        &\ \int_0^{\si(T)}\int_{\mathbb{R}^{3}}t^{2-\eta}|\mathcal{G}_1|^4\,dxdt\nonumber\\
        \leq&\ C\int_0^{\si(T)}t^{\frac{2\eta-1}{2}}
        \Big(t^{1-\eta}\int_{\mathbb{R}^{3}}|\mathcal{G}_1|^2\,dx\Big)^{\frac{1}{2}}
        \Big(t^{\frac{4-3\eta}{3}}\int_{\mathbb{R}^{3}}|\nabla\mathcal{G}_1|^2\,dx\Big)^{\frac{3}{2}}\,dt\nonumber\\
        \leq&\ C\Big(\sup_{0\leq t\leq \si(T)}
        t^{1-\eta}\int_{\mathbb{R}^{3}}|\mathcal{G}_1|^2\,dx\Big)^\frac{1}{2}
        \Big(\sup_{0\leq t\leq \si(T)}t^{2-\eta}
        \int_{\mathbb{R}^{3}}|\nabla\mathcal{G}_1|^2\,dx\Big)^\frac{1}{2}
        \int_0^{\si(T)}\int_{\mathbb{R}^{3}}t^{1-\eta}
        |\nabla\mathcal{G}_1|^2\,dxdt\nonumber\\
        \leq&\ C+\frac{1}{2}\sup_{0\leq t\leq
            \si(T)}t^{2-\eta}\int_{\mathbb{R}^{3}}|\mathcal{P}\dot {u}|^2\,dx.
    \end{align}
The estimate of \eqref{4.4.1-0406} follows from \eqref{4.4.14-0312}-\eqref{4.4.16-0312}.
The proof of Lemma \ref{le4.4} is completed.
\end{proof}

The required energy estimate of $A_1(T)$ is derived in the following lemma.
\begin{Lemma}\label{le4.3}
Under the conditions of Proposition \ref{prop4.1}, there exists a positive constant $\de_3\leq\de_2$ such that if $\mathcal{E}_0\leq\de_3$, then
    \begin{equation}
        \label{4.3.1}
        A_1(T)\leq \mathcal{E}_0^{\frac{1}{4}}.
    \end{equation}
\end{Lemma}

\begin{proof}
First, we prove that
    \begin{equation}
        \label{4.3.2}
        \int_0^T\int_{\mathbb{R}^{n}}\sigma^n|\nabla u|^4\,dxdt
        \leq C\mathcal{E}_0^{\frac{1}{2}}.
    \end{equation}
It follows from $\eqref{1.1}_2$ and \eqref{324} that
    \begin{equation}
        \label{4.3.3}
        (\mathbb F-\mathbb I)_t+u\cdot\na(\mathbb F-\mathbb I)
        +\frac{1}{\mu}(\mathbb F-\mathbb I)
        =\frac{1}{\mu}\mathcal{G}\mathbb F+\frac{1}{\mu}\mathcal{R}\mathbb F
        -\frac{1}{\mu}(\mathbb F-\mathbb I)^2.
    \end{equation}
Multiplying \eqref{4.3.3} by $4\sigma^n|\mathbb F-\mathbb I|^2(\mathbb F-\mathbb I)$, and integrating the resulting equation over $\mathbb{R}^n\times[0, T]$, by Young's inequality, the $L^p$-boundedness of the Riesz transform, \eqref{prop4.1.1}, and \eqref{324-0311}, we obtain
    \begin{align}
        \label{4.3.4}
        &\ \int_{\mathbb{R}^{n}}\sigma^n(T)|\mathbb F(x,T)-\mathbb I|^4\,dx
        -n\int_0^T\int_{\mathbb{R}^{n}}\sigma^{n-1}\si'|\mathbb F-\mathbb I|^4\,dxdt
        +\frac{4}{\mu}\int_0^T\int_{\mathbb{R}^{n}}\sigma^n|\mathbb F-\mathbb I|^4\,dxdt\nonumber\\
        \leq&\ C\int_0^T\int_{\mathbb{R}^{n}}\sigma^n
        \Big[(|\mathcal{G}|+|\mathcal{R}|+|\mathbb F-\mathbb I|)|\mathbb F-\mathbb I|^4+(|\mathcal{G}|+|\mathcal{R}|)|\mathbb F-\mathbb I|^3\Big]\,dxdt\nonumber\\
        \leq&\ \frac{1}{2\mu}\int_0^T\int_{\mathbb{R}^{n}}\sigma^n|\mathbb F-\mathbb I|^4\,dxdt
        +C\int_0^T\int_{\mathbb{R}^{n}}\sigma^n\Big(|\mathcal{G}|^4+|\mathcal{R}|^4
        +|\mathbb F-\mathbb I|^5+|\mathbb F-\mathbb I|^{\frac{16}{3}}\Big)\,dxdt\nonumber\\
        \leq&\ \frac{1}{2\mu}\int_0^T\int_{\mathbb{R}^{n}}\sigma^n|\mathbb F-\mathbb I|^4\,dxdt+C\int_0^T\int_{\mathbb{R}^{n}}\sigma^n|\mathcal{G}_2|^4\,dxdt\nonumber\\
        &\ +C\sup\limits_{0\leq t\leq T}\Big(
        \|\mathbb{F}-\mathbb{I}\|_{L^\infty(\mathbb{R}^n)}^4
        +\|\mathbb{F}-\mathbb{I}\|_{L^\infty(\mathbb{R}^n)}
        +\|\mathbb{F}-\mathbb{I}\|_{L^\infty(\mathbb{R}^n)}^{\frac{4}{3}}\Big)
        \int_0^T\int_{\mathbb{R}^{n}}\sigma^n|\mathbb F-\mathbb I|^4\,dxdt\nonumber\\
        \leq&\ \Big(\frac{1}{2\mu}+C_4\mathcal{E}_0^{\frac{\be}{2}}\Big)
        \int_0^T\int_{\mathbb{R}^{n}}\sigma^n|\mathbb F-\mathbb I|^4\,dxdt
        +C\int_0^T\int_{\mathbb{R}^{n}}\sigma^n|\mathcal{G}_2|^4\,dxdt.
    \end{align}
Note that
    \begin{align*}
        &\ n\int_0^T\int_{\mathbb{R}^{n}}\sigma^{n-1}\si'|\mathbb F-\mathbb I|^4\,dxdt
        =n\int_0^{\si(T)}\int_{\mathbb{R}^{n}}\sigma^{n-1}|\mathbb F-\mathbb I|^4\,dxdt\\
        \leq&\ C\sup\limits_{0\leq t\leq \si(T)}\Big(\|\mathbb{F}-\mathbb{I}\|_{L^\infty(\mathbb{R}^n)}^2
        \|\mathbb{F}-\mathbb{I}\|_{L^2(\mathbb{R}^n)}^2\Big)
        \leq C\mathcal{E}_0.
    \end{align*}
Plugging the above estimate into \eqref{4.3.4} gives
    \begin{equation}
        \label{4.3.5}
        \int_0^T\int_{\mathbb{R}^{n}}\sigma^n|\mathbb F-\mathbb I|^4\,dxdt
        \leq
        C\mathcal{E}_0
        +C\int_0^T\int_{\mathbb{R}^{n}}\sigma^n|\mathcal{G}_2|^4\,dxdt,
    \end{equation}
provided that
    \begin{equation*}
        \mathcal{E}_0
        \leq\de_{3,1}
        \triangleq \min\Big\{\de_2,(2\mu C_4)^{-\frac{2}{\be}}\Big\}.
    \end{equation*}
By the definition of ${\mathcal G}_1$ (see \eqref{ck0}) and \eqref{4.3.5}, one has
    \begin{align}
        \label{4.3.6}
        \int_0^T\int_{\mathbb{R}^{n}}\sigma^n|\nabla u|^4\,dxdt
        \leq&\ C\Big(\int_0^T\int_{\mathbb{R}^{n}}\sigma^n|(-\Delta)^{-1}\nabla \mathcal{P}\operatorname{div}(\mathbb F\mathbb F^\top-\mathbb I)|^4\,dxdt+\int_0^T\int_{\mathbb{R}^{n}}\sigma^n|\mathcal{G}_1|^4\,dxdt\Big)\nonumber\\
        \leq&\ C\Big(\int_0^T\int_{\mathbb{R}^{n}}\sigma^n|\mathbb F-\mathbb I|^4\,dxdt+\int_0^T\int_{\mathbb{R}^{n}}\sigma^n|\mathcal{G}_1|^4\,dxdt\Big)\nonumber\\
        \leq&\ C\mathcal{E}_0
        +C\int_0^T\int_{\mathbb{R}^{n}}
        \sigma^n\Big(|\mathcal{G}_1|^4+|\mathcal{G}_2|^4\Big)\,dxdt.
    \end{align}
It follows from Lemma \ref{le2.3} that for $i=1, 2$,
    \begin{align}
        \label{4.3.7}
        &\ \int_0^T\int_{\mathbb{R}^{n}}\sigma^n|\mathcal{G}_i|^4\,dxdt
        \leq C\int_0^T\sigma^n\Big(\int_{\mathbb{R}^{n}}|\mathcal{G}_i|^2\,dx\Big)^\frac{4-n}{2}
        \Big(\int_{\mathbb{R}^{n}}|\nabla\mathcal{G}_i|^2\,dx\Big)^\frac{n}{2}\,dt\nonumber\\
        \leq&\
        \begin{cases}
			\ \displaystyle C\Big(\sup_{0\leq t\leq T}
        \sigma\int_{\mathbb{R}^{2}}|\mathcal{G}_i|^2\,dx\Big)
        \int_0^T\int_{\mathbb{R}^{2}}\sigma|\nabla\mathcal{G}_i|^2\,dxdt,& n=2,\smallskip\nonumber\\
			\ \displaystyle C\Big(\sup_{0\leq t\leq T}
        \sigma\int_{\mathbb{R}^{3}}|\mathcal{G}_i|^2\,dx\Big)^\frac{1}{2}
        \Big(\sup_{0\leq t\leq T}\sigma^{3}\int_{\mathbb{R}^{3}}|\nabla\mathcal{G}_i|^2\,dx\Big)^\frac{1}{2}
        \int_0^T\int_{\mathbb{R}^{3}}\sigma|\nabla\mathcal{G}_i|^2\,dxdt, & n=3
		\end{cases}\nonumber\\
        \leq&\ C\mathcal{E}_0^{\frac{1}{2}},
    \end{align}
where in the last inequality one has used the following estimates
    \begin{align}
        \label{4.3.8}
        \sup_{0\leq t\leq T}\sigma\int_{\mathbb{R}^{n}}|\mathcal{G}_i|^2\,dx
        \leq&\ C\sup_{0\leq t\leq T}\int_{\mathbb{R}^{n}}\Big(|\na u|^2+|\mathbb F\mathbb F^\top-\mathbb I|^2\Big)\,dx\nonumber\\
        \leq&\ C\sup_{0\leq t\leq T}\Big[\int_{\mathbb{R}^{n}}|\nabla u|^2\,dx
        +\Big(\|\mathbb F-\mathbb I\|_{L^\infty(\mathbb{R}^n)}^2+1\Big)
        \int_{\mathbb{R}^{n}}|\mathbb F-\mathbb I|^2\,dx\Big]\nonumber\\
        \leq&\ C[A_1(T)+(B(T)+1)\mathcal{E}_0]
        \leq C\mathcal{E}_0^\frac{1}{4},
    \end{align}
    \begin{equation*}
        \int_0^T\int_{\mathbb{R}^{n}}\sigma|\nabla\mathcal{G}_i|^2\,dxdt
        \leq C\int_0^T\int_{\mathbb{R}^{n}}\sigma|\mathcal{P}\dot{u}|^2\,dxdt
        \leq CA_1(T)
        \leq C\mathcal{E}_0^\frac{1}{4},
    \end{equation*}
and
    \begin{equation}
        \label{4.3.10}
        \sup_{0\leq t\leq T}\sigma^n\int_{\mathbb{R}^{n}}|\nabla \mathcal{G}_i|^2\,dx
        \leq C\sup_{0\leq t\leq T}\sigma^n\int_{\mathbb{R}^{n}}|\mathcal{P}\dot{u}|^2\,dx
        \leq CA_1(T)
        \leq C\mathcal{E}_0^\frac{1}{4}
    \end{equation}
due to \eqref{ck0}-\eqref{ck1}.
Putting \eqref{4.3.7} into \eqref{4.3.6} yields \eqref{4.3.2}.

Next, we show that
    \begin{equation}
        \label{4.3.11.1-0311}
        \left|\int_0^T\int_{\mathbb{R}^{n}}\sigma \na u:(\na u\na u)\,dxdt\right|
        \leq C\mathcal{E}_0^{\frac{5}{16}}.
    \end{equation}
When $n=2$, it follows Cauchy's inequality and \eqref{4.3.2} that
    \begin{equation*}
        \left|\int_0^T\int_{\mathbb{R}^{2}}\sigma \na u:(\na u\na u)\,dxdt\right|
        \leq C\Big(\int_0^T\int_{\mathbb{R}^{2}}|\nabla u|^2\,dxdt+\int_0^T\int_{\mathbb{R}^{2}}\sigma^2|\nabla u|^4\,dxdt\Big)
        \leq C\mathcal{E}_0^{\frac{1}{2}}.
    \end{equation*}
When $n=3$ and $T>1$, utilizing Cauchy's inequality and \eqref{4.3.2}, one gets
    \begin{equation*}
        \left|\int_{1}^T\int_{\mathbb{R}^{3}}\sigma \na u:(\na u\na u)\,dxdt\right|
        \leq C\Big(\int_1^T\int_{\mathbb{R}^{3}}|\nabla u|^2\,dxdt+\int_1^T\int_{\mathbb{R}^{3}}|\nabla u|^4\,dxdt\Big)
        \leq C\mathcal{E}_0^{\frac{1}{2}}.
    \end{equation*}
Therefore, it is sufficient to prove
    \begin{align}
        \label{4.3.11}
        \left|\int_0^{\sigma(T)}\int_{\mathbb{R}^{3}}t \na u:(\na u\na u)\,dxdt\right|
        \leq C\mathcal{E}_0^{\frac{5}{16}}.
    \end{align}
By the definition of $\mathcal{G}_1$ (see \eqref{ck0}), integration by parts, and Cauchy's inequality, we have
    \begin{align}
        \label{4.3.12}
        &\ \left|\int_0^{\sigma(T)}\int_{\mathbb{R}^{3}}t \na u:(\na u\na u)\,dxdt\right|\nonumber\\
        =&\ \frac{1}{\mu^2}\left|\int_0^{\sigma(T)}\int_{\mathbb{R}^{3}}t \na u:\big[\big(\mathcal{G}_1+(-\De)^{-1}\na{\mathcal{P}}
        \diver(\mathbb{F}\mathbb{F}^{\top}-\mathbb{I})\big)
        \big(\mathcal{G}_1+(-\De)^{-1}\na{\mathcal{P}}
        \diver(\mathbb{F}\mathbb{F}^{\top}-\mathbb{I})\big)\big]\,dxdt\right|\nonumber\\
        \leq&\
        \frac{1}{\mu^2}\left|\int_0^{\sigma(T)}\int_{\mathbb{R}^{3}}t \na u:(\mathcal{G}_1\mathcal{G}_1)\,dxdt\right|
        +C\int_0^{\sigma(T)}\int_{\mathbb{R}^{3}}t|\nabla u|^2
        |(-\Delta)^{-1}\nabla{\mathcal{P}}\diver(\mathbb{F}\mathbb{F}^{\top}-\mathbb I)|\,dxdt\nonumber\\
        &\ +C\int_0^{\sigma(T)}\int_{\mathbb{R}^{3}}t|\nabla u|
        |(-\Delta)^{-1}\nabla{\mathcal{P}}\diver(\mathbb{F}\mathbb{F}^{\top}-\mathbb I)|^2\,dxdt\nonumber\\
        \triangleq&\ \sum_{i=1}^3 Z_i.
    \end{align}
For $Z_1$, it follows from H\"{o}lder's and Young's inequalities, \eqref{prop4.1.1}, \eqref{4.1.1}, \eqref{4.4.13}, and \eqref{4.4.14} that
    \begin{align}
        \label{4.3.13}
        Z_1
        \leq&\ C\int_0^{\sigma(T)}t
        \Big(\int_{\mathbb{R}^{3}}|\nabla u|^2\,dx\Big)^\frac12
        \Big(\int_{\mathbb{R}^{3}}|\mathcal{G}_1|^3\,dx\Big)^\frac13
        \Big(\int_{\mathbb{R}^{3}}|\mathcal{G}_1|^2\,dx\Big)^\frac16\,dt\nonumber\\
        \leq&\ C\int_0^{\sigma(T)}t^{\eta-\frac14}
        \Big(\int_{\mathbb{R}^{3}}|\nabla u|^2\,dx\Big)^\frac14
        \Big(t\int_{\mathbb{R}^{3}}|\nabla u|^2\,dx\Big)^\frac14
        \Big(t^{1-\eta}\int_{\mathbb{R}^{3}}|\mathcal{G}_1|^2\,dx\Big)^\frac14
        \Big(t^{1-\eta}\int_{\mathbb{R}^{3}}|\na\mathcal{G}_1|^2\,dx\Big)^\frac34\,dt\nonumber\\
        \leq&\ CA_1(T)^{\frac{1}{4}}
        \Big(\int_0^{\sigma(T)}\int_{\mathbb{R}^{3}}|\nabla u|^2\,dxdt\Big)^\frac14
        \Big(\int_0^{\sigma(T)}t^{1-\eta}\int_{\mathbb{R}^{3}}
        |\na\mathcal{G}_1|^2\,dxdt\Big)^\frac34\nonumber\\
        \leq&\ C\mathcal{E}_0^{\frac{5}{16}}.
    \end{align}
For $Z_2$, by H\"{o}lder's inequality, \eqref{4.1.1}, \eqref{4.3.13.1}, and \eqref{4.3.2}, one has
    \begin{align}
        \label{4.3.14}
        Z_2
        \leq&\ C\Big(\int_0^{\sigma(T)}\int_{\mathbb{R}^{3}}|\nabla u|^2\,dxdt\Big)^\frac12
        \Big(\int_0^{\sigma(T)}\int_{\mathbb{R}^{3}}\sigma^4|\nabla u|^4\,dxdt\Big)^\frac14
        \Big(\int_0^{\sigma(T)}\int_{\mathbb{R}^{3}}|\mathbb{F}\mathbb{F}^{\top}-\mathbb I|^4\,dxdt\Big)^\frac14\nonumber\\
        \leq&\ C\mathcal{E}_0^{\frac{3}{4}}
        \Big(\int_0^{\sigma(T)}\int_{\mathbb{R}^{3}}\sigma^3|\nabla u|^4\,dxdt\Big)^\frac14\nonumber\\
        \leq&\ C\mathcal{E}_0^{\frac{7}{8}}.
    \end{align}
For $Z_3$, it follows from the $L^p$-boundedness of the Riesz transform, Cauchy's inequality, \eqref{prop4.1.1}, and \eqref{4.3.13.1} that
    \begin{align}
        \label{4.3.15}
        Z_3
        =&\ \int_0^{\sigma(T)}\int_{\mathbb{R}^{3}}|\nabla u||(-\Delta)^{-1}\nabla{\mathcal{P}}\diver(\mathbb{F}\mathbb{F}^{\top}-\mathbb I)|^2\,dxdt\nonumber\\
        \leq&\ C\Big(\int_0^{\sigma(T)}\int_{\mathbb{R}^{3}}|\nabla u|^2\,dxdt
        +\int_0^{\sigma(T)}\int_{\mathbb{R}^{3}}|\mathbb{F}\mathbb{F}^{\top}-\mathbb I|^4\,dxdt\Big)\nonumber\\
        \leq&\ C\mathcal{E}_0.
    \end{align}
Plugging \eqref{4.3.13}-\eqref{4.3.15} into \eqref{4.3.12} implies \eqref{4.3.11},
and hence we obtain \eqref{4.3.11.1-0311}.
Substituting \eqref{4.3.2} and \eqref{4.3.11.1-0311} into \eqref{4.1.2} yields
    \begin{equation*}
        A_1(T)\leq C_5\mathcal{E}_0^{\frac{5}{16}},
    \end{equation*}
which gives \eqref{4.3.1} provided that
    \begin{equation*}
        \mathcal{E}_0
        \leq\de_{3}
        \triangleq \min\Big\{\de_{3,1},C_5^{-16}\Big\}.
    \end{equation*}
The proof of Lemma \ref{le4.3} is completed.
\end{proof}


To finish the proof of Proposition \ref{prop4.1}, it remains to bound $B(T)$.
\begin{Lemma}\label{le4.5}
Under the conditions of Proposition \ref{prop4.1}, there exists a positive constant $\de\leq\de_3$ such that if $\mathcal{E}_0\leq\de$, then
    \begin{equation}
        \label{4.5.1}
        B(T)\leq \mathcal{E}_0^{\be}.
    \end{equation}
\end{Lemma}

\begin{proof}
Applying the homogeneous paraproduct operator $\dot\Delta_k$ (see \eqref{bony-1}) to the equation \eqref{4.3.3}, we obtain from \eqref{0} and \eqref{bony} that
    \begin{equation}
        \label{4.5.2}
        (\dot\Delta_k\mathbb {E}^{ij})_t
        +\dot S_{k-1}u\cdot\na(\dot\Delta_k\mathbb {E}^{ij})
        +\frac{1}{\mu}\dot\Delta_k\mathbb {E}^{ij}
        =\frac{1}{\mu}\dot\Delta_k(\mathcal{G}\mathbb F)^{ij}
        +\frac{1}{\mu}\dot\Delta_k(\mathcal{R}\mathbb F)^{ij}
        -\frac{1}{\mu}\dot\Delta_k(\mathbb E\mathbb E)^{ij}
        +f_k,
    \end{equation}
where $\mathbb {E}=\mathbb {F}-\mathbb {I}$ and
    \begin{align*}
        f_k
        \triangleq&\diver(\dot S_{k-1}u\dot\Delta_k\mathbb {E}^{ij}
        -\dot\Delta_k\dot T_u\mathbb {E}^{ij})
        -\dot\Delta_k\diver(\dot T_{\mathbb {E}^{ij}}u
        +\dot R(u,\mathbb {E}^{ij}))\nonumber\\
        =&\diver\Big(\sum_{|k'-k|\leq4}(\dot S_{k-1}u-\dot S_{k'-1}u)\dot\Delta_k\dot\Delta_{k'}\mathbb {E}^{ij}
        -[\dot\Delta_k,\dot S_{k'-1}u]\dot\Delta_{k'}\mathbb {E}^{ij}\Big)
        -\dot\Delta_k\diver(\dot T_{\mathbb {E}^{ij}}u+\dot R(u,\mathbb {E}^{ij})).
    \end{align*}
Then, one deduces from Lemmas \ref{ber}-\ref{jh} that
    \begin{align}
        \label{4.5.2-0312}
        \sum_{k\in\mathbb Z}2^{\frac npk}\|f_k\|_{L^p(\mathbb{R}^n)}
        \leq&\ C\bigg(\|\nabla{u}\|_{L^\infty(\mathbb{R}^n)}\|\mathbb {E}\|_{\dot{B}_{p,1}^{\frac np}(\mathbb{R}^n)}
        +\|\dot T_{\mathbb {E}^{ij}}{u}+\dot R({u},\mathbb {E}^{ij})\|_{\dot{B}_{p,1}^{1+\frac np}(\mathbb{R}^n)}\bigg)\nonumber\\
        \leq&\ C\|\nabla{u}\|_{L^\infty(\mathbb{R}^n)\cap\dot{B}_{p,1}^{\frac np}(\mathbb{R}^n)}
        \|\mathbb E\|_{L^\infty(\mathbb{R}^n)\cap\dot{B}_{p,1}^{\frac np}(\mathbb{R}^n)} \nonumber\\
        \leq&\ C\|\nabla{u}\|_{\dot{B}_{p,1}^{\frac np}(\mathbb{R}^n)}\|\mathbb {E}\|_{\dot{B}_{p,1}^{\frac np}(\mathbb{R}^n)}.
    \end{align}
For $p$ satisfying \eqref{th1.2.0}-\eqref{th1.2.0-0413}, multiplying \eqref{4.5.2} by $2^{\frac npk}|\dot\Delta_k\mathbb {E}^{ij}|^{p-2}\dot\Delta_k\mathbb {E}^{ij}$, summing with respect to $i,j ,k$, integrating the resulting equation over $\mathbb R^n$, by \eqref{4.5.2-0312} and \eqref{prop4.1.1}, we obtain
    \begin{align}
        \label{4.5.3}
        \frac{d}{dt}\|\mathbb E\|_{\dot{B}_{p,1}^{\frac{n}{p}}(\mathbb{R}^n)}
        +\frac{1}{\mu}\|\mathbb E\|_{\dot{B}_{p,1}^{\frac{n}{p}}(\mathbb{R}^n)}
        \leq&\ C\bigg(\|\nabla{u}\|_{\dot{B}_{p,1}^{\frac np}(\mathbb{R}^n)}\|\mathbb {E}\|_{\dot{B}_{p,1}^{\frac np}(\mathbb{R}^n)}
        +\|\mathcal{G}\mathbb{F}\|_{\dot{B}_{p,1}^{\frac{n}{p}}(\mathbb{R}^n)}
        +\|\mathcal{R}\mathbb{F}\|_{\dot{B}_{p,1}^{\frac{n}{p}}(\mathbb{R}^n)}+\|\mathbb E\mathbb E\|_{\dot{B}_{p,1}^{\frac{n}{p}}(\mathbb{R}^n)}\bigg)\nonumber\\
        \leq&\ C\bigg(\|\mathcal G_1\|_{\dot{B}_{p,1}^{\frac np}(\mathbb{R}^n)}\|\mathbb {E}\|_{\dot{B}_{p,1}^{\frac np}(\mathbb{R}^n)}
        +\|\mathcal G_2\|_{\dot{B}_{p,1}^{\frac np}(\mathbb{R}^n)}
        +\|\mathbb E\|_{\dot{B}_{p,1}^{\frac{n}{p}}(\mathbb{R}^n)}^2\bigg)\nonumber\\
        \leq&\ C\bigg(\|\mathcal G_1\|_{\dot{B}_{p,1}^{\frac np}(\mathbb{R}^n)}
        +\|\mathcal G_2\|_{\dot{B}_{p,1}^{\frac np}(\mathbb{R}^n)}\bigg)
        +C_6\mathcal{E}_0^{\frac{\be}{2}}\|\mathbb E\|_{\dot{B}_{p,1}^{\frac{n}{p}}(\mathbb{R}^n)},
    \end{align}
where in the second inequality we have used the following estimates
    \begin{equation*}
        \|\nabla{u}\|_{\dot{B}_{p,1}^{\frac np}(\mathbb{R}^n)}
        \leq C\bigg(\|\mathcal G_1\|_{\dot{B}_{p,1}^{\frac np}(\mathbb{R}^n)}
        +\|(-\Delta)^{-1}\na\mathcal P\diver(\mathbb F\mathbb F^\top-\mathbb I)\|_{\dot{B}_{p,1}^{\frac np}(\mathbb{R}^n)}\bigg)
        \leq C\bigg(\|\mathcal G_1\|_{\dot{B}_{p,1}^{\frac np}(\mathbb{R}^n)}
        +\|\mathbb E\|_{\dot{B}_{p,1}^{\frac{n}{p}}(\mathbb{R}^n)}\bigg),
    \end{equation*}
and
    \begin{equation*}
        \|\mathcal{G}\mathbb{F}\|_{\dot{B}_{p,1}^{\frac{n}{p}}(\mathbb{R}^n)}
        +\|\mathcal{R}\mathbb{F}\|_{\dot{B}_{p,1}^{\frac{n}{p}}(\mathbb{R}^n)}
        \leq C\bigg(\|\mathcal{G}\|_{\dot{B}_{p,1}^{\frac{n}{p}}(\mathbb{R}^n)}
        +\|\mathcal{R}\|_{\dot{B}_{p,1}^{\frac{n}{p}}(\mathbb{R}^n)}\bigg)
        \leq C\bigg(\|\mathcal G_2\|_{\dot{B}_{p,1}^{\frac np}(\mathbb{R}^n)}
        +\|\mathbb E\|_{\dot{B}_{p,1}^{\frac{n}{p}}(\mathbb{R}^n)}^2\bigg)
    \end{equation*}
by the definitions of $\mathcal{G}_1$, $\mathcal{G}$ and $\mathcal{R}$ (see \eqref{ck0} and \eqref{324-0311}), Lemma \ref{TB2}, and the boundedness of the Riesz transform in Besov spaces $\dot B_{p,1}^{\frac np}(\mathbb{R}^n)$ $(1< p<\infty)$.
It follows from \eqref{4.5.3} that
    \begin{equation}
        \label{4.5.4}
        \frac{d}{dt}\|\mathbb E\|_{\dot{B}_{p,1}^{{\frac{n}{p}}}(\mathbb{R}^n)}
        +\frac{1}{2\mu}\|\mathbb E\|_{\dot{B}_{p,1}^{{\frac{n}{p}}}(\mathbb{R}^n)}
        \leq C\bigg(\|\mathcal G_1\|_{\dot{B}_{p,1}^{\frac np}(\mathbb{R}^n)}+\|\mathcal G_2\|_{\dot{B}_{p,1}^{\frac np}(\mathbb{R}^n)}\bigg)
    \end{equation}
if
    \begin{equation*}
        \mathcal{E}_0
        \leq\de_{4}
        \triangleq \min\Big\{\de_3,(2\mu C_6)^{-\frac{2}{\be}}\Big\}.
    \end{equation*}
Employing
Lemma \ref{ri}, Lemma \ref{le2.3}, and \eqref{ck0}-\eqref{ck1}, one has for $i=1, 2$,
    \begin{align}
        \label{4.5.5}
        \|\mathcal{G}_i\|_{\dot{B}_{p,1}^{\frac{n}{p}}(\mathbb{R}^n)}
        \leq&\ C\|\mathcal{G}_i\|_{L^p(\mathbb{R}^n)}^{\frac{p-n}{p}}
        \|\nabla\mathcal{G}_i\|_{L^p(\mathbb{R}^n)}^{\frac{n}{p}}\nonumber\\
        \leq&\ C\Big(\|\mathcal{G}_i\|_{L^2(\mathbb{R}^n)}^{\frac{(2-n)p+2n}{2p}}
        \|\na\mathcal{G}_i\|_{L^2(\mathbb{R}^n)}^{\frac{n(p-2)}{2p}}\Big)
        ^{\frac{p-n}{p}}
        \Big(\|\nabla\mathcal{G}_i\|_{L^2(\mathbb{R}^n)}^{\frac{(2-n)p+2n}{2p}}
        \|\nabla^2\mathcal{G}_i\|_{L^2(\mathbb{R}^n)}^{\frac{n(p-2)}{2p}}\Big)
        ^{\frac{n}{p}}
        \nonumber\\
        \leq&\ C(\|\na u\|_{L^2(\mathbb{R}^n)}+\|\mathbb F\mathbb F^\top-\mathbb I\|_{L^2(\mathbb{R}^n)})^{\frac{(2-n)p^2+n^2p-2n^2}{2p^2}}
        \|\mathcal P\dot u\|_{L^2(\mathbb{R}^n)}^{\frac{np^2-2n^2p+4n^2}{2p^2}}
        \|\na\mathcal P \dot{{u}}\|_{L^2(\mathbb{R}^n)}^{\frac{n^2p-2n^2}{2p^2}}.
    \end{align}

Let
    \begin{equation*}
        \tilde{p}\triangleq\frac{4p^2}{n^2p-2n^2}.
    \end{equation*}
For $1<t\leq T$, multiplying \eqref{4.5.4} by $\|\mathbb E\|_{\dot{B}_{p,1}^{\frac{n}{p}}(\mathbb{R}^n)}^{\tilde{p}-1}$, and integrating the resulting inequality over $(1,t]$, we obtain
    \begin{align}
        \label{4.5.6}
        &\ \|\mathbb E\|_{\dot{B}_{p,1}^{\frac{n}{p}}(\mathbb{R}^n)}^{\tilde{p}}(t)
        +\frac{1}{2\mu}\int_1^{t}\|\mathbb E\|_{\dot{B}_{p,1}^{\frac{n}{p}}(\mathbb{R}^n)}^{\tilde{p}}\,d\tau\nonumber\\
        \leq&\ \|\mathbb E\|_{\dot{B}_{p,1}^{\frac{n}{p}}(\mathbb{R}^n)}^{\tilde{p}}(1)
        +\int_1^{t}\Big(\|\mathcal G_1\|_{\dot{B}_{p,1}^{\frac{n}{p}}(\mathbb{R}^n)}
        +\|\mathcal G_2\|_{\dot{B}_{p,1}^{\frac{n}{p}}(\mathbb{R}^n)}\Big)\|\mathbb E\|_{\dot{B}_{p,1}^{\frac{n}{p}}(\mathbb{R}^n)}^{\tilde{p}-1}\,d\tau\nonumber\\
        \leq&\ \|\mathbb E\|_{\dot{B}_{p,1}^{\frac{n}{p}}(\mathbb{R}^n)}^{\tilde{p}}(1)
        +\frac{1}{4\mu}\int_1^{t}\|\mathbb E\|_{\dot{B}_{p,1}^{\frac{n}{p}}(\mathbb{R}^n)}^{\tilde{p}}\,d\tau
        +C\int_1^{t}\Big(\|\mathcal G_1\|_{\dot{B}_{p,1}^{\frac{n}{p}}(\mathbb{R}^n)}^{\tilde{p}}
        +\|\mathcal G_2\|_{\dot{B}_{p,1}^{\frac{n}{p}}(\mathbb{R}^n)}^{\tilde{p}}\Big)\,d\tau.
    \end{align}
It follows from \eqref{4.5.5}, \eqref{4.3.1}, and \eqref{4.1.1} that for $i=1, 2$,
    \begin{align*}
        \int_1^{t}\|\mathcal G_i\|_{\dot{B}_{p,1}^{\frac{n}{p}}(\mathbb{R}^n)}^{\tilde{p}}\,d\tau
        \leq&\ C\int_1^t\Big(\|\na u\|_{L^2(\mathbb{R}^n)}^2+\|\mathbb F\mathbb F^\top-\mathbb I\|_{L^2(\mathbb{R}^n)}^2+\|\mathcal P\dot u\|_{L^2(\mathbb{R}^n)}^2\Big)^{\frac{\tilde{p}-2}{2}}
        \|\na\mathcal P\dot u\|_{L^2(\mathbb{R}^n)}^2\,d\tau\\
        \leq&\ C\sup\limits_{1\leq t\leq T}\Big(\|\na u\|_{L^2(\mathbb{R}^n)}^2+\|\mathbb F-\mathbb I\|_{L^2(\mathbb{R}^n)}^2+\|\mathcal P\dot u\|_{L^2(\mathbb{R}^n)}^2\Big)^{\frac{\tilde{p}-2}{2}}
        \int_1^T\|\na\mathcal P\dot u\|_{L^2(\mathbb{R}^n)}^2\,dt\\
        \leq&\ C\Big(A_1(T)^{\frac{\tilde{p}-2}{2}}
        +\mathcal{E}_0^{\frac{\tilde{p}-2}{2}}\Big)A_1(T)
        \leq C\mathcal{E}_0^{\frac{\tilde{p}}{8}}.
    \end{align*}
Substituting the above inequality into \eqref{4.5.6} yields
    \begin{equation}
        \label{4.5.7}
        B(T)\leq C\Big(B(1)+\mathcal{E}_0^{\frac{1}{4}}\Big),\quad T>1.
    \end{equation}

For $0\leq t\leq \si(T)$, integrating \eqref{4.5.4} over $[0,t]$, one gets
    \begin{equation}
        \label{4.5.8}
        \|\mathbb E\|_{\dot{B}_{p,1}^{\frac{n}{p}}(\mathbb{R}^n)}(t)
        +\frac{1}{2\mu}\int_0^{t}\|\mathbb E\|_{\dot{B}_{p,1}^{\frac{n}{p}}(\mathbb{R}^n)}\,d\tau
        \leq \|\mathbb E_0\|_{\dot{B}_{p,1}^{\frac{n}{p}}(\mathbb{R}^n)}
        +C\int_0^{t}\Big(\|\mathcal G_1\|_{\dot{B}_{p,1}^{\frac{n}{p}}(\mathbb{R}^n)}
        +\|\mathcal G_2\|_{\dot{B}_{p,1}^{\frac{n}{p}}(\mathbb{R}^n)}\Big)\,d\tau.
    \end{equation}
Utilizing H\"{o}lder's inequality, \eqref{4.1.1}, \eqref{4.4.1-0406}, and \eqref{4.3.1}, we deduce that for $i=1, 2$,
    \begin{align}
        \label{4.5.9}
        &\ \int_0^{t}\|\mathcal G_i\|_{\dot{B}_{p,1}^{\frac{n}{p}}(\mathbb{R}^n)}\,d\tau\nonumber\\
        \leq&\ C\int_0^{\si(T)}
        \Big(t\|\na u\|_{L^2(\mathbb{R}^n)}^2+t\|\mathbb F\mathbb F^\top-\mathbb I\|_{L^2(\mathbb{R}^n)}^2\Big)
        ^{\frac{(2-n)p^2+n^2p-2n^2}{4p^2}}
        \Big(t^{2-\eta}\|\mathcal P\dot u\|_{L^2(\mathbb{R}^n)}^2\Big)
        ^{\frac{np^2-2n^2p+4n^2}{4p^2}-\tilde{\be}}
        \nonumber\\
        &\ \times\Big(t\|\mathcal P\dot u\|_{L^2(\mathbb{R}^n)}^2\Big)
        ^{\tilde{\be}}
        \Big(t^{2-\eta}\|\na\mathcal P\dot u\|_{L^2(\mathbb{R}^n)}^2\Big)
        ^{\frac{n^2p-2n^2}{4p^2}}
        t^{-\frac{[2+(1-\eta)(n-4\tilde{\be})]p^2-(1-\eta)n^2p+2(1-\eta)n^2}{4p^2}}\,dt\nonumber\\
        \leq&\ C\sup\limits_{0\leq t\leq \si(T)}\Big[\Big(t\|\na u\|_{L^2(\mathbb{R}^n)}^2+t\|\mathbb F\mathbb F^\top-\mathbb I\|_{L^2(\mathbb{R}^n)}^2\Big)
        ^{\frac{(2-n)p^2+n^2p-2n^2}{4p^2}}\Big(t^{2-\eta}\|\mathcal P\dot u\|_{L^2(\mathbb{R}^n)}^2\Big)
        ^{\frac{np^2-2n^2p+4n^2}{4p^2}-\tilde{\be}}\Big]\nonumber\\
        &\ \times\Big(\int_0^{\si(T)}
        \Big(t\|\mathcal P\dot u\|_{L^2(\mathbb{R}^n)}^2\Big)
        ^{\frac{4\tilde{\be}p^2}{4p^2-n^2p+2n^2}}
        t^{-\frac{[2+(1-\eta)(n-4\tilde{\be})]p^2-(1-\eta)n^2p+2(1-\eta)n^2}{4p^2-n^2p+2n^2}}\,dt\Big)
        ^{\frac{4p^2-n^2p+2n^2}{4p^2}}\nonumber\\
        &\ \times\Big(\int_0^{\si(T)}t^{2-\eta}\|\na\mathcal P\dot u\|_{L^2(\mathbb{R}^n)}^2\,dt\Big)
        ^{\frac{n^2p-2n^2}{4p^2}}\nonumber\\
        \leq&\ C\Big(\int_0^{\si(T)}t\|\mathcal P\dot u\|_{L^2(\mathbb{R}^n)}^2\,dt\Big)^{\tilde{\be}}
        \Big(\int_0^{1}t^{-\frac{[2+(1-\eta)(n-4\tilde{\be})]p^2-(1-\eta)n^2p+2(1-\eta)n^2}{4(1-\tilde{\be})p^2-n^2p+2n^2}}\,dt\Big)^{\frac{4(1-\tilde{\be})p^2-n^2p+2n^2}{4p^2}},
    \end{align}
where $\eta$ and $p$ satisfy \eqref{th1.2.0}-\eqref{th1.2.0-0413}, and
    \begin{equation*}
        \tilde{\be}\triangleq
        \begin{cases}
            \ \frac{1}{4},&n=2,\smallskip\\
            \ \frac{(15\eta-8)}{48\eta},&n=3.
		\end{cases}
    \end{equation*}
By delicate computations, one has
    \begin{equation*}
        -\frac{[2+(1-\eta)(n-4\tilde{\be})]p^2-(1-\eta)n^2p+2(1-\eta)n^2}{4(1-\tilde{\be})p^2-n^2p+2n^2}
        >-1,
    \end{equation*}
which ensures that the last integral on the right-hand side of \eqref{4.5.9} is convergent.
Therefore, putting \eqref{4.5.9} into \eqref{4.5.8} and using \eqref{4.3.1}, we obtain
    \begin{equation*}
        B(\si(T))
        \leq C\Big(B(0)+A_1(T)^{2\tilde{\be}}\Big)
        \leq C\mathcal{E}_0^{2\be},
    \end{equation*}
which together with \eqref{4.5.7} leads to
    \begin{equation*}
        B(T)\leq C_7\mathcal{E}_0^{2\be}.
    \end{equation*}
From the above estimate, we arrive at \eqref{4.5.1} provided that
    \begin{equation*}
        \mathcal{E}_0
        \leq\de
        \triangleq \min\Big\{\de_{4},C_7^{-\frac{1}{\be}}\Big\}.
    \end{equation*}
The proof of Lemma \ref{le4.5} is finished.
\end{proof}

With the a priori estimates in Proposition \ref{prop4.1}, the local existence result in Proposition \ref{prop1.1}, and the blow-up criteria in Theorem \ref{th1.1}, we are in a position to prove Theorem \ref{th1.2}.

\begin{proof}[\textbf{Proof of Theorem \ref{th1.2}}]
According to Proposition \ref{prop1.1}, there exists a $T_*>0$ such that the Cauchy problem \eqref{1.1}-\eqref{1.2} has a unique classical solution $(u,\mathbb{F})$ in $\mathbb{R}^n\times[0,T_*]$.
Now, we extend the local solution to a global one.

Fix $\de$ in Proposition \ref{prop4.1} and define
    \begin{align}
        \label{4.1}
        T^*\triangleq
        &\ \sup\Big\{T>0\ \big|\
        (u,\mathbb{F})\ \text{is the classical solution of \eqref{1.1}-\eqref{1.2}}\ \text{on}\ [0,T]\ \text{satisfying \eqref{prop4.1.1}}\Big\}.
	\end{align}
We first explain that $T^*>0$ makes sense.
It follows from the definition of $\mathcal{E}_0$ (see \eqref{th1.2.2}), \eqref{A1}-\eqref{B}, and \eqref{4.2.7}-\eqref{4.2.7-0403} that
    \begin{equation*}
        A_1(0)=0,\quad
        A_2(0)\leq \mathcal{E}_0^{\al},\quad
        B(0)\leq \mathcal{E}_0^{\be},
    \end{equation*}
provided $\mathcal{E}_0\leq\de$.
Furthermore, there exists a $T_1\in(0,T_*]$ such that \eqref{prop4.1.1} holds for $T=T_1$.
Therefore, it holds that $T^*\geq T_1>0$ as long as $\mathcal{E}_0\leq\de$.

Next, we claim that
    \begin{equation}
        \label{4.2}
        T^*=\infty
    \end{equation}
under the assumption $\mathcal{E}_0\leq\de$.
Otherwise, suppose $T^*<\infty$.
By virtue of Proposition \ref{prop4.1},
we derive
    \begin{equation}
        \label{4.3}
        A_1(T)\leq \mathcal{E}_0^{\frac14},\quad
        A_2(\si(T))\leq \mathcal{E}_0^{\al},\quad
        B(T)\leq \mathcal{E}_0^{\be},\quad
        \forall\, T\in[0,T^*).
    \end{equation}
Moreover, the Sobolev embedding theorem gives for any $T\in[0,T^*)$,
    \begin{equation}
        \label{4.3-0409-1}
        \|\mathbb{F}\|_{L^\infty(0,T;L^\infty(\mathbb{R}^n))}
        \leq C(1+\|\mathbb{F}-\mathbb{I}\|_{L^\infty(0,T;L^\infty(\mathbb{R}^n))})
        \leq C(1+B(T)),
	\end{equation}
and
    \begin{align}
        \label{4.3-0409-2}
        &\ \int_0^T\|u\|_{L^6(\mathbb{R}^3)}^4\,dt
        = \int_0^{\si(T)}\|u\|_{L^6(\mathbb{R}^3)}^4\,dt
        +\int_{\si(T)}^T\|u\|_{L^6(\mathbb{R}^3)}^4\,dt\nonumber\\
        \leq&\  C\int_0^{\si(T)}\|u\|_{L^3(\mathbb{R}^3)}
        \|u\|_{L^9(\mathbb{R}^3)}^3\,dt
        +C\int_{\si(T)}^T\|\na u\|_{L^2(\mathbb{R}^3)}^4\,dt\nonumber\\
        \leq&\  C\sup\limits_{0\leq t\leq \si(T)}\|u\|_{L^3(\mathbb{R}^3)}
        \int_0^{\si(T)}\big\|\na|u|^{\frac{3}{2}}\big\|_{L^2(\mathbb{R}^3)}^2\,dt
        +C\sup\limits_{\si(T)\leq t\leq T}\|\na u\|_{L^2(\mathbb{R}^3)}^2
        \int_{\si(T)}^T\|\na u\|_{L^2(\mathbb{R}^3)}^2\,dt\nonumber\\
        \leq&\  C\Big(A_2(\si(T))^{\frac{4}{3}}+A_1(T)\mathcal{E}_0\Big),
	\end{align}
where one has used the following estimate
    \begin{equation*}
        \|u\|_{L^9(\mathbb{R}^3)}^3
        =\big\||u|^{\frac{3}{2}}\big\|_{L^6(\mathbb{R}^3)}^2
        \leq C\big\|\na|u|^{\frac{3}{2}}\big\|_{L^2(\mathbb{R}^3)}^2.
	\end{equation*}
From the estimates of \eqref{4.3}-\eqref{4.3-0409-2}, we obtain
    \begin{equation*}
        \begin{cases}
			\ \lim\limits_{T\to T^*}
        \|\mathbb{F}\|_{L^\infty(0,T;L^\infty(\mathbb{R}^2))}<\infty,&n=2,\\
            \ \lim\limits_{T\to T^*}
        \left(\|\mathbb{F}\|_{L^\infty(0,T;L^\infty(\mathbb{R}^3))}
        +\|u\|_{L^4(0,T;L^6(\mathbb{R}^3))}\right)<\infty,&n=3,
		\end{cases}
	\end{equation*}
which contradicts the blow-up criteria \eqref{th1.1.2} and \eqref{th1.1.2-2} in Theorem \ref{th1.1}.
Thus, it can be concluded that $T^*$ is not the maximum
existence time of the solution $(u,\mathbb{F})$.
Furthermore, since the energy functionals $A_1$, $A_2$, and $B$ are continuous with regard to time $t$, we deduce from \eqref{4.3} that there exists a small $a>0$ such that
    \begin{equation*}
        A_1(T)\leq \frac{3}{2}\mathcal{E}_0^{\frac14},\quad
        A_2(\si(T))\leq \frac{3}{2}\mathcal{E}_0^{\al},\quad
        B(T)\leq \frac{3}{2}\mathcal{E}_0^{\be},\quad
        \forall\, T\in[0,T^*+a].
    \end{equation*}
This contradicts the definition of $T^*$ (see \eqref{4.1}).
Hence, \eqref{4.2} holds.

In order to complete the proof of Theorem \ref{th1.2}, it remains to prove the large-time behavior \eqref{th1.2.3} of the solution $(u,\mathbb{F})$.

Multiplying \eqref{4.3.3} by $4|\mathbb F-\mathbb I|^2(\mathbb F-\mathbb I)$, and integrating the resulting equation over $\mathbb{R}^n$, by Young's inequality, the $L^p$-boundedness of the Riesz transform, \eqref{prop4.1.1}, \eqref{324-0311}, and \eqref{4.5.1}, we obtain
    \begin{align}
        \label{4.4}
        \frac{d}{dt}\|\mathbb F-\mathbb I\|_{L^4(\mathbb{R}^n)}^4
        \leq&\ C\int_{\mathbb{R}^{n}}
        \Big[(1+|\mathcal{G}|+|\mathcal{R}|+|\mathbb F-\mathbb I|)|\mathbb F-\mathbb I|^4+(|\mathcal{G}|+|\mathcal{R}|)|\mathbb F-\mathbb I|^3\Big]\,dx\nonumber\\
        \leq&\ C\int_{\mathbb{R}^{n}}\Big(|\mathcal{G}|^4+|\mathcal{R}|^4
        +|\mathbb F-\mathbb I|^4+|\mathbb F-\mathbb I|^5+|\mathbb F-\mathbb I|^{\frac{16}{3}}\Big)\,dx\nonumber\\
        \leq&\ C\|\mathcal{G}_2\|_{L^4(\mathbb{R}^n)}^4
        +C\Big(1+\|\mathbb{F}-\mathbb{I}\|_{L^\infty(\mathbb{R}^n)}^4
        +\|\mathbb{F}-\mathbb{I}\|_{L^\infty(\mathbb{R}^n)}
        +\|\mathbb{F}-\mathbb{I}\|_{L^\infty(\mathbb{R}^n)}^{\frac{4}{3}}\Big)
        \|\mathbb F-\mathbb I\|_{L^4(\mathbb{R}^n)}^4\nonumber\\
        \leq&\ C\Big(\|\mathcal{G}_2\|_{L^4(\mathbb{R}^n)}^4+\|\mathbb F-\mathbb I\|_{L^4(\mathbb{R}^n)}^4\Big).
    \end{align}
Moreover, it follows from \eqref{4.3.5} and \eqref{4.3.7} that
    \begin{equation}
        \label{4.5}
        \int_1^\infty\Big(\|\mathcal{G}_2\|_{L^4(\mathbb{R}^n)}^4+\|\mathbb F-\mathbb I\|_{L^4(\mathbb{R}^n)}^4\Big)\,dt
        \leq C.
    \end{equation}
This together with \eqref{4.4} implies
    \begin{equation}
        \label{4.6}
        \int_1^\infty\left|\frac{d}{dt}\|\mathbb F-\mathbb I\|_{L^4(\mathbb{R}^n)}^4\right|\,dt
        \leq C.
    \end{equation}
From \eqref{4.5} and \eqref{4.6}, we conclude that
    \begin{equation*}
        \|\mathbb F(\cdot,t)-\mathbb I\|_{L^4(\mathbb{R}^n)}\to0\quad\text{as}\ t\to\infty,
    \end{equation*}
which, along with \eqref{4.1.1}, \eqref{4.5.1} and the standard $L^p$-interpolation inequality, leads to
    \begin{equation*}
        \lim_{t\to\infty}\|\mathbb F(\cdot,t)-\mathbb I\|_{L^{p_1}(\mathbb{R}^n)}=0
    \end{equation*}
for all $p_1\in(2,\infty)$.

By the arguments similar to those used in \eqref{4.1.4-0319} and \eqref{4.1.6}, we deduce from Cauchy's inequality that
    \begin{align*}
        \frac{d}{dt}\|\na u\|_{L^2(\mathbb{R}^n)}^2
        =&\ 2\int_{\mathbb R^n} \na u_t:\na u\,dx\\
        =&\ 2\int_{\mathbb R^n}\na\mathcal P\dot u:\na u\,dx
        -2\int_{\mathbb R^n}\na(u\cdot\na u):\na u\,dx
        -2\int_{\mathbb R^n}\mathcal Q(u\cdot\na u)\cdot\De u\,dx\\
        =&\ 2\int_{\mathbb R^n}\na\mathcal P\dot u:\na u\,dx
        -2\int_{\mathbb R^n}(\na u\na u):\na u\,dx\\
        \leq&\ C\Big(\|\na\mathcal P\dot u\|_{L^2(\mathbb{R}^n)}^2+\|\na u\|_{L^2(\mathbb{R}^n)}^2+\|\na u\|_{L^4(\mathbb{R}^n)}^4\Big),
    \end{align*}
which combined with \eqref{4.1.1}, \eqref{4.3.1}, and \eqref{4.3.2} gives
    \begin{equation}
        \label{4.7}
        \int_1^\infty\left|\frac{d}{dt}\|\na u\|_{L^2(\mathbb{R}^n)}^2\right|\,dt
        \leq C.
    \end{equation}
It thus follows from \eqref{4.1.1} and \eqref{4.7} that
    \begin{equation}
        \label{4.8}
        \|\na u(\cdot,t)\|_{L^2(\mathbb{R}^n)}\to0\quad\text{as}\ t\to\infty.
    \end{equation}
For $t\geq1$, using the definition of ${\mathcal G}_1$ (see \eqref{ck0}), the $L^p$-boundedness of the Riesz transform, the Sobolev embedding theorem, \eqref{4.1.1}, \eqref{4.3.8}, \eqref{4.3.10}, and \eqref{4.5.1}, we have
    \begin{equation}
        \label{4.9}
        \|\na u\|_{L^k(\mathbb{R}^n)}
        \leq C(\|\mathcal{G}_1\|_{L^k(\mathbb{R}^n)}
        +\|\mathbb{F}\mathbb{F}^\top-\mathbb{I}\|_{L^k(\mathbb{R}^n)}
        \leq C(\|\mathcal{G}_1\|_{H^1(\mathbb{R}^n)}
        +\|\mathbb{F}\mathbb{F}^\top-\mathbb{I}\|_{L^2(\mathbb{R}^n)\cap L^\infty(\mathbb{R}^n)})
        \leq C
    \end{equation}
with
    \begin{equation*}
        \begin{cases}
			\ k\in[2,\infty),&n=2,\\
            \ k\in[2,6],&n=3.
		\end{cases}
	\end{equation*}
It follows \eqref{4.1.1}, \eqref{4.8}, \eqref{4.9}, the standard $L^p$-interpolation inequality, and Lemma \ref{le2.3} that
    \begin{equation*}
        \lim_{t\to\infty}\big(\|u(\cdot,t)\|_{L^{p_2}(\mathbb{R}^n)}+\|\na u(\cdot,t)\|_{L^{p_3}(\mathbb{R}^n)}\big)=0
    \end{equation*}
for all $p_2\in(2,\infty]$ and
    \begin{equation*}
        \begin{cases}
			\ p_3\in[2,\infty),&n=2,\\
            \ p_3\in[2,6),&n=3.
		\end{cases}
	\end{equation*}
The proof of Theorem \ref{th1.2} is finished.
\end{proof}

\section*{Acknowledgements}
Yong Wang would like to thank The Institute of Mathematical Sciences of The Chinese University of Hong Kong for the hospitality during his visit.

Yong Wang was partially supported by National Key R\&D Program of China (No. 2021YFA1002900) and Guangzhou City Basic and Applied Basic Research Fund (No. 2024A04J6336).
Zhouping Xin was partially supported by the RGC of Hong Kong-SAR (Nos. 14300917, 14300819, 14302819, 14301421), the Basic and Applied Basic Research Foundation of Guangdong Province (No. 2020131515310002), and the Key Project of National Nature Science Foundation of China (No. 12131010).

\bigskip

{\bf Data Availability:} Data sharing is not applicable to this article.

\bigskip

{\bf Conflict of Interest:} The authors declare that they have no conflict of interest.

\end{document}